\documentclass[reqno]{amsart}
\usepackage{enumerate}
\usepackage{amssymb,latexsym}
\usepackage{graphicx}
\usepackage{xcolor} 
\usepackage{verbatim}
\usepackage{imakeidx}
\usepackage{hyperref}
\hypersetup{
    colorlinks=true, 
    linktoc=all,     
    linkcolor=blue,  
}

\newtheorem{theorem}{Theorem}[section]

\newtheorem{lemma}[theorem]{Lemma}
\newtheorem{proposition}[theorem]{Proposition}

\makeindex

\begin{document}
\title[Hypersurfaces with Genuine Deformations in Codimension Two]
{Euclidean Hypersurfaces  with Genuine Conformal Deformations in Codimension Two.}
\author{Sergio Chion and Ruy Tojeiro}

\begin{abstract}
In this paper we classify Euclidean hypersurfaces $f\colon M^n \rightarrow \mathbb{R}^{n+1}$ 
with a  principal curvature of multiplicity $n-2$ that admit a genuine conformal 
deformation $\tilde{f}\colon M^n \rightarrow \mathbb{R}^{n+2}$. That 
$\tilde{f}\colon M^n \rightarrow \mathbb{R}^{n+2}$ is a genuine  conformal deformation 
of $f$  means that it is a conformal immersion for which there exists no  open subset 
$U \subset M^n$ such that the restriction $\tilde{f}|_U$ is a composition 
$\tilde f|_U=h\circ f|_U$ of $f|_U$ with a conformal immersion $h\colon V\to \mathbb{R}^{n+2}$ 
of an open subset $V\subset  \mathbb{R}^{n+1}$ containing $f(U)$.
\end{abstract}

\maketitle

\section{Introduction}

Euclidean hypersurfaces $f\colon M^n \rightarrow \mathbb{R}^{n+1}$ that are free
of flat (respectively, conformally flat) points and admit an isometric 
(respectively, conformal) deformation $g\colon M^n \rightarrow \mathbb{R}^{n+1}$
that is not isometrically congruent (respectively, conformally congruent) to $f$ on 
any open subset of $M^n$  are called \emph{Sbrana-Cartan hypersurfaces}
(respectively,  \emph{Cartan hypersurfaces}).  These two types of hypersurfaces 
have been classified in the beginning of the twentieth century: in the isometric case 
 by Sbrana  \cite{vS1909} and Cartan \cite{eC1916} for $n\geq 3$, and in the conformal 
 one by  Cartan  \cite{eC1917} for $n\geq 5$. The most interesting classes of Sbrana-Cartan 
 (respectively, Cartan) hypersurfaces are envelopes of certain two-parameter congruences 
 of affine hyperplanes (respectively, hyperspheres), which may admit either a one-parameter 
 family of isometric (respectively, conformal) deformations, or a single one. 
 Partial results on Cartan hypersurfaces of dimensions four and three were also obtained 
 by Cartan in \cite{eC1918} and \cite{eC1920}, respectively.

The classification of Sbrana-Cartan hypersurfaces  was extended to the case of 
nonflat ambient space forms by Dajczer-Florit-Tojeiro  \cite{mD1998}. Moreover,
among other things, in that paper it was given an
 affirmative answer to the question of  whether Sbrana-Cartan hypersurfaces that allow 
a single  deformation do exist, which  was  not addressed  by Sbrana and  Cartan.  
 
 A nonparametric description of Cartan hypersurfaces of dimension $n\geq 5$ 
  was given in  \cite{mD2000}, where it was shown that any such  hypersurface arises 
  by intersecting the light-cone $\mathbb{V}^{n+2}$
 in Lorentzian space $\mathbb{L}^{n+3}$ with a flat space-like submanifold of  
 codimension two of $\mathbb{L}^{n+3}$.  We  refer to \cite{mD1998} and \cite{mD2000},
 respectively, for  modern accounts of the classifications of Sbrana-Cartan 
 and Cartan hypersurfaces.

When studying isometric or conformal deformations of  a Euclidean submanifold with
codimension greater than one, one has to 
take into account that any submanifold of a deformable submanifold already 
possesses the isometric or conformal deformations induced by the latter. 
Therefore,  it is necessary to restrict the study to
those deformations that are ``genuine'', that is, those that are not induced by 
  deformations of an ``extended'' submanifold. It is also of 
interest to consider  deformations of a submanifold that take place in 
a possibly different codimension.  These ideas have been made precise   in \cite{mD2004} 
in the isometric case, and extended to the conformal realm in \cite{rT2010} as follows.

Let $f\colon M^n \to \mathbb{R}^{n+p}$ be a conformal immersion of an $n$-dimensional 
Riemannian manifold $M^{n}$ into Euclidean space.  A conformal immersion 
$\tilde{f}\colon M^{n} \to \mathbb{R}^{n+q}$ is said to be a 
\emph{genuine conformal deformation} of $f$ if there exists no open subset  $U\subset M^{n}$
such that the restrictions $f|_U$ and $\tilde{f}|_U$ are compositions $f|_U = F \circ j$ and 
$\tilde{f}|_U = \tilde{F}\circ j$ of a conformal embedding 
$j\colon U \to N^{n+\ell}$ into a Riemannian manifold $N^{n+\ell}$, with $\ell>0$, 
and conformal immersions $F\colon N^{n+\ell}\to \mathbb{R}^{n+p}$ and 
$\tilde{F}\colon N^{n+\ell}\to \mathbb{R}^{n+q}$:

\medskip 
 
\begin{picture}(80,84) 
\put(95,31){$U\subset M^n$} 
\put(169,31){$N^{n+\ell}$} 
\put(202,62){$\mathbb{R}^{n+p}$}
\put(202,0){$\mathbb{R}^{n+q}$} 
\put(155,59){${}_f|_U$} 
\put(155,9){${}_{\tilde f|_U}$} 
\put(200,46){${}_F$} 
\put(199,26){${}_{\tilde F}$} 
\put(153,40.5){${}_j$} 
\put(188,42){\vector(1,1){16}} 
\put(188,28){\vector(1,-1){16}} 
\put(137,44){\vector(3,1){60}} 
\put(137,26){\vector(3,-1){60}} 
\put(145,34){\vector(1,0){21}} 
\put(145,36){\oval(7,4)[l]} 
\end{picture} 
\bigskip

     In this work we are interested in the particular case in which $p=1$ and $q=2$.  
     In the isometric realm, from the assumption that 
     $f\colon M^n \to \mathbb{R}^{n+1}$ admits a genuine isometric deformation 
     $\tilde f\colon M^n \to \mathbb{R}^{n+2}$, it follows from  
     Theorem 1 in \cite{mD1996} that $\text{rank}\,f $, that is, the rank of the 
     shape operator of $f$, must be at most $3$ at any point.  
     The case in which $\text{rank}\,f = 2$ was solved in \cite{mD2013}.   
     In the conformal instance, from Theorem 1 of \cite{mD1997} it follows that 
     a Euclidean hypersurface $f\colon M^n \to \mathbb{R}^{n+1}$ must have a 
     principal curvature $\lambda$ with multiplicity greater than or equal to $n-3$ 
     at any point if it admits a genuine conformal deformation 
     $\tilde f\colon M^n \to \mathbb{R}^{n+2}$.  We will study the particular case 
     in which the multiplicity is $n-2$.  For the case $n-3$, it seems  better to start 
     by attempting to solve the analogous problem in the isometric realm,
      which is also still open.

Hypersurfaces $f\colon M^n \to \mathbb{R}^{n+1}$ that have  a principal curvature $\lambda$ 
of multiplicity $n-2$ are envelopes of two-parameter congruences of hyperspheres. 
These are given by a focal function $h\colon L^2 \to \mathbb{R}^{n+1}$, the locus of centers of
the hyperspheres of the congruence, and a 
radius function $r \in C^\infty(L)$, where $L^2 = M^n / \Delta$ 
is the quotient space of leaves of the  eigendistribution distribution $\Delta$ of $\lambda$.
In terms of the  model of Euclidean space $\mathbb{R}^{n+1}$ as a hypersurface of the light-cone 
$\mathbb{V}^{n+2} \subset \mathbb{L}^{n+3}$, the congruence of hyperspheres $(h,r)$ can be 
represented by a surface $s\colon L^2 \to \mathbb{S}^{n+2}_{1,1} \subset \mathbb{L}^{n+3}$ 
in the de Sitter space.   With the aid of the conformal Gauss parametrization, 
the hypersurface $f$ can be recovered back from the surface $s$. 
Our approach is to determine which such surfaces 
give rise to hypersurfaces $f\colon M^n \to \mathbb{R}^{n+1}$ that admit genuine conformal 
deformations $\tilde f\colon M^n \to \mathbb{R}^{n+2}$.

In the proof, we follow similar steps to those of the isometric case.
 We show in Section \ref{ch:thetriple} that the existence of a 
genuine conformal deformation $\tilde f\colon M^n \to \mathbb{R}^{n+2}$ of a hypersurface 
$f\colon M^n \to \mathbb{R}^{n+1}$ with  a principal curvature  of multiplicity $n-2$  
can be encoded by a triple $(D_1,D_2,\psi)$ satisfying several conditions, 
where $D_i \in \Gamma(\text{End}(\Delta^\perp))$, $1\leq i\leq 2$, 
and $\psi$ is a one-form on $M^n$. This requires the preliminary algebraic step 
of determining the structure of the second fundamental form
of the isometric light-cone representative of a  genuine conformal deformation 
$\tilde f\colon M^n \to \mathbb{R}^{n+2}$ of $f$, 
which is carried out in Section \ref{ch:defincodtwo}.

   The next step is to prove that the triple $(D_1,D_2,\psi)$ can be projected down  
   to a triple $(\bar{D}_1,\bar{D}_2,\bar{\psi})$ on the 
 surface  $s\colon L^2 \to \mathbb{S}^{n+2}_{1,1}$, 
 and to express the conditions on  $(D_1,D_2,\psi)$ in terms of simpler ones  
 on $(\bar{D}_1,\bar{D}_2,\bar{\psi})$ (see Section \ref{ch:reduction}). 
This is one of the main differences with respect to the approach used 
in the isometric case in \cite{mD2013}, where this reduction process 
was carried out in terms of the Gauss map and the support function of the hypersurface.

 The last step is  to characterize the surfaces $s\colon L^2 \to \mathbb{S}^{n+2}_{1,1}$ 
 that carry a triple $(\bar{D}_1,\bar{D}_2,\bar{\psi})$ 
 satisfying the aforementioned conditions. 
 This is done in Section~\ref{ch:subset}. 
 For the proof of the classification of Euclidean hypersurfaces 
 $f\colon M^n \to \mathbb{R}^{n+1}$ that admit genuine conformal deformations 
 $\tilde f\colon M^n \to \mathbb{R}^{n+2}$ in Section \ref{ch:classification}, 
 all that was needed was to put together the steps accomplished in the previous sections.

 The main theorem of this article is, as far as we know, the first classification result 
 for a class of submanifolds admitting genuine conformal deformations, apart from the 
 classical one by Cartan of the  hypersurfaces $f\colon M^n \to \mathbb{R}^{n+1}$ 
 that admit genuine conformal deformations $\tilde f\colon M^n \to \mathbb{R}^{n+1}$. 
 In the isometric realm, besides the isometric version of our result in \cite{mD2013}, 
 isometric immersions  $f\colon M^n \to \mathbb{R}^{n+2}$ of rank two  that admit 
 genuine isometric deformations $\tilde f\colon M^n \to \mathbb{R}^{n+2}$ 
 have been classified in \cite{mD20042}, \cite{mD2009} and \cite{lF2017}.

\section[Preliminaries]{Preliminaries}

Two Riemannian metrics $\langle\,,\,\rangle$ and 
$\langle\,,\,\rangle'$ on a manifold $M^n$ 
are \emph{conformal} if there exists 
a positive function $\varphi\in C^\infty(M)$ such that 
$
\langle\,,\,\rangle'=\varphi^2\langle\,,\,\rangle.
$
The function  $\varphi$ is called the 
\emph{conformal factor} of $\langle\,,\,\rangle'$ 
with respect to $\langle\,,\,\rangle$. 
An immersion $f\colon M^n\to\bar M^m$ between Riemannian
manifolds is  \emph{conformal}
if  its induced 
metric $\langle\,,\,\rangle_f$ is conformal 
to the Riemannian metric $\langle\,,\,\rangle$ of $M^n$, and the 
\emph{conformal factor of $f$} is the conformal 
factor of $\langle\,,\,\rangle_f$ 
with respect to $\langle\,,\,\rangle$.
\vspace{1ex}

Let $\mathbb{L}^{m+2}$ be the $(m+2)$--dimensional Minkowski space, that
is, $\mathbb{R}^{m+2}$ endowed with a  Lorentz scalar product of signature
$(-, +,\ldots,+)$, and let
$$
\mathbb{V}^{m+1}= \{p\in\mathbb{L}^{m+2}\colon\langle p,p\rangle=0, \, p\neq 0\}
$$
denote the light cone in $\mathbb{L}^{m+2}$. Then
$$
\mathbb{E}^{m}=\mathbb{E}^{m}_w
=\{p\in\mathbb{V}^{m+1}\colon\langle p,w\rangle=1\}
$$
is a model of $m$--dimensional Euclidean space for any
$w\in\mathbb{V}^{m+1}$. Namely, if $p_0\in \mathbb{E}^{m}$ and 
$C\colon\,\mathbb{R}^{m}\to \mbox{span}\{p_0,w\}^\perp\subset
\mathbb{L}^{m+2}$ is a
linear isometry,  the triple $(p_0,w,C)$ gives rise to an isometric embedding
$\Psi=\Psi_{p_0,w,C}\colon\,\mathbb{R}^{m}\to 
\mathbb{L}^{m+2}$ defined by
$$
\Psi(x)=p_0+Cx-\frac{1}{2}\|x\|^2w
$$
that has $\mathbb{E}^{m}$ as image and whose second fundamental form is
\begin{equation}\label{eq:sffpsi}
\alpha^{\Psi}(Z,W)=-\langle Z,W\rangle w \quad \text{for all} \;\;Z, W\in \mathfrak{X}(\mathbb{R}^{m}).
\end{equation}

Hyperspheres can be nicely described in the model $\mathbb{E}^{m}$ 
of $m$-dimensional Euclidean space: given a
hypersphere $S\subset\mathbb{R}^{m}$  with (constant) mean curvature $H$ with
respect to a unit normal vector field  $N$ along $S$, then
$v=H\Psi+\Psi_*N\in \mathbb{L}^{m+2}$ is a constant space-like vector
of unit length, as follows by differentiating the right-hand-side.
Moreover,  $\langle v,\Psi(q)\rangle=0$ for all
$q\in S$, and hence $\Psi(S)=\mathbb{E}^{m}\cap \{v\}^\perp$. 

In this way, (oriented) hyperspheres of $\mathbb{R}^{n+1}$ are in 
one-to-one correspondence with points of the Lorentzian sphere
$\mathbb{S}_{1,1}^{n+2}= \{p\in\mathbb{L}^{n+3}\colon\langle p,p\rangle=1\}$.
Therefore, given an oriented hypersurface $f\colon\,M^n\to \mathbb{R}^{n+1}$ and 
smooth maps $h\colon M^n\to\mathbb{R}^m$ and $R\in C^{\infty}(M)$, $R> 0$, a
sphere congruence $x\in M^n\mapsto S(h(x),R(x))$, with radius function
$R$ and $h$ as the locus of centers, which is enveloped by $f$,
that is, $$
f(x)\in S(h(x),r(x))\;\;\mbox{and}\;\;f_*T_xM\subset T_{f(x)}S(h(x),r(x))
$$
for all $x\in M^n$, 
can be identified with the map
$S\colon\,M^n\to\mathbb{S}_{1,1}^{n+2}$ given by
\begin{equation}\label{s} 
S(q)=\frac{1}{R(q)}\Psi(f(q))+\Psi_*(f(q))N(q),
\end{equation}
where $N$ is a unit normal vector field along $f$.

If a hypersurface $f\colon M^n\to\mathbb{R}^{n+1}$ envelops a \emph{$k$-parameter} 
congruence of hyperspheres $S\colon M^n\to\mathbb{S}_{1,1}^{n+2}$, 
$1\leq k\leq n-1$, that is,  the map $S$ has rank $k$, 
then $f$ has a  principal curvature $\lambda$ such that 
$\mbox{ker}\, S_*(x)\subset E_\lambda(x)$ for all $x\in M^n$,
with $\mbox{ker}\, S_*(x)= E_\lambda(x)$ for all $x$ in an open 
dense subset of $M^n$, on which $\lambda$ is constant along $E_\lambda$. 
Conversely, any hypersurface $f\colon M^n\to\mathbb{R}^{n+1}$ 
that carries a Dupin principal curvature of multiplicity $n-k$ envelops 
a $k$-parameter congruence of hyperspheres. Recall that a principal curvature 
$\lambda$ is \emph{Dupin} if $\lambda$ is constant along $E_\lambda$, 
which is always the case if the multiplicity of $\lambda$ is at least two. 
Therefore, in this case the congruence of hyperspheres $S$ gives rise to a  map 
$s\colon L^k\to\mathbb{S}_{1,1}^{n+2}$
defined on the quotient space $L^k$ of leaves of $E_\lambda$.

Let us fix $w=(w_0, \ldots, w_{n+2})\in \mathbb{V}^{m+1}\subset \mathbb{L}^{m+2}$ 
with $w_0<0$, so that $\mathbb{E}^{m}$ is contained in the upper half 
$\mathbb{V}^{m+1}_+$ of $\mathbb{V}^{m+1}$. 
Then, any conformal immersion $f\colon M^n\to\mathbb{R}^m$ 
with conformal factor $\varphi\in C^\infty(M)$  gives rise to an 
isometric immersion $\mathcal{I}(f)\colon M^n\to\mathbb{V}_+^{m+1}$ 
into the light-cone of $\mathbb{L}^{m+2}$, given by
$$
\mathcal{I}(f)=\frac{1}{\varphi}\Psi\circ f,
$$
called its
\emph{isometric light-cone representative}. 
Conversely, any isometric immersion  $F\colon M^n\to\mathbb{V}_+^{m+1}\smallsetminus\mathbb{R} w$
gives rise to a conformal immersion $\mathcal{C}(F)\colon M^n\to\mathbb{R}^m$  
with conformal factor $1/\langle F,w\rangle$ given by
$$ 
\Psi\circ \mathcal{C}(F) = \Pi\circ F,
$$
where $\Pi\colon\mathbb{V}_+^{m+1}\smallsetminus\mathbb{R} w \to\mathbb{E}^m$,  
$\mathbb{R} w=\{tw:t<0\}$, denotes the projection 
onto $\mathbb{E}^m$ given by $ \Pi(u)=u/\langle u,w\rangle$.
Moreover, for any conformal immersion $f\colon M^n\to\mathbb{R}^m$ and for any 
isometric immersion $F\colon M^n\to\mathbb{V}_+^{m+1}\smallsetminus\mathbb{R} w$ one has
$$
\mathcal{C}(\mathcal{I}(f))=f\;\;\mbox{and}\;\;\mathcal{I}(\mathcal{C}(F))=F.
$$

Two  immersions $f,g\colon M^n\to\mathbb{R}^m$  are said to be 
\emph{conformally congruent} if $g=\tau\circ f$ for some conformal 
transformation $\tau$ of $\mathbb{R}^m$. The next result is well-known.

\begin{proposition}\label{p:conformalisometric} Two conformal immersions 
$f,g\colon M^n\to\mathbb{R}^m$ are conformally congruent if and only 
if their isometric light-cone representatives $\mathcal{I}(f),
\mathcal{I}(g)\colon M^n\to\mathbb{V}_+^{m+1}\subset \mathbb{L}^{m+2}$ 
are isometrically congruent.
\end{proposition}

\section[Light-cone representatives of conformal deformations]
{Light-cone representatives of conformal deformations}
\label{ch:defincodtwo}

In this section we show how nongenuine conformal deformations 
$\tilde{f}\colon M^n \to \mathbb{R}^{n+2}$ of a hypersurface 
$f\colon M^n \to \mathbb{R}^{n+1}$ can be characterized in terms 
of their isometric light-cone representatives, and study the 
structure of the second fundamental form of the isometric 
light-cone representative of a genuine conformal deformation.

\subsection{Characterizing nongenuine conformal deformations}

 Given conformal immersions $f\colon M^n \to \mathbb{R}^{n+1}$ and 
 $\tilde{f}\colon M^n \to \mathbb{R}^{n+p}$, the following result 
  characterizes,  in terms of their isometric light-cone representatives, 
 when $\tilde f$ is the composition $\tilde f=h\circ f$ of $f$ 
 with a conformal  immersion $h\colon V\to \mathbb{R}^{n+p}$ 
 of an open subset $V \subset \mathbb{R}^{n+1}$ containing $f(M^n)$.

\begin{proposition}\label{p:conformalchar}
Let $f\colon M^n \to \mathbb{R}^{n+1}$ and 
$\tilde f\colon M^n \to \mathbb{R}^{n+p}$ 
be conformal immersions.  Endow $M^n$ with the metric induced by $f$ 
and let $F\colon M^n \to \mathbb{V}^{n+2}\subset\mathbb{L}^{n+3}$ and 
$\tilde F \colon M^n \to \mathbb{V}^{n+p+1}\subset\mathbb{L}^{n+p+2}$ 
be the light-cone  representatives of $f$ and $\tilde f$, respectively.  
Given an open set $U \subset M^n$, there exists a conformal immersion 
$h\colon V\to \mathbb{R}^{n+p}$ of an open subset 
$V \supset f(U)$ of $\mathbb{R}^{n+1}$ such that 
$\tilde f|_U=h \circ f|_U$ if and only if there exists an 
isometric immersion $H\colon W \to \mathbb{V}^{n+p+1}$ of an 
open subset $W \supset F(U)$ of $\mathbb{V}^{n+2}$ such that 
$\tilde F|_U = H \circ F|_U$.
\end{proposition}

\proof 
 Assume first that $H\colon W\to \mathbb{V}^{n+p+1}$ is an  
 isometric immersion  of an open subset  $W \supset F(U)$ of 
 $\mathbb{V}^{n+2}$ such that $\tilde F|_U = H \circ F|_U$. 
 Define $V = \Psi^{-1}(W)$ and consider 
 $H \circ \Psi\colon V \to \mathbb{V}^{n+p+1}$.  
 Then $h = \mathcal{C}(H\circ \Psi)\colon V\to \mathbb{R}^{n+p}$ 
 is a conformal immersion and  
$$\tilde f|_U=\mathcal{C}(\tilde F|_U) = \mathcal{C}(H \circ F|_U) = 
\mathcal{C}(H \circ \Psi) \circ f|_U = h \circ f|_U.$$

Conversely, let $h\colon V\to \mathbb{R}^{n+p}$ be a conformal immersion 
of  an open subset $V \supset f(U)$ of $\mathbb{R}^{n+1}$ such that 
$\tilde f|_U=h \circ f|_U$.  
Let $H\colon \Psi(V) \to \mathbb{V}^{n+p+1}\subset\mathbb{L}^{n+p+2}$ 
be defined by $\mathcal{I}(h) = H \circ \Psi$. 
Then 
$$\mathcal{C}(H \circ F|_U) = \mathcal{C}(H \circ \Psi) \circ f|_U 
= h \circ f|_U = \tilde f|_U,$$
hence $\tilde F|_U = H \circ F|_U$ by Proposition \ref{p:conformalisometric}.  
Now extend $H$ to an isometric immersion 
$H\colon W\subset \mathbb{V}^{n+2}\to \mathbb{V}^{n+p+1}$ 
by setting $ H(t\Psi(x))=tH(\Psi(x))$ for any $x\in V$.\vspace{1ex}\qed

In order to apply Proposition \ref{p:conformalchar}, one must have 
sufficient conditions on a pair of isometric immersions 
$F\colon M^n \to \mathbb{V}^{n+2}\subset\mathbb{L}^{n+3}$ 
and $\tilde F \colon M^n \to \mathbb{V}^{n+p+1}\subset\mathbb{L}^{n+p+2}$ 
which imply the existence of an  isometric immersion 
$H\colon W \to \mathbb{V}^{n+p+1}$ 
of an open subset $W \supset F(M^n)$ of $\mathbb{V}^{n+2}$ 
such that $\tilde F = H \circ F$.
This is the content of  the next lemma in the case of 
interest to us in this work, namely, the case $p=2$.  

\begin{lemma} \label{le:blem} 
Let $F\colon M^n \to \mathbb{V}^{n+2}\subset \mathbb{L}^{n+3}$ 
and $\tilde F\colon M^n\to \mathbb{V}^{n+3} \subset \mathbb{L}^{n+4}$ 
be isometric immersions, and suppose that $F$ is an embedding.   
Assume that there exist $\xi \in \Gamma(N_{\tilde F}M)$ of unit length,
with $\langle \xi,\tilde F \rangle = 0$, 
$\text{rank}\, A_\xi^{\tilde F} = 1$ and 
$^{\tilde F}\nabla_Z^\perp \xi = 0$ for all $Z \in \ker A_\xi^{\tilde F}$,
and a  parallel vector bundle isometry $T\colon N_{F}M \to L = \{\xi\}^\perp$  
with respect to the induced connection on $L$  such that  $TF=\tilde F$ and
$$\alpha_{\tilde F} = T \circ \alpha_F + \big<A_\xi^{\tilde F}, \big>\xi.$$
Then, there exists an isometric immersion $H\colon W \to \mathbb{V}^{n+3}$ 
of an open subset $W\subset \mathbb{V}^{n+2}$ containing $F(M)$ 
such that $\tilde F = H \circ F$.
\end{lemma}
\proof
Let $Y \in \left(\ker A_\xi\right)^\perp$ be an eigenvector of $A_\xi$ having $\beta$ 
as the unique non-zero eigenvalue. Then 
$$W = \left\{\left(\tilde{\nabla}_X \xi\right)_{{\tilde F}_*TM \oplus L} : X\in \mathfrak{X}(M)\right\}$$
is a line subbundle of $\mathbb{R}({\tilde F}_*Y) \oplus L$ spanned 
by the vector field  $-\beta {\tilde F}_*Y + \nabla_Y^\perp \xi$. Its 
 orthogonal complement $\Gamma$  in $\mathbb{R}({\tilde F}_*Y) \oplus L$ 
 is a rank-$3$ subbundle   such that   $\Gamma \cap {\tilde F}_*TM = \{0\}$ and 
 $\tilde{\nabla}_X \delta \in {\tilde F}_*TM \oplus L$ for any  section $\delta$ of $\Gamma$.
Moreover, since the position vector field $\tilde F$ is parallel in the normal connection 
and is everywhere orthogonal to $\xi$ by assumption, it is a section of  $\Gamma$. 
Define a vector-bundle isometry
$\mathcal{T}\colon F_*TM \oplus N_FM \to {\tilde F}_*TM \oplus L$
 by setting 
 $$\mathcal{T}(F_*Y+\eta)={\tilde F}_*Y+T\eta$$ for all $Y\in \mathfrak{X}(M)$ 
 and $\eta \in \Gamma(N_FM)$. The vector subbundle  
 $\Omega = \mathcal{T}^{-1}(\Gamma)$ is transversal to $F_*TM$, 
because $\Gamma \cap {\tilde F}_*TM = \{0\}$.  Also,   the position vector field  $F$ 
is a section of $\Omega$, for $TF=\tilde F$.   Since  $F$ is an embedding,  
the map $G\colon \Omega \to \mathbb{L}^{n+3}$ defined by
$$G(e) = F(x) + e,$$
where $\pi\colon \Omega\to M^n$ is the projection and $x=\pi(e)$, parametrizes a 
tubular neighborhood of $F(M^n)$ if restricted to a neighborhood ${U}$ of the 
$0$-section of $\Omega$.  Endow ${U}$ with the  Lorentzian metric induced by $G$.  
For a vertical vector $Z \in T_{e}\Omega$ we have
$G_*(e)Z = Z$. On the other hand, any nonvertical vector 
$Z \in T_{e}\Omega$ can be written as 
$Z=\zeta_*X$ for some $\zeta\in \Gamma(\Omega)$ with $\zeta(x)=e$ and $X\in T_xM$. 
Writing $\zeta=F_*Y+\eta$, with $Y\in \mathfrak{X}(M)$ and $\eta \in \Gamma(N_FM)$, 
we have
\begin{align*}
G_*(e)Z &= F_*X + \tilde{\nabla}_{X}\left(F_*Y + \eta\right) \\
&= F_*\left(X + \nabla_{X} Y -  A_\eta^F X\right) + \alpha^F(X,Y) + ^F\nabla_{X}^\perp \eta.
\end{align*}  
We claim that the map $\tilde{G}\colon \Omega \to \mathbb{L}^{n+4}$, defined by
$$\tilde{G}(e) = \tilde{F}(x) + \mathcal{T}(e),$$
with $x=\pi(e)$, is an isometric immersion on ${U}$, that is,
$||\tilde{G}_*(e)Z|| = ||{G}_*(e)Z||$
for all $e\in U$ and $Z\in T_eU$.  To prove this, it suffices to show that
\begin{equation}\label{eq:tilG} 
\tilde{G}_*(e)Z =T{G}_*(e)Z
\end{equation} for all $e\in U$ and $Z\in T_eU$, for then the claim follows from 
the fact that $T$ is a vector bundle isometry.

For any vertical $Z\in T_eU$,  \eqref{eq:tilG} follows from 
$\tilde{G}_*(e)Z = TZ=T{G}_*(e)Z$. If $Z=\zeta_*X$ for some 
$\zeta=F_*Y+\eta\in \Gamma(\Omega)$, with $\zeta(x)=e$, $X\in T_xM$,  
 $Y\in \mathfrak{X}(M)$ and $\eta \in \Gamma(N_FM)$, 
 since $\mathcal{T}\zeta\in \Gamma$ then \eqref{eq:tilG} follows from 
\begin{align*}
\tilde{G}_*(e)Z &= \tilde{F}_*X + \tilde{\nabla}_{X} (\tilde{F}_*Y + T\eta)\\
&= \tilde{F}_*\left(X + \nabla_{X}Y - A_{T\eta}^{\tilde{F}}X\right) + \alpha^{\tilde{F}}_L(X,Y) 
+ (^{\tilde{F}}\nabla_{X}^\perp T\eta	)_L\\
&= T\left(F_*X + \nabla_{X}Y - A_{\eta}^FY  + \alpha^F(X,Y) + ^F\nabla_{X}^\perp \eta\right)\\
&=T{G}_*(e)Z.
\end{align*}

Now define $H\colon G({U}) \subset \mathbb{L}^{n+3} \to \mathbb{L}^{n+4}$  
by $H(G(e))=\tilde{G}(e)$ for any $e\in U$.
Then $H$ is an isometric immersion and  $\tilde F = H \circ F$. 
Define an open set in $\mathbb{V}^{n+2}$ by $W = G({U}) \cap \mathbb{V}^{n+2}$.  
Because $F(M^n) \subset G({U})$ and 
$F\colon M^n \to \mathbb{V}^{n+2}\subset \mathbb{L}^{n+3}$, 
it is clear that $F(M^n) \subset W$.  The only thing left to prove is that 
$H(W) \subset \mathbb{V}^{n+3}$.  
To see this, choose local sections $\delta_1$, $\delta_2$ of $\Gamma$ such that 
$\{\tilde{F}, \delta_1, \delta_2\}$ is a frame for $\Gamma$.  
Then $\{F, \bar{\delta}_1, \bar{\delta}_2\}$, where $\mathcal{T}(\bar{\delta}_i) = \delta_i$, 
is a frame for $\Omega$. From the definition of $G$ and because $G(U)$ 
is a tubular neighborhood of $F(M^n)$, 
we may write $G\colon U\times I^3 \to \mathbb{L}^{n+3}$ as
$$G(x,t,s_1,s_2) = (1+t)F(x) + s_1\bar{\delta}_1+s_2\bar{\delta}_2$$
and $\tilde{G}\colon U\times I^3 \to \mathbb{L}^{n+4}$ as
$$\tilde{G}(x,t,s_1,s_2) = (1+t)\tilde{F}(x) + s_1\delta_1+s_2\delta_2,$$
where $I$ is an interval containing zero.  Since $\tilde{G} = H \circ G$, we have
$$\big<\delta_1,\delta_2\big> = \big<\tilde{G}_*\partial_{s_1},\tilde{G}_*\partial_{s_2}\big> 
= \big<G_*\partial_{s_1},G_*\partial_{s_2}\big> = \big<\bar{\delta}_1,\bar{\delta}_2\big>$$
and
$$\big<\tilde{F},\delta_i\big> = \big<\tilde{G}_*\partial_t,\tilde{G}_*\partial_{s_i}\big> = 
\big<G_*\partial_t,G_*\partial_{s_i}\big> = \big<F,\delta_i\big>.$$
Hence, $ \big<H(G),H(G)\big>=\big<\tilde{G},\tilde{G}\big> = \big<G,G\big>$, 
which implies that $H(W) \subset \mathbb{V}^{n+3}$. \vspace{1ex}
\qed

We will also need the following slightly more general version of Lemma \ref{le:blem}.

\begin{lemma} \label{le:blem2} Let $\tilde{F}\colon M^n\to \mathbb{V}^{n+3} \subset \mathbb{L}^{n+4}$ 
be an isometric immersion of a Riemannian manifold.   Assume that there exists 
$\xi \in \Gamma(N_{\tilde F}M)$ of unit length such that  $\big<\xi,\tilde{F}\big> = 0$,  
$\text{rank}\, A_\xi^{\tilde{F}} = 1$ and  $^{\tilde{F}}\nabla_Z^\perp \xi = 0$ 
for all $Z \in \ker A_\xi^{\tilde{F}}$.
 Suppose further that the vector subbundle $L = \{\xi\}^\perp$,  the connection on $L$  
 induced by the normal connection of $\tilde F$, and the $L$-valued symmetric bilinear form 
 $\alpha_L = \pi_L \circ \alpha^{\tilde{F}}$ satisfy the Gauss, Codazzi and Ricci equations 
 for an isometric immersion of $M^n$ into  $\mathbb{L}^{n+3}$.  Then there exist, locally, 
 isometric immersions $F\colon M^n \to \mathbb{V}^{n+2} \subset \mathbb{L}^{n+3}$ 
 and $H\colon W\subset \mathbb{V}^{n+2} \to \mathbb{V}^{n+3}$ with $F(M) \subset W$, 
 such that $\tilde{F} = H \circ F$.
\end{lemma}

\proof  Since the assertion is of local nature, we may assume 
that $M^n$ is simply connected.
By the Fundamental Theorem of Submanifolds, there exist  an isometric immersion 
$F\colon M^n \to \mathbb{L}^{n+3}$ and a vector bundle isometry 
$\phi\colon L \to N_FM$ such that 
\begin{equation}\label{e:compositionsecondfund}
\alpha^F = \phi \circ \alpha^{\tilde{F}}_L \quad \text{and} \quad ^F\nabla^\perp \phi
 = \phi (^{\tilde F}\nabla^\perp)_L.
\end{equation}
Since $\big<\xi,\tilde{F}\big> = 0$,
the position vector field $\tilde F$ is a section of $L$.  Hence
$$\tilde{\nabla}_X \phi(\tilde F) = -F_* A^F_{\phi(\tilde F)}X 
+ ^F\nabla_X^\perp \phi(\tilde F) = F_*X.$$
Therefore, the section $F - \phi(\tilde F)$ is constant, say,  
$F - \phi(\tilde F)=P_0 \in \mathbb{L}^{n+3}$.  
Since $\phi$ is a vector bundle isometry and $\tilde F$ is a light-like section, 
it follows that $F - P_0 \in \mathbb{V}^{n+2}$.  
Without loss of generality we may assume that $P_0=0$, and so $\phi(\tilde{F})=F$.  

Define $T\colon N_FM \to L$ by $T \circ \phi = I$. 
Since $N_FM$ and $L$ have the same dimension and 
$T\colon N_FM \to L$, $\phi\colon L \to N_FM$ 
are vector bundle isometries with $T \circ \phi = I$, 
we have $\phi \circ T = I$.  Then
$$\phi(^{\tilde F}\nabla^\perp T)_L = 
\,^F\nabla^\perp (\phi \circ T) = \,^F\nabla^\perp$$
and $TF = \tilde{F}$.  Moreover, applying $T$ 
to both sides of the last equation, we get
$$(^{\tilde F}\nabla^\perp T)_L = T ( ^F\nabla^\perp),$$
which means that $T$ is parallel in the induced connection.  
From \eqref{e:compositionsecondfund} we get
$$\alpha^{\tilde{F}}(X,Y) = \pi_L \circ \alpha^{\tilde{F}}(X,Y) 
+ \left<A_\xi X,Y\right>\xi = T \circ \alpha^F(X,Y) +  \left<A_\xi X,Y\right>\xi.$$  
We finish by applying the previous lemma to $F|_V$, 
where $V \subset M^n$ is an open neighborhood 
of a given point of $M^n$ such that  $F|_V$ is an embedding.
\qed

\subsection{Structure of the second fundamental form}

Let $f\colon M^n \to \mathbb{R}^{n+1}$ be an isometric immersion  
with a nowhere vanishing principal curvature $\lambda$ of multiplicity $n-2$. 
Assume that $f$ is not a Cartan hypersurface and admits a 
genuine conformal deformation $\tilde{f}\colon M^n \to \mathbb{R}^{n+2}$. 
Our aim in this section is to describe the structure of the second fundamental 
form of the isometric light-cone representative of $\tilde{f}$. 

 We will make use of the following basic result on 
 flat bilinear forms known as the Main Lemma (see \cite{CD1987}). 
 Recall that a bilinear form $\beta\colon V \times V \to W$ is \emph{flat} 
 with respect to an inner product on $W$ if for all $X,Y,Z, W\in V$ we have
\begin{align*}
\left<\beta(X,Y),\beta(Z,W)\right> &- \left<\beta(X,W),\beta(Z,Y)\right> 
= 0.
\end{align*}
  
\begin{lemma}
\label{l:mainlemma}
Let $\beta\colon V^n \times V^n \to W^{p,q}$ be a symmetric flat 
bilinear form such that $\mathcal{S}(\beta) = W^{p,q}$.  
If $p \leq 5$ and $p+q < n$, then
$$\dim \mathcal{N}(\beta) \geq \dim V - \dim W = n - p - q,$$
where $\mathcal{N}(\beta)=\{Y\in V\,:\,\beta(X,Y)=0\;\;\mbox{for all}\;\; X\in V\}$.  
\end{lemma}

     The remaining of this section is devoted to proving the following result.

\begin{proposition}\label{p:commoneigen}
Let	$f\colon M^{n} \to \mathbb{R}^{n+1}$, $n \geq 6$, be an oriented hypersurface 
with a nowhere vaninhing principal curvature  $\lambda$ 
of constant multiplicity $n-2$ 
that is  not  a Cartan hypersurface on any open subset of $M^n$. 
Assume that $f$  admits a genuine conformal deformation 
$\tilde{f}\colon M^{n} \to \mathbb{R}^{n+2}$ and let 
$\tilde{F}=\mathcal{I}(\tilde{f})\colon M^{n} \to \mathbb{V}^{n+3} \subset \mathbb{L}^{n+4}$ 
be its isometric light-cone representative. Then,
 for each $x \in  M^n$ there exist  a space-like vector $\mu \in N_{\tilde{F}}M(x)$ 
 of  unit length and a flat bilinear form 
 $\gamma\colon T_xM \times T_xM \to \text{span}\{\mu\}^\perp$ such that
\begin{equation}\label{e:formatF}
\alpha^{\tilde{F}}(X,Y) = \left<AX,Y\right>\mu + \gamma(X,Y) 
\end{equation}
for all $X$,$Y \in T_xM$.  Moreover, $\lambda = -\big<\mu,\tilde{F}\big>^{-1}$  
and $\mathcal{N}(\gamma)$ coincides with the eigenspace 
$E_\lambda=\mbox{ker}\, (A-\lambda I)$ of $\lambda$ at $x$. 
\end{proposition}
\proof
Differentiating $\tilde{F} = \varphi^{-1}(\Psi \circ \tilde{f})$ we get 
$$\tilde{F}_*X = X(\varphi^{-1})(\Psi \circ \tilde{f}) + \varphi^{-1} \Psi_* \tilde{f}_* X.$$
Thus, the normal bundle $N_{\tilde{F}}M$ of $\tilde{F}$ splits orthogonally as
\begin{equation}\label{eq:ortdec}
N_{\tilde{F}}M = \Psi_{*}N_{\tilde{f}}M \oplus \mathbb{L}^2
\end{equation}
where $\mathbb{L}^2$ is a Lorentzian plane bundle having the position vector field 
$\tilde{F}$ as a section.  Hence, there exist unique sections $\xi$ and $\eta$ 
of $\mathbb{L}^2$ such that
$\left<\xi,\xi\right> = -1$, $\left<\xi,\eta\right>=0$, $\left<\eta,\eta\right> = 1$
and $\tilde{F}=\xi + \eta$.
At any $x \in M^n$, endow $W(x) = N_{f}M(x) \oplus N_{\tilde{F}}M(x)$ with the 
indefinite metric of type $(2,3)$ given by
$$\left<\!\left<\,,\right>\!\right>_{W(x)} = 
\left<\,,\right>_{N_fM(x)} - \left<\,,\right>_{N_{\tilde{F}}M(x)}.$$ 
Define a symmetric bilinear form by
$$\beta = \alpha^{f}\oplus \alpha^{\tilde{F}}\colon T_xM \times T_xM \to W(x).$$
From
\begin{equation}\label{eq:atildef}
\big<\alpha^{\tilde{F}}(X,Y),\tilde{F}\big> = -\left<X,Y\right>
\end{equation}
we deduce that $\mathcal{N}(\alpha^{\tilde{F}}) = \{0\}$, 
and hence $\mathcal{N}(\beta)=\{0\}$, 
for $\mathcal{N}(\beta) \leq \mathcal{N}(\alpha^{\tilde{F}})$.  
Moreover, the Gauss equations for $f$ and $\tilde{F}$ imply that
$\beta$ is flat with respect to $\left<\!\left<\,,\right>\!\right>$.

From Lemma \ref{l:mainlemma} for the case $(p,q) = (2,3)$, and since  $n\geq 6$, 
it follows that $\mathcal{S}(\beta)$ is degenerate,  
that is, the isotropic vector subspace
$\Omega = \mathcal{S}(\beta) \cap \mathcal{S}(\beta)^\perp$
is non-trivial.  Since the inner-product $\left<\!\left<\,,\right>\!\right>$ 
is positive definite on $W_1 = \text{span}\{N,\xi\}$ and negative definite 
on $W_2 = \Psi_{*}N_{\tilde{f}}M \oplus \text{span}\{\eta\}$,
 the orthogonal projections $P_1\colon W \to W_1$ and 
 $P_2\colon W \to W_2$ map $\Omega$ isomorphically onto 
 $P_1(\Omega)$ and $P_2(\Omega)$, respectively. Since
$\dim \mathcal{S}(\beta) + \dim \mathcal{S}(\beta)^\perp = 5,$
it follows that $\dim \Omega = 1$ or $\dim \Omega = 2$. 
Our first step is to show that our assumption that $\tilde f$ is a genuine 
conformal deformation of $f$ implies that the second possibility 
can not occur at any point of $M^n$.

Assume first that there is an open subset $U\subset M^n$ where $\dim \Omega = 2$ 
and that $\beta$ is null, that is, 
$\mathcal{S}(\beta)\subset \mathcal{S}(\beta)^\perp$.
Since $P_1|_\Omega$ is an isomorphism onto $W_1$ along $U$, due to dimensional reasons, 
there exists $\zeta \in \Omega$ be such that $P_1(\zeta) = \xi$.  Therefore $\zeta$ 
is a light-like vector in $\mathcal{S}(\alpha^{\tilde{F}})^\perp$.  
Moreover,  $\tilde{F}$ and  $\zeta_2 = \big<\zeta,\tilde{F}\big>^{-1}\zeta$  are 
linearly independent  by \eqref{eq:atildef}, with $\big<\zeta_2,\tilde{F}\big> = 1$.  
Let $\nu \in \Omega$ be such that $P_1(\nu) = N$.  Then $\nu = N + \tilde{\mu}$, 
where $\tilde{\mu} \in N_{\tilde{F}}U$ is a space-like vector of unit length.  From
$$0 = \big<\beta(X,Y), N + \tilde{\mu}\big> = \big<\alpha^f(X,Y),N\big> 
- \big<\alpha^{\tilde{F}}(X,Y),\tilde{\mu}\big>,$$ 
we conclude that $A=A_N$ coincides with $A_{\tilde{\mu}}^{\tilde{F}}$.  
Because $\nu$, $\zeta \in \Omega$, we have 
$0=\left<\nu,\zeta\right>= \left<\tilde{\mu},\zeta\right>= \left<\tilde{\mu},\zeta_2\right>$.
Define $\mu = \tilde{\mu} - \big<\tilde{\mu},\tilde{F}\big>\zeta_2$ 
and choose a space-like vector $\zeta_1 \in \{\mu,\zeta_2,\tilde{F}\}^\perp$ 
of unit length. Then  $\{\mu,\zeta_1,\zeta_2,\tilde{F}\}$ is a pseudo-orthonormal frame 
with respect to which the second fundamental  of $\tilde{F}$ is given by 
\begin{equation}\label{eq:sfftilF}
\alpha^{\tilde{F}}(X,Y) = \left<AX,Y\right>\mu 
+ \left<A_{\zeta_1}X,Y\right>\zeta_1 - \left<X,Y\right>\zeta_2.
\end{equation}  
Since $\beta$ is null,  we must have  $A_{\zeta_1} = 0$.
From the Codazzi equations of $f$ and $\tilde{F}$ for $A = A_\mu$ we get
$$\left<\nabla_X^\perp \mu,\zeta_2\right>Y = \left<\nabla_Y^\perp \mu,\zeta_2\right>X.$$
Hence $\left<\nabla_X^\perp \mu,\zeta_2\right> = 0$.  
From the Codazzi equation for $A_{\zeta_1} = 0$, we arrive to
$$\left<\nabla_X^\perp\zeta_1,\mu\right>AY - 
\left<\nabla_X^\perp\zeta_1,\zeta_2\right>Y =
\left<\nabla_Y^\perp\zeta_1,\mu\right>AX - \left<\nabla_Y^\perp\zeta_1,\zeta_2\right>X.$$  
Picking an orthonormal frame of eigenvectors $X_1,\cdots,X_n$ of $A$ 
correspondent to the principal curvatures
$\lambda_1,\cdots,\lambda_n$, respectively, with 
$\lambda_1= \cdots=\lambda_{n-2} = \lambda \neq 0$, 
we obtain for $i \neq j$ that
$\lambda_j\left<\nabla_{X_i}^\perp\zeta_1,\mu\right> 
= \left<\nabla_{X_i}^\perp\zeta_1,\zeta_2\right>,$
 hence $\left<\nabla_{X_i}^\perp\zeta_1,\zeta_2\right>=0
=\left<\nabla_{X_i}^\perp\zeta_1,\mu\right>$ 
for  $i=1,\cdots,n$.  Therefore $\mu$,$\zeta_1$,$\zeta_2$ 
and $\tilde{F}$ are parallel normal sections. 

Let $\bar f\colon U\to \mathbb{R}^{n+2}$ be the composition of 
$f|_U$ with a totally geodesic inclusion 
$i\colon  \mathbb{R}^{n+1}\to \mathbb{R}^{n+2}$. Then the second fundamental form of 
its isometric light-cone representative 
$\bar F \colon U\to \mathbb{V}^{n+3} \subset \mathbb{L}^{n+4}$ is 
$$\alpha^{\bar F}(X,Y)=\big<AX, Y\big>\Psi_*i_*N-\big<X,Y\big>w.$$
Let $\bar N$ be a unit normal vector field to $i$ along $f|_U$. 
Then, the vector bundle isometry $\tau\colon N_{\bar F}U\to N_{\tilde F}U$ given by
$$\tau\Psi_*i_*N=\mu, \,\,\,\tau \Psi_*\bar N=\zeta_1, 
\,\,\,\tau \bar F=\tilde F\,\,\,\mbox{and}\,\,\,\tau w=\zeta_2$$
is parallel and satisfies $\tau \alpha^{\bar F}=\alpha^{\tilde F|_U}$. 
It follows from the Fundamental Theorem of Submanifolds  
that $\tilde F|_U$ and $\bar F$ are congruent, 
and hence  $\tilde f|_U$ is conformally congruent to 
$\bar f = i \circ f|_U$ by Proposition \ref{p:conformalisometric}, 
which contradicts the assumption that $\tilde f$ is a genuine conformal deformation of $f$.

Now assume that there is an open subset $U\subset M^n$ where 
$\dim \Omega = 2$ and $\beta$ is not null. 
As in the previous case, there exists a pseudo-orthonormal frame 
$\{\mu,\zeta_1,\zeta_2,\tilde{F}\}$ with respect to which the 
second fundamental form of $\tilde{F}$ is given by \eqref{eq:sfftilF}, but now, 
since the bilinear form $\langle A_{\zeta_1}\;, \;\rangle$ is flat and $\beta$ is not null, 
we must have  $\dim \ker A_{\zeta_1} = n-1$.    
 From the Codazzi equation for $A=A_\mu$ we get
$$\left<\nabla_X^\perp \mu, \zeta_1\right>A_{\zeta_1}Y - \left<\nabla_X^\perp \mu,\zeta_2\right>Y 
= \left<\nabla_Y^\perp \mu, \zeta_1\right>A_{\zeta_1}X -\left<\nabla_Y^\perp \mu,\zeta_2\right>X.$$
For $X$, $Y \in \ker A_{\zeta_1}$ we conclude that $\ker A_{\zeta_1} \leq \ker \omega_2$, 
where $\omega_i$, $i=1,2$,  are the one-forms defined by $\omega_i(Y) 
= \left<\nabla_Y^\perp\mu,\zeta_i\right>$.   If $X \in \ker A_{\zeta_1}$ and 
$Y$ is an eigenvector of $A_{\zeta_1}$ with respect to the  unique non-zero eigenvalue, we get
$$\left<\nabla_X^\perp \mu, \zeta_1\right>A_{\zeta_1}Y 
 = -\left<\nabla_Y^\perp \mu,\zeta_2\right>X.$$
Therefore, $\omega_2 = 0$ and $\ker A_{\zeta_1} \leq \ker \omega_1$.

Let $F\colon M^n \to \mathbb{V}^{n+2}\subset \mathbb{L}^{n+3}$ be the isometric light-cone 
representative of the hypersurface $f\colon M^n\to \mathbb{R}^{n+1}$, 
whose second fundamental form is given by
$$\alpha^F(X,Y)=\left<AX,Y\right>\Psi_*N -\left<X, Y\right>w$$
for all $X, Y\in \mathfrak{X}(M)$. Define a vector bundle isometry  
$T\colon N_{F}U \to L = \{\zeta_1\}^\perp$ by 
$$T(F)=\tilde F, \;\;\;T(\Psi_*N)=\mu\;\;\;\mbox{and}\;\;\;T(w)=\zeta_2.$$
Then the second fundamental forms of $F|_U$ and $\tilde F|_U$  are related by
$$\alpha^{\tilde F}=T\circ \alpha^{F}+\left<A_{\zeta_1}\cdot, \cdot\right>\zeta_1.$$
Moreover, using that $\omega_2 = 0$ one can easily check that  $T$ is parallel 
with respect to the induced connection on $L$.  Since $\ker A_{\zeta_1} \leq \ker \omega_1$, 
it follows from Lemma \ref{le:blem} that, restricted to any open  subset $U_1\subset U$
where $F$ is an embedding, $\tilde F|_{U_1}$ is a composition 
$\tilde F|_{U_1}=H\circ F|_{U_1}$ of $F|_{U_1}$ with an isometric immersion 
$H\colon W\subset \mathbb{V}^{n+2} \to \mathbb{V}^{n+3}$ with $F(U_1) \subset W$.  
By Proposition \ref{p:conformalchar}, there exists a conformal immersion 
$h\colon V\to \mathbb{R}^{n+p}$ of an open subset $V \supset f(U_1)$ of 
$\mathbb{R}^{n+1}$ such that $\tilde f|_{U_1}=h \circ f|_{U_1}$,
contradicting the assumption that $\tilde f$ is a genuine conformal deformation of $f$.

In summary,  the subspace $\Omega$ must be one-dimensional at any point of $M^n$. 
The next step is to show that $\beta$ can not be  null at any point of $M^n$. 
Assume otherwise that $\beta$ is null at  $x \in M^n$. If 
$\Omega=\mathcal{S}(\beta)$  projects onto $\text{span}\{\xi\}$ under $P_1$, 
then $A=0$, a contradiction. Suppose now that $P_1(\Omega) \neq \text{span}\{\xi\}$.  
This is equivalent to requiring that the orthogonal projection 
$\Pi_1\colon W \to N_fM$ map $\Omega$ isomorphically onto $N_fM$, say, 
${N = \Pi_1(\nu)}$ for some $\nu \in \Omega$.  Set $\mu = \Pi_2(\nu)$, 
where $\Pi_2\colon W \to N_{\tilde{F}}M$ is 
the orthogonal projection onto $N_{\tilde{F}}M$.  
Then  $A=A_\mu^{\tilde{F}}$, for $N + \mu = \nu \in \Omega = 
\mathcal{S}(\beta) \subset \mathcal{S}(\beta)^\perp$, and hence
$$\beta(X,Y) = \big(\alpha^f(X,Y),\alpha^{\tilde{F}}(X,Y)\big)= 
\big(\left<AX,Y\right>N,\left<AX,Y\right>\mu\big).$$
Therefore,
$$-\big<X,Y\big> = \big<\alpha^{\tilde{F}}(X,Y),\tilde{F}\big> = 
\big<AX,Y\big>\big<\mu,\tilde{F}\big>,$$
contradicting the fact that the principal curvature 
$\lambda$ has multiplicity $n-2$. Thus $\beta$ is not null.

We now show that there is no point of $M^n$ where  
$P_1(\Omega) = \text{span}\{\xi\}$.
Suppose otherwise that  $P_1(\Omega) = \text{span}\{\xi\}$ at $x$. 
Then,  a light-like vector $\zeta$ spanning $\Omega$ belongs to 
$\mathcal{S}(\alpha^{\tilde{F}})^\perp$, and from \eqref{eq:atildef}
 it follows that $\tilde{F} \notin \Omega$. Thus
 $\zeta_2 = \big<\zeta,\tilde{F}\big>^{-1}\zeta$  and $\tilde{F}$ 
 form a  pseudo-orthonormal frame 
 of a Lorentzian plane $L$, and the $L$-component  of the 
 second fundamental form of $\tilde{F}$ is given by
  $\alpha^{\tilde{F}}_L(X,Y) = - \left<X,Y\right>\zeta_2$.  Hence
  $$\alpha^{\tilde{F}}(X,Y) = \left<A_{\zeta_0}X,Y\right>\zeta_0 
  + \left<A_{\zeta_1}X,Y\right>\zeta_1- \left<X,Y\right>\zeta_2,$$
  where $\{\zeta_0,\zeta_1,\zeta_2,\tilde{F}\}$ 
  is a pseudo-orthonormal basis of $N_{\tilde{F}}M(x)$.  
  Since $\dim \Omega =1$, the bilinear form
  $\hat{\beta}\colon T_xM\times T_xM \to \text{span}\{N,\zeta_0,\zeta_1\}$
  defined by 
  $$\hat{\beta} = \alpha^f \oplus \big<\alpha^{\tilde{F}},\zeta_0\big>\zeta_0
   \oplus \big<\alpha^{\tilde{F}},\zeta_1\big>\zeta_1$$
  is flat and nondegenerate, hence $\dim \mathcal{N}(\hat{\beta}) \geq n-3$ 
  by Lemma \ref{l:mainlemma}.  From
  $$\mathcal{N}(\hat{\beta}) = \ker A \cap \ker A_{\zeta_0} \cap \ker A_{\zeta_1}$$
  it follows that  $\lambda$  must be zero, contradicting the assumption.

Therefore  $P_1(\Omega) \neq \text{span}\{\xi\}$ at any point of $M^n$.  
Then, as in the case when $\beta$ was 
assumed to be null, there exists $\nu \in \Omega$ such that $\nu=N+\mu$,
 with $\mu$ of  unit length and $A = A_\mu$. Hence
  $$\alpha^{\tilde{F}}(X,Y) = \left<AX,Y\right>\mu + \gamma(X,Y) $$
  for  $\gamma\colon T_xM\times T_xM \to \{\mu\}^\perp$
 a  flat nondegenerate bilinear form.  Thus  ${\mathcal{N}(\gamma) \geq n-3}$ 
 by Lemma \ref{l:mainlemma}.  If $T \in \mathcal{N}(\gamma)$, then
  $$-\big<T,Y\big> = \big<\alpha^{\tilde{F}}(T,Y),\tilde{F}\big> 
  = \big<AT,Y\big>\big<\mu,\tilde{F}\big>,$$
  hence $\big<\mu,\tilde{F}\big>$ is non-zero and 
  $\lambda = -\big<\mu,\tilde{F}\big>^{-1}$, with $\mathcal{N}(\gamma) \leq E_\lambda$. 
  To complete the proof, it remains to show that  $\mathcal{N}(\gamma) = E_\lambda$, 
  that is, $\dim \mathcal{N}(\gamma) = n-2$.

  Assume, by contradiction, that  $\Delta= \mathcal{N}(\gamma)$ has dimension $n-3$ 
  on some open subset $U\subset M^n$.  
  We will prove that $\tilde{f}|_U = h \circ g$, 
  where $g\colon U \to \mathbb{R}^{n+1}$ is a genuine conformal deformation of $f$ and 
  $h\colon V\subset \mathbb{R}^{n+1} \to \mathbb{R}^{n+2}$ is a conformal immersion 
  of an open subset containing $f(U)$. In particular, it will follow that 
  $f|_U$ is a  Cartan hypersurface, contradicting our assumption. 
  
  Defining $\zeta = \lambda \tilde{F} + \mu$, we have $\left<\zeta,\zeta\right> = -1$, 
  $\left<\zeta,\mu\right> = 0$ and  $A_\zeta = A - \lambda I$.
  Therefore, if $T \in E_\lambda \cap \Delta^\perp$, then 
  $0 \neq \gamma_T(T_xM)  \leq \text{span}\{\mu, \zeta\}^\perp$ at any $x\in U$. 
 We claim that $\gamma_T(T_xM)$ has dimension one. 
 Assume otherwise,  and let $X \in \ker \gamma_T \cap \Delta^\perp$. Then
  $$0 = \left<\gamma(T,X),\gamma(Z,W)\right> = \left<\gamma(T,W),\gamma(X,Z)\right>$$
  for all $Z$, $W \in T_xM$ by  the flatness of  $\gamma$, and hence 
  $\gamma_X(T_xM) \leq \text{span}\{\zeta\}$.  Notice that 
   $\gamma_X(T_xM)$ can not be trivial, for $X\in \Delta^\perp$, 
   thus $\gamma_X(T_xM) = \text{span}\{\zeta\}$.
   Using again the flatness of $\gamma$, we obtain  that
   $\gamma_Y(T_xM)\leq \{\zeta\}^\perp$, or equivalently, 
   $Y\in \mbox{ker}\, A_\zeta=E_\lambda$, for all 
   $Y\in \mbox{ker}\, \gamma_X$. 
   This contradicts the fact that $\lambda$ 
   has multiplicity $n-2$ and proves the claim.
   
   Let $\{\mu, \zeta_1,\zeta_2,\zeta\}$ be an orthonormal frame of 
   $N_{\tilde F}U$ with $\gamma_T(T_xM) = \text{span}\{\zeta_1\}$ for all $x\in U$.
   Flatness of $\gamma$ now implies that $X\in \text{ker} \gamma_T$ 
   if and only if $\gamma_X(T_xM)\leq \{\zeta_1\}^\perp$, that is,  
   if and only if $X\in \text{ker} A_{\zeta_1}$. Thus $\text{rank}\,A_{\zeta_1} = 1$. 
   Moreover, since  $\text{rank}\,A_\zeta = \text{rank}\,(A-\lambda I) = 2$ 
   and $\gamma$ is nondegenerate, we must have 
  \begin{equation}\label{e:cartancond}
  A_{\zeta_2} \neq \pm A_\zeta.
  \end{equation}
  
  Define the symmetric bilinear form 
  $$\hat{\gamma} = \gamma - \left<\gamma,\zeta_1\right>\zeta_1 = 
  \left<\gamma,\zeta_2\right>\zeta_2 
  - \left<\gamma,\zeta\right>\zeta\colon T_xM \times T_xM \to \text{span}\{\zeta_2,\zeta\}.$$
  Using that $\text{rank}\,A_{\zeta_1} = 1$, 
  from the flatness and nondegeneracy of $\gamma$ it follows easily 
  that $\hat{\gamma}$ is also flat and nondegenerate.  
  By  Lemma \ref{l:mainlemma},  we have that $\dim \mathcal{N}(\hat{\gamma}) \geq n-2$, 
  and since $\mathcal{N}(\hat{\gamma}) \leq \ker A_{\zeta_2}$, 
  it follows that  $\text{rank}\,A_{\zeta_2}\leq 2$.   
  If $\text{rank}\,A_{\zeta_2}\leq 1$, then  
  $\hat{\gamma} - \left<\gamma,\zeta_2\right>\zeta_2 = -\left<\gamma,\zeta\right>\zeta$ 
  would be flat.  Also, it is nondegenerate, because $\zeta$ is a time-like unit vector. 
  Thus,   Lemma \ref{l:mainlemma} would imply that $\dim \mathcal{N}(A_\zeta)\geq n-1$, 
  which is impossible, because $A_\zeta = A- \lambda I$ has rank two.  
  Therefore, $\text{rank}\,A_{\zeta_2} = 2$.  Also,
  since  $\mathcal{N}(\hat{\gamma}) 
  = \ker A_{\zeta_2} \cap \ker A_\zeta$ and $\dim \ker A_{\zeta_2} 
  = \dim \ker A_{\zeta} =n-2$,  we must have  $\ker A_{\zeta_2} = \ker A_{\zeta}$.  
  Observe also that $\ker A_{\zeta_2}$ can not be contained in $\ker A_{\zeta_1}$, 
  because $\Delta=\ker A_{\zeta_2} \cap \ker A_{\zeta_1}$ 
  has dimension $n-3$. Equivalently, $\text{Img}\,A_{\zeta_1} \cap \text{Img}\,A_{\zeta_2} = \{0\}$.

 From the Codazzi equation for $A_\mu = A$ we have that
  $A_{\nabla_X^\perp \mu}Y = A_{\nabla_Y^\perp \mu}X$, and taking into consideration that
   $\nabla_X^\perp\zeta = X(\lambda) \tilde{F} + \nabla_X^\perp\mu$ we get
  \begin{align*}
  \left<\nabla_X^\perp\mu,\zeta_1\right>A_{\zeta_1}Y 
  &+ \left<\nabla_X^\perp \mu,\zeta_2\right>A_{\zeta_2}Y - \lambda^{-1}X(\lambda)A_{\zeta}Y \\
  &=\left<\nabla_Y^\perp\mu,\zeta_1\right>A_{\zeta_1}X
   + \left<\nabla_Y^\perp \mu,\zeta_2\right>A_{\zeta_2}X - \lambda^{-1}Y(\lambda)A_{\zeta}X.
  \end{align*}
  For $Y=R \in \Delta$ and $X \in \ker A_{\zeta_1}\cap E_\lambda$, the preceding equation gives
  \begin{equation}\label{e:chypmuzeta2delta}
  \left<\nabla_R^\perp \mu,\zeta_2\right> = 0 \quad \text{for}\,\, R \in \Delta.
  \end{equation}
  For $Y=R \in \Delta$, $X \in (\ker A_{\zeta_1})^\perp$ 
  and using  \eqref{e:chypmuzeta2delta} we obtain
  \begin{equation}\label{e:chypmuzeta1delta}
  \left<\nabla_R^\perp\mu,\zeta_1\right>= 0, \quad \text{for}\,\,R \in \Delta.
  \end{equation}

  Using that $\left<\nabla_X^\perp \zeta_1,\mu\right> = \left<\nabla_X^\perp \zeta_1,\zeta\right>$ 
  and $A_{\zeta} = A - \lambda I$,
  the Codazzi equation for $A_{\zeta_1}$ gives
  \begin{align*}
  \nabla_X &A_{\zeta_1} Y - A_{\zeta_1}\nabla_X Y - 
  \lambda\left<\nabla_X^\perp \zeta_1,\mu\right>Y - 
  \left<\nabla_X^\perp \zeta_1,\zeta_2\right>A_{\zeta_2}Y 	\\
  &= \nabla_Y A_{\zeta_1}X - A_{\zeta_1}\nabla_Y X - 
  \lambda\left<\nabla_Y^\perp \zeta_1,\mu\right>X - 
  \left<\nabla_Y^\perp \zeta_1,\zeta_2\right>A_{\zeta_2}X.  \end{align*}
  For $Y=R \in \Delta$ and $X \in \ker A_{\zeta_1}$ and using \eqref{e:chypmuzeta1delta} we get
  $$-A_{\zeta_1}\nabla_X R - \lambda\left<\nabla_X^\perp\zeta_1,\mu\right>R = - A_{\zeta_1}\nabla_R X -  \left<\nabla_R^\perp\zeta_1,\zeta_2\right>A_{\zeta_2}X,$$
 hence
  \begin{equation}\label{e:chypzeta1mu}
  \left<\nabla_X^\perp\zeta_1,\mu\right> = 0 \quad \text{for}\,\, X \in \ker A_{\zeta_1}.
  \end{equation} 
  Now, for $X$, $Y \in  \ker A_{\zeta_1}$, and using \eqref{e:chypzeta1mu}, we have
  $$- A_{\zeta_1}\nabla_X Y- \left<\nabla_X^\perp \zeta_1,\zeta_2\right>A_{\zeta_2}Y 
  =  - A_{\zeta_1}\nabla_Y X - \left<\nabla_Y^\perp \zeta_1,\zeta_2\right>A_{\zeta_2}X,$$
  thus
  \begin{equation}\label{e:chypzeta1zeta2}
  \left<\nabla_X^\perp \zeta_1,\zeta_2\right> = 0 \quad \text{for}\,\, X \in \ker A_{\zeta_1}.
  \end{equation}
  It follows from  \eqref{e:chypzeta1mu}, \eqref{e:chypzeta1zeta2} and 
  $\left<\nabla_X^\perp \zeta_1,\mu\right> = \left<\nabla_X^\perp \zeta_1,\zeta\right>$  
  that $\zeta_1$ is parallel along $\ker A_{\zeta_1}$.  
  
  Define the rank-$3$ subbundle $L$ by $L = \{\zeta_1\}^\perp$.  
  Since $A_{\zeta_1}$ has rank $1$, the $L$-component $\alpha_L^{\tilde{F}}$ 
  satisfies the Gauss equations for an isometric immersion of 
  $U$ into $\mathbb{L}^{n+3}$. We now show that 
  $(\alpha^{\tilde{F}}_L, (\nabla^\perp)_L)$ also satisfies the  Codazzi and Ricci equations.
    
  The Codazzi equation for $A_\mu = A$ with respect to $(\nabla^\perp)_L$ reduces to 
  $$\left<\nabla_X^\perp \mu, \zeta_2\right>A_{\zeta_2}Y 
  - \left<\nabla_X^\perp\mu,\zeta\right>A_\zeta Y 
  = \left<\nabla_Y^\perp \mu, \zeta_2\right>A_{\zeta_2}X 
  - \left<\nabla_Y^\perp\mu,\zeta\right>A_\zeta X.$$
  Because $A_{\nabla_X^\perp \mu}Y = A_{\nabla_Y^\perp \mu}X$, it suffices to show that
  $$\left<\nabla_X^\perp\mu,\zeta_1\right>A_{\zeta_1}Y 
  = \left<\nabla_Y^\perp\mu,\zeta_1\right>A_{\zeta_1}X.$$
  But this holds because   $\dim\ker A_{\zeta_1} = n-1$ 
  and $\zeta_1$ is parallel along $\ker A_{\zeta_1}$.
  The other Codazzi equations are proved in a similar way.

  Let us move on to the Ricci equations.  
  Using the Ricci equation for $\tilde F$ involving $\mu$ and $\zeta_2$, 
  the corresponding one for the pair $(\alpha^{\tilde{F}}_L, (\nabla^\perp)_L)$ reduces to   
  $$\left<\nabla_X^\perp \zeta_1,\zeta_2\right>\left<\nabla_Y^\perp\mu,\zeta_1\right>- \left<\nabla_X^\perp\mu,\zeta_1\right>\left<\nabla_Y^\perp\zeta_1,\zeta_2\right> = 0,$$
  which is true because $\dim \ker A_{\zeta_1} = n-1$ and $\zeta_1$ is parallel along 
  $\ker A_{\zeta_1}$.  The remaining Ricci equations for 
  $(\alpha^{\tilde{F}}_L, (\nabla^\perp)_L)$   follow in a similar way.
  
  By Lemma \ref{le:blem2},  there exist, locally,  isometric immersions 
  $G\colon  M^n \to \mathbb{V}^{n+2} \subset \mathbb{L}^{n+3}$ and 
  $H\colon W\subset \mathbb{V}^{n+2} \to \mathbb{V}^{n+3}$ with 
  $G(M) \subset W$, such that $\tilde{F} = H \circ G$. 
  By Lemma \ref{p:conformalchar}, there exist, locally,   
  conformal immersions $g\colon M^n\to \mathbb{R}^{n+1}$ 
  and  $h\colon V\to \mathbb{R}^{n+2}$, of an open subset  
  $V\subset \mathbb{R}^{n+1}$ containing $g(M)$, such that $\tilde f=h\circ g$.

We  now argue that $g$ is a genuine conformal deformation of $f$.  
Suppose, on the contrary,  that $f$ and $g$ are conformally congruent.  
Then, from Proposition \ref{p:conformalchar}, their isometric light-cone 
representatives $F$ and $G$ are isometrically congruent, that is, 
there exist an isometry $T:\mathbb{L}^{n+3} \to \mathbb{L}^{n+3}$ 
such that $G = T \circ F$.  Since the second fundamental form of $G$ is
$$\alpha^G(X,Y) = \left<AX,Y\right>\mu + \left<A_{\zeta_2}X,Y\right>\zeta_2 
- \left<A_\zeta X,Y\right>\zeta,$$
and that of $F$ is
$$\alpha^F(X,Y) = \left<AX,Y\right>\Psi_*N - \left<X,Y\right>w,$$
it is easy to see that the condition $\alpha^G = T \circ \alpha^F$ would imply that 
$A_{\zeta_2} = \pm A_\zeta$, a contradiction with  \eqref{e:cartancond}. \qed

\section{The triple $(D_1, D_2, \psi)$}
\label{ch:thetriple}

The aim of this section is to show that, for a hypersurface $f\colon M^n \to \mathbb{R}^{n+1}$ 
that carries a nowhere vanishing principal curvature  of constant multiplicity $n-2$ 
and is not a Cartan hypersurface on any open subset of $M^n$, the existence of a genuine 
conformal deformation ${\tilde{f}\colon M^n \to \mathbb{R}^{n+2}}$ is equivalent to $f$ being 
a  hyperbolic or elliptic hypersurface on which one can define  a pair of tensors $D_1, D_2$ 
and a one-form $\psi$  satisfying certain conditions. Before giving a precise statement 
(Proposition~\ref{p:equivalentcartanhyp} below), we need some definitions.  

Let $f\colon M^n \to \mathbb{R}^{n+1}$ be a hypersurface that carries a principal curvature 
of multiplicity $n-2$, let $\Delta$ denote the corresponding eigenbundle, and let
$$C\colon \Gamma(\Delta) \to \Gamma(\text{End}(\Delta^\perp))$$
be its splitting tensor, defined by
$$C_T X = -\nabla_X^h T$$
for all $T \in \Gamma(\Delta)$ and $X\in \Gamma(\Delta^\perp)$, where the superscript 
$h$ denotes taking the component in $\Delta^\perp$.  The hypersurface $f$ is said to be 
\emph{hyperbolic}\index{hyperbolic} (respectively, \emph{parabolic} or \emph{elliptic}) 
if there exists $J \in \Gamma(\text{End}(\Delta^\perp))$ satisfying the following conditions:
\begin{enumerate}[(i)]
\item $J^2 = I$ (respectively, $J^2 = 0$, with $J\neq 0$,  and $J^2 = -I$).
\item $\nabla_T^hJ = 0$ for all $T\in \Gamma(\Delta)$.
\item $C_T \in \text{span}\{I,J\}$ for all $T\in \Gamma(\Delta)$.
\end{enumerate}

A hypersurface $f\colon M^n\to \mathbb{R}^{n+1}$, $n \geq 3$, is called
\emph{conformally surface-like}\index{conformally surface-like} if $f(M)$ is the image 
by a M\"{o}bius transformation of $\mathbb{R}^{n+1}$ of an open subset of one of the following:
\begin{enumerate}
\item a cylinder $M^2 \times \mathbb{R}^{n-2}$ over a surface $M^2 \subset \mathbb{R}^3$;
\item a cylinder $CM^2 \times \mathbb{R}^{n-3}$, where $CM^2 \subset \mathbb{R}^4$ 
denotes the cone over $M^2 \subset \mathbb{S}^3$;
\item a rotation hypersurface over a surface $M^2 \subset \mathbb{R}^3_+$.
\end{enumerate}

We will need the following characterization of conformally surface-like hypersurfaces, 
which is a consequence of a more general result in \cite{mD2001} (see also \cite{rT2006}). 

\begin{proposition}\label{c:confsurfacelike}
A hypersurface ${f\colon M^n \to \mathbb{R}^{n+1}}$ is \emph{conformally surface-like} 
if and only if it has a principal curvature $\lambda$ of multiplicity $n-2$ 
whose eigendistribution $\Delta=\mbox{ker}(A-\lambda I)$ 
has the property that the distribution $\Delta^\perp$ is umbilical. 
\end{proposition}

  In the remaining of this section we prove the following result.

\begin{proposition}\label{p:equivalentcartanhyp}
Let $f\colon M^n \to \mathbb{R}^{n+1}$ be an oriented  hypersurface with a nowhere vanishing 
principal curvature $\lambda$ of constant multiplicity $n-2$. 
Assume that $f$ is not a  Cartan hypersurface 
on any open subset of $M^n$ and that it  admits a genuine conformal deformation 
${\tilde{f}\colon M^n \to \mathbb{R}^{n+2}}$. Then, on each connected component 
of an open dense subset, $f$ is either hyperbolic or elliptic with respect to a tensor 
$J \in \Gamma(\text{End}(\Delta^\perp))$, where $\Delta=\ker (A-\lambda I)$,  
and there exists a unique (up to signs and permutation) pair $(D_1,D_2)$ of tensors in 
$\Gamma(\text{End}(\Delta^\perp))$, with $D_i \in \text{span}\{I,J\}$ for $i=1,2$,
and a unique one-form $\psi$ on $M^n$ satisfying the following conditions:
	\begin{enumerate}[(i)]
	\item $\Delta \leq \ker \psi$,
	\item $\det D_i = \frac{1}{2}$,
	\item $\nabla_T^h D_i = 0 = [D_i,C_T]$ for all $T\in \Delta$,
	\item $(\nabla_X (A-\lambda I)D_i)Y - (\nabla_Y (A-\lambda I)D_i)X \\ 
	= (X\wedge Y)D^t_i \text{grad}\,\lambda + (-1)^j(A-\lambda I)\left(\psi(X)D_jY 
	- \psi(Y)D_jX\right)$,
	\item $\left<(\nabla_Y D_i)X - (\nabla_X D_i)Y,\text{grad}\,\lambda\right> 
	+ \text{Hess}\,\lambda(D_i X,Y) - \text{Hess}\,\lambda(X,D_iY) \\
	+ (-1)^j\psi(X)\left<D_jY,\text{grad}\,\lambda\right> 
	- (-1)^j\psi(Y)\left<D_jX,\text{grad}\,\lambda\right> \\
	= \lambda\left( \left<AX,(A-\lambda I)D_iY\right> 
	- \left<(A-\lambda I)D_iX,AY\right> \right)$,  
	\item $\text{d}\psi(Z,T)=0$ for all $Z \in \mathfrak{X}(M)$ and $T \in \Delta$,
	\item $\text{d}\psi(X,Y) = \left<[(A-\lambda I)D_1,(A-\lambda I)D_2]X,Y\right>$.  
	\item $D_2^2 \neq \pm D_1^2$.
	\item $\text{rank}\, (D_1^2+D_2^2-I)=2$.
	\end{enumerate}
	
	Conversely, let $f\colon M^n \to \mathbb{R}^{n+1}$ be a simply connected hypersurface 
	that is not conformally surface-like and carries a nowhere vanishing principal curvature 
	of constant multiplicity $n-2$.    If $f$ is  hyperbolic or elliptic  with respect to 
	${J \in \text{End}(\Delta^\perp)}$, where $\Delta=\mbox{ker}(A-\lambda I)$, 
	and there exist a triple $(D_1,D_2,\psi)$	satisfying items (i)-(ix), 
	with $D_i \in \text{span}\{I,J\}$ for $i=1,2$, then  $f$  admits a genuine conformal 
	deformation 	$\tilde{f}\colon M^n \to \mathbb{R}^{n+2}$.  Moreover, distinct triples 
	(up to sign and permutation) yield non conformally congruent conformal deformations.  
\end{proposition} 

\proof Let $\tilde{F}\colon M^n \to \mathbb{V}^{n+3}\subset \mathbb{L}^{n+4}$ be the 
isometric light-cone representative of $\tilde f$.  
For each $x \in M^n$, let $\mu \in N_{\tilde{F}}M(x)$ and 
$\gamma:T_xM \times T_xM \to \text{span}\{\mu\}^\perp$ be 
given by Proposition \ref{p:commoneigen}. 
Then, the vector field $\zeta = \lambda \tilde{F} + \mu$ satisfies
$$\left<\zeta,\zeta\right>=-1,\;\;\; \left<\zeta,\mu\right>=0
\;\;\;\mbox{and}\;\;\;A_\zeta = A - \lambda I.$$
Consider the Riemannian plane-bundle $\mathbb{P} = \{\zeta,\mu\}^\perp$.  
For each $\xi \in \Gamma(\mathbb{P})$, define
$$D_\xi = (A- \lambda I)^{-1}A_\xi = A_\zeta^{-1}A_\xi \in \Gamma(\text{End}(\Delta^\perp))$$ 
where all endomorphisms are considered restricted to $\Delta^\perp$, and let 
$$W = \text{span}\{D_\xi: \xi \in \Gamma(\mathbb{P})\}.$$  

\begin{lemma}\label{l:dimW}
The map $\xi \in \mathbb{P}(x) \to D_\xi\in W(x)$ is an isomorphism 
for all $x$ in  an open dense subset of $M^n$. 
\end{lemma}

\proof 
Suppose there exists a non-trivial $\tilde{\rho} \in \Gamma(\mathbb{P})$ 
on an open subset $U\subset M^n$ such that $D_{\tilde{\rho}}=0$, 
and hence $A_{\tilde{\rho}} = 0$.
Decompose $\tilde{\rho} = \Psi_{*} \rho + \rho_1$,  
with $\rho \in \Gamma(N_{\tilde{f}}U)$ and $\rho_1 \in \Gamma(\mathbb{L}^2)$,  
according to  the orthogonal decomposition \eqref{eq:ortdec} of $N_{\tilde F}U$. 
Since $\tilde{\rho}$ and $\tilde{F}$ are orthogonal, we have  
$\rho_1=\big<\rho_1, \tilde \zeta\big>\tilde F $, where $\{\tilde \zeta, \tilde F\}$ 
is a pseudo-orthonormal frame of $\mathbb{L}^2$ 
with $\langle\tilde \zeta, \tilde F\rangle=1$.  Because the 
$\Psi_*N_{\tilde{f}}U$-component of $\alpha^{\tilde{F}}$ is 
$\varphi^{-1}\Psi_*\alpha^f$, from $A_{\tilde{\rho}}= 0$ we get
$$0 = \varphi^{-1}\big<A_\rho X,Y\big> - \big<X,Y\big>\big<\tilde{\zeta},\rho_1\big>,$$  
for all $X$, $Y \in \mathfrak{X}(U)$.  In particular, since $\tilde{\rho}$ is not trivial, 
the normal vector field $\rho$ can not be trivial either.  
We conclude that $A_\rho = \beta I$, with $\beta=\varphi \big<\tilde{\zeta},\rho_1\big>$. 
If $\rho$ is parallel in the normal connection, then $\tilde{f}(U)$ is contained 
in either an affine hyperplane or a hypersphere of $\mathbb{R}^{n+2}$, according 
to whether $\beta$ vanishes or not. But this implies  $f$ to be a Cartan hypersurface, 
contrary to our assumption. Otherwise,  $U$ is conformally flat by Theorem 14 
in \cite{mD1996} if $\beta \neq 0$, and flat by an elementary computation using the 
Codazzi equation if $\beta=0$. Both possibilities  contradict the assumption that 
$\lambda$ is nowhere vanishing and  has multiplicity  $n-2$. \vspace{1ex}\qed

 We will need the following properties of the tensors $D_\xi$.

\begin{lemma}\label{l:propofD}
The following holds:
\begin{enumerate}[(i)]
\item $[D_\xi,C_T] = 0$ for all $T \in \Gamma(\Delta)$.
\item $\nabla_T^h D_\xi = 0$ for all $T \in \Gamma(\Delta)$ if 
$\xi \in \Gamma(N_{\tilde{F}}M)$ is parallel along $\Delta$.
\end{enumerate}
\end{lemma}

\proof
Using the Codazzi equation we obtain
\begin{equation}\label{e:covderhorA}
(\nabla_T^hA)X=(A - \lambda I)C_T X
\end{equation}
and 
\begin{equation}\label{e:covderhorAxi}
(\nabla_T^hA)(X,\xi)=A_\xi C_T X
\end{equation}
for all $X \in \Gamma(\Delta^\perp)$.  
In particular, $(A - \lambda I)C_T$ and $A_\xi C_T$ are symmetric. Therefore
$$
(A-\lambda I)D_\xi C_T = A_\xi C_T = C_T^t A_\xi
= C_T^t (A-\lambda I)D_\xi = (A-\lambda I)C_T D_\xi,	
$$
which proves $(i)$, because  $A - \lambda I$ is an isomorphism on $\Delta^\perp$.  
If $\xi \in \Gamma(N_{\tilde{F}}M)$ is parallel along $\Delta$, then
 $$
(A-\lambda I)D_\xi C_T = A_\xi C_T = \nabla_T^h A_\xi 
= \nabla_T^h (A-\lambda I)D_\xi = \nabla_T^h AD_\xi - \lambda \nabla_T^h D_\xi.
$$
On the other hand, from  \eqref{e:covderhorA} we also have
$(A-\lambda I)C_T D_\xi = (\nabla_T^hA)D_\xi$. 
We get $(ii)$ by subtracting the preceding identities and using $(i)$:
$$
0=(A-\lambda I)[D_\xi,C_T] = A \nabla_T^h D_\xi - \lambda \nabla_T^h D_\xi 
= (A - \lambda I)\nabla_T^h D_\xi. \qed
$$

\begin{lemma}\label{l:generatorW}
There exists $J\in \Gamma(\text{End}(\Delta^\perp))$ such that 
$J^2 = \epsilon I$, $\epsilon \in \{1,0,-1\}$, and 
$$\rm{span}\{I\} < C(\Gamma(\Delta)) \leq \rm{span}\{I,J\} = W.$$	
\end{lemma}
\proof Since $f$ is not conformally surface-like on any open subset of $M^n$, 
otherwise it would be a Cartan hypersurface on that subset, 
by Corollary \ref{c:confsurfacelike} the distribution 
$\Delta^\perp$ is not umbilical, and hence  
$C(\Gamma(\Delta)) \neq \text{span}\{I\}$. Let 
$$S = \{A \in \text{End}(\Delta^\perp): AB=BA \,\, \text{for} \, B \in W\}.$$
 Part (i) of Lemma \ref{l:propofD} says that $C(\Gamma(\Delta)) \leq S$.  
 Since $\dim W = 2$ by  Lemma \ref{l:dimW}, we must have  
 $I \in W$, for otherwise  we would have $S = \text{span}\{I\}$, 
 a contradiction.  Therefore, $W=\rm{span}\{I, J\}$, where $J$ is a  
 tensor on $\Delta^\perp$ satisfying $J^2=\epsilon I$,
$\epsilon\in \{-1, 1, 0\}$.
 In particular, $W\subset S$
and, on the other hand,  the fact that any element of $S$ commutes with  $J$  
implies that the dimension of $S$ is at most two. 
Hence $W=S$ and $C(\Gamma(\Delta))\subset S=\rm{span}\,\{I, J\}$. \vspace{1ex}\qed

Now consider  any orthonormal frame $\{\tilde\xi_1,\tilde\xi_2\}$ of $\mathbb{P}$ 
and define the one-forms
\begin{equation*}
\tilde{\psi}(X) = \left<\nabla_X^\perp \tilde\xi_1,\tilde\xi_2\right>, \quad \tilde{\omega}_1(X) 
= \left<\nabla_X^\perp \tilde\xi_1,\mu\right> \quad \text{and} \quad \tilde{\omega}_2(X) 
= \left<\nabla_X^\perp \tilde\xi_2,\mu\right>.
\end{equation*}
Using that	
$\nabla_X^\perp (\zeta - \mu) = X(\lambda) \tilde{F} = \lambda^{-1}X(\lambda)(\zeta-\mu)$
 for all $X \in \mathfrak{X}(M)$, we obtain
\begin{equation}\label{e:norcovderxi}
\nabla_X^\perp \tilde\xi_1 =  \tilde{\omega}_1(X) (\mu - \zeta)+ \tilde{\psi}(X)\tilde\xi_2, 
\end{equation}
\begin{equation}\label{e:norcovderxi2}
\nabla_X^\perp \tilde\xi_2 = \tilde{\omega}_2(X)(\mu - \zeta) - \tilde{\psi}(X)\tilde\xi_1,
\end{equation}
\begin{equation}\label{e:norcovdermu}
\nabla_X^\perp \mu = -\tilde{\omega}_1(X)\tilde\xi_1 - 
\tilde{\omega}_2(X)\tilde\xi_2 -\lambda^{-1}X(\lambda)\zeta 
= \nabla_X^\perp \zeta -\lambda^{-1}X(\lambda)(\zeta-\mu).
\end{equation}

Straightforward computations using \eqref{e:norcovderxi}, \eqref{e:norcovderxi2}, 
\eqref{e:norcovdermu} and the Codazzi and Ricci equations of $F$ show that, for all 
$X, Y\in \Gamma(\Delta^\perp)$, 
\begin{equation}\label{e:codmu2}
(X\wedge Y)\text{grad}\,\lambda = D_{\tilde\xi_1}\left(\lambda\tilde{\omega}_1(X)Y 
- \lambda\tilde{\omega}_1(Y)X\right)
 + D_{\tilde\xi_2}\left(\lambda\tilde{\omega}_2(X)Y- \lambda\tilde{\omega}_2(Y)X \right),
\end{equation} 
while, for all $X$, $Y \in \mathfrak{X}(M)$  and $1 \leq i \neq j \leq 2$,  
\begin{equation}\label{e:codxi}\begin{array}{l}
	(\nabla_X A_{\tilde\xi_i})Y - (\nabla_Y A_{\tilde\xi_i})X 
	=\vspace{1ex}\\
	\hspace{12ex}
	 \lambda \left(\tilde{\omega}_i(X)Y -  \tilde{\omega}_i(Y)X\right) 
	 + (-1)^j\left(\tilde{\psi}(X)A_{\tilde\xi_j}Y - \tilde{\psi}(Y)A_{\tilde\xi_j}X\right),
\end{array}
\end{equation}
\begin{align}\label{e:riccimuxi}
\left<[A_\mu,A_{\tilde\xi_i}]X,Y\right> &= 	
-d\tilde{\omega}_i(X,Y) + (-1)^j\tilde{\omega}_j(Y)\tilde{\psi}(X) 
- (-1)^{j}\tilde{\omega}_j(X)\tilde{\psi}(Y) \\
&\quad+ \lambda^{-1}Y(\lambda)\tilde{\omega}_i(X) 
- \lambda^{-1}X(\lambda)\tilde{\omega}_i(Y),\nonumber
\end{align}
\begin{equation}\label{e:riccixixi}
\left<[A_{\tilde\xi_1},A_{\tilde\xi_2}]X,Y\right> = \text{d}\tilde{\psi}(X,Y)
\end{equation}
and
\begin{align}\label{e:riccixizeta}
\left<[A_{\tilde\xi_i},A_\zeta]X,Y\right> &= 	d\tilde{\omega}_i(X,Y) 
+ \lambda^{-1}\tilde{\omega}_i(Y)X(\lambda) + (-1)^{j}\tilde{\psi}(Y)\tilde{\omega}_j(X)\\
&\quad- \lambda^{-1}\tilde{\omega}_i(X)Y(\lambda)
 - (-1)^{j}\tilde{\psi}(X)\tilde{\omega}_j(Y).\nonumber
\end{align}

\begin{lemma}\label{l:proptensorw}
For any orthonormal frame $\{\tilde\xi_1,\tilde\xi_2\}$ of $\mathbb{P}$ we have
\begin{itemize}
\item[(i)] 
$1 = \det D_{\tilde\xi_1} + \det D_{\tilde\xi_2}$.
\item[(ii)] $\Delta \leq \ker \tilde{\omega}_1 \cap \ker \tilde{\omega}_2$.
\item[(iii)]  $D_i = D_{\tilde\xi_i}$, for $i=1,2$, satisfy
 $D_2^2 \neq -D_1^2$.
 \end{itemize}
\end{lemma}

\proof
$(i)$ Flatness of  $\gamma$ means that
$\det A_\zeta = \det A_{\tilde\xi_1} + \det A_{\tilde\xi_2}$.\vspace{1ex}\\
$(ii)$ Using \eqref{e:norcovdermu}, 
the Codazzi equation
$0=A_{\nabla^\perp_X\mu}Y- A_{\nabla^\perp_Y\mu}X$
applied to  $X \in \Gamma(\Delta^\perp)$  and $Y=T \in \Gamma(\Delta)$ yields
$\tilde{\omega}_1(T)D_{\tilde\xi_1} + \tilde{\omega}_2(T)D_{\tilde\xi_2} = 0.$	
Thus $\tilde{\omega}_1(T) = 0=\tilde{\omega}_1(T)$, for   $D_{\tilde\xi_1}$ 
and $D_{\tilde\xi_2}$ are linearly independent by Lemma~\ref{l:dimW}.\vspace{1ex}\\
$(iii)$ Suppose, by contradiction, that $D_2^2 = - D_1^2$.  
In view of Lemma \ref{l:generatorW},  we may write $D_1 = aI + bJ$ and $D_2 = cI + dJ$ for some
$a,b,c,d \in C^{\infty}(M)$.  Then 
$$(c^2 + \epsilon d^2)I + 2cdJ = - (a^2 + \epsilon b^2)I - 2abJ.$$  
Thus $a=b=c=d=0$ if $\epsilon = 1$, and  $a=c=0$ if $\epsilon = 0$,  a contradiction.  

If $\epsilon=-1$, denote by $\hat{D}_i=\theta_i I+\bar{\theta}_i \hat J$ 
the complex linear extension of $D_i$, $1\leq i\leq 2$, where $\hat J$ 
is the complex linear extension of $J$.  
From $D_2^2 = - D_1^2$ we get $\theta_2^2 = - \theta_1^2$, and we may assume that 
$\theta_2 = i \theta_i$.  
From part $(i)$ we get $1= 2|\theta_1|^2$, so we can write 
$$\sqrt{2}\hat{D}_1 =\theta I +\bar{\theta} \hat J\;\;\;\mbox{and}\;\;\;\sqrt{2}\hat{D}_2 
=i\theta I -i\bar{\theta} \hat J$$
for some $\theta\in \mathbb{S}^1$.
 Writing $\theta = e^{i\beta}$, we have
 $$ \sqrt{2}D_1 =\cos \beta I + \sin \beta J \;\;\;\mbox{and}\;\;\;\sqrt{2}\hat{D}_2 
 = -\sin \beta I + \cos \beta J.$$
Then, the orthonormal frame $\{\xi, \eta\}$ of $\mathbb{P}$ defined by
$$\xi = \cos \beta \tilde\xi_1 -  \sin \beta \tilde\xi_2 
\quad \text{and} \quad \eta= 
\sin \beta \tilde\xi_1 + \cos \beta \tilde\xi_2$$
satisfies 
$\sqrt{2}D_\xi = I$ and $\sqrt{2}D_\eta = J$.  
 Using  \eqref{e:codxi} with $\tilde\xi_1 = \xi$ and 
 $\tilde\xi_2 = \eta$ yields
\begin{equation}\label{e:contradeq}
Y(\lambda)X-X(\lambda)Y=\sqrt{2}\,\lambda(\tilde{\omega}_1(X)Y-\tilde{\omega}_1(Y)X)
+\sqrt{2}\,\tilde{\psi}(X)A_{\eta}Y-\sqrt{2}\,\tilde{\psi}(Y)A_{\eta}X
\end{equation}
for all $X, Y\in \mathfrak{X}(M)$. 
For  $Y = T \in \Gamma(\Delta)$ and $X\in \Gamma(\Delta^\perp)$,  using part $(ii)$ 
we obtain
$$\sqrt{2}\tilde{\psi}(T)A_\eta X  = \big(X(\lambda) 
+ \sqrt{2}\lambda\tilde{\omega}_1(X)\big)T,$$
hence
$$\Delta \leq \ker \tilde{\psi} \quad \text{and} \quad X(\lambda) 
+ \sqrt{2}\lambda\tilde{\omega}_1(X) = 0, \quad \text{for}\, \, X \in \Delta^\perp.$$
Substituting the last identity in  \eqref{e:contradeq} for $X$ and 
$Y \in \Delta^\perp$ gives $\tilde{\psi}(Y)A_\eta X = \tilde{\psi}(X)A_\eta Y$,
hence $\tilde{\psi} = 0$.  From  \eqref{e:riccixixi} we obtain
$$ \left<[A_\xi,A_\eta]X,Y\right> = \text{d}\tilde{\psi}(X,Y) = 0,$$
hence $[(A-\lambda I),(A- \lambda I)J] = 0$.  This means that $A$ and $J$ commute, 
hence $A - \lambda I = \beta I$ in $\Delta^\perp$,
with $\beta \neq 0$.  Using the identity $(A-\lambda I)C_T = \nabla_T^h A$, we get 
\begin{align*}
\beta C_T &= \nabla_T^h (\beta + \lambda)I= T(\beta + \lambda) I,
\end{align*}
 a contradiction because $f$ is not conformally surface-like.  
 Thus $D_2^2 \neq -D_1^2$.  \vspace{1ex}\qed

The next lemma shows that the Riemannian plane bundle $\mathbb{P}$ 
has a distinguished orthonormal frame $\{\xi_1,\xi_2\}$.

\begin{lemma}\label{l:particularframe}
There exists a unique (up to sign and permutation) orthonormal frame 
$\{\xi_1,\xi_2\}$ of $\mathbb{P}$ such that $D_i = D_{\xi_i}$,  $i=1,2$, satisfy
$$\det D_1 = \frac{1}{2} = \det D_2.$$
Moreover,  $\xi_1$ and $\xi_2$ are parallel along $\Delta$.  	
\end{lemma}

\proof 
Pick an arbitrary orthonormal frame $\{\xi,\eta\}$ for $\mathbb{P}$.  
Since $1 = \det D_\xi + \det D_\eta$ by  part $(i)$ of Lemma \ref{l:proptensorw}, 
we are done  if either $D_\xi$ or $D_\eta$ has determinant $1/2$.  
So, suppose  that $\det D_\xi < 1/2$ and $\det D_\eta > 1/2$.  Define 
 $\xi_1(\theta) = \cos \theta \,\xi + \sin \theta \,\eta$ and $\xi_2(\theta)
  = -\sin \theta \,\xi + \cos \theta \,\eta$, $\theta\in [0, \pi/2]$.  Since
$$ \det D_\xi =\det D_{\xi_1(0)} < \det D_{\xi_1(\frac{\pi}{2})} = \det D_\eta,$$
  existence follows  by continuity.  
  Uniqueness follows using part $(iii)$ of Lemma~\ref{l:proptensorw}.

 We now show that $\xi_1$ and $\xi_2$ are parallel along $\Delta$. 
Given $x\in M^n$, $T\in \Delta$ and an integral curve $\gamma$ of $T$ 
starting at $x$, let $\hat\xi_i(t)$  denote the parallel transport of
$\xi_i(x)$ along $\gamma$ at $\gamma(t)$. By Lemma \ref{l:propofD}, we have that
$\nabla_{\gamma'(t)}D_{\hat\xi_i(t)}=0$, hence $\det D_{\hat\xi_i(t)}=1/2$.
Since $\xi_1$ and $\xi_2$ are unique (up to signs and permutation) with this property,
by continuity we must have $\hat \xi_i(t)=\xi_i(\gamma(t))$ for any $t$. It
follows that $\nabla_T^\perp\xi_i=0$ for any $T\in \Delta$, $i=1,2$. 
\vspace{1ex}\qed

 From now on, we fix the privileged orthonormal frame $\{\xi_1,\xi_2\}$ of $\mathbb{P}$ 
 given by the  above lemma and omit the tilde notation in 
 $\omega_1$, $\omega_2$ and $\psi$ when  using this frame.  
 Also,  from now on $D_i$ stands for $D_{\xi_i}$, $i=1,2$.  
We will show that the pair $(D_1, D_2)$ and the one-form $\psi$ satisfy conditions 
$(i)$-$(ix)$ in the statement.

From Lemma \ref{l:proptensorw}, and because $\xi_1$ and $\xi_2$ are parallel along $\Delta$, we have
\begin{equation}\label{e:kerdelta}
\Delta \leq \ker \psi \cap \ker \omega_1 \cap \ker\omega_2.
\end{equation}
Thus, condition $(i)$ is satisfied.  Conditions $(ii)$ and $(iii)$ follow from  
Lemma~\ref{l:particularframe}  and Lemma \ref{l:propofD}, respectively.

From   \eqref{e:codxi} for $Y=T \in \Gamma(\Delta)$, a unit length section, 
and $X \in \Gamma(\Delta^\perp)$, we get
\begin{equation}\label{eq:lambom}
0 
= \lambda \omega_i(X)	T  + A_{\xi_i}\nabla_X T +\nabla_T A_{\xi_i}X - A_{\xi_i}\nabla_T X.
\end{equation}
Using that $\Delta$ is an umbilical distribution whose mean curvature vector field
$\delta$ is given by $(\lambda I -A)\delta = \text{grad}\,\lambda$ 
(see Eq. $2$ in the proof of Proposition $8$ of \cite{mD2001}), we obtain
$$
 \left<A_{\xi_i} X, \nabla_T T\right>\\
=  \left<(A-\lambda I)D_{i} X, \delta \right>\\
= -\left<D_iX, \text{grad}\, \lambda\right>.
$$
Therefore, taking the inner product with $T$ of both sides of (\ref{eq:lambom}) yields
\begin{equation}\label{e:omegaequivdef}
\omega_i(X) = -\frac{1}{\lambda}\left<D_i X,\text{grad}\, \lambda\right>.
\end{equation}
For $X$, $Y \in \Gamma(\Delta^\perp)$, we obtain from \eqref{e:codxi} that 
\begin{align*}
(\nabla_X (A-\lambda I)D_i)Y &- (\nabla_Y (A-\lambda I)D_i)X  \\
&= \lambda\left(\omega_i(X)Y - \omega_i(Y)X\right) 
+ (-1)^j\left(\psi(X)A_{\xi_j}Y - \psi(Y)A_{\xi_j}X\right),
\end{align*}
  From  \eqref{e:omegaequivdef} we get
\begin{align*}
	\lambda\left(\omega_i(X)Y - \omega_i(Y)X\right)
	&= \left<D_iY,\text{grad}\,\lambda\right>X - 
	\left<D_iX,\text{grad}\,\lambda\right>Y \\
	&= (X \wedge Y)D^t_i \text{grad}\,\lambda.
\end{align*}
Because $A_{\xi_j} = (A - \lambda I)D_j$, combining the last two equations 
gives item $(iv)$. Differentiating  \eqref{e:omegaequivdef} yields
$$
Y\omega_i(X) = -\lambda^{-1}Y(\lambda)	\omega_i(X) 
- \lambda^{-1}\left<\nabla_Y D_iX,\text{grad}\,\lambda\right> 
- \lambda^{-1}\text{Hess}\,\lambda(D_iX,Y).
$$
  Therefore, 
\begin{align*}
\text{d}\omega_i&(X,Y) - \lambda^{-1}Y(\lambda)\omega_i(X) 
+ \lambda^{-1}X(\lambda)\omega_i(Y) \\
&\quad =\text{d}\omega_i(X,Y) + Y\omega_i(X) 
+ \lambda^{-1}\left<\nabla_Y D_iX,\text{grad}\,\lambda\right> 
+ \lambda^{-1}\text{Hess}\,\lambda(D_iX,Y)
\\
&\quad \quad - X\omega_i(Y) 
- \lambda^{-1}\left<\nabla_X D_iY,\text{grad}\,\lambda\right>
 - \lambda^{-1}\text{Hess}\,\lambda(D_iY,X)\\
&\quad =\frac{1}{\lambda}\left( \left<(\nabla_Y D_i)X 
- (\nabla_X D_i)Y,\text{grad}\,\lambda\right> 
+ \text{Hess}\,\lambda(D_i X,Y) - \text{Hess}\,\lambda(X,D_iY)\right).
\end{align*}
Substituting the preceding expression 
in  \eqref{e:riccimuxi} and using again  
\eqref{e:omegaequivdef} yields $(v)$. Applying  \eqref{e:riccixixi} 
to $Y = T \in \Gamma(\Delta)$ yields $(vi)$, whereas  item $(vii)$ 
follows from the same equation applied to $X, Y\in \Gamma(\Delta^\perp)$.
We have from part $(iii)$ of Lemma \ref{l:proptensorw} that $D_2^2 \neq - D_1^2$.  
It is easily checked that $D_1$ and $D_2$ would be linearly dependent 
if  $D_2^2 = D_1^2$, so $(viii)$ is proven. \vspace{1ex}

The next lemma completes the proof that $f$ is hyperbolic, parabolic or elliptic 
with respect to $J \in \Gamma(\text{End}(\Delta^\perp))$ 
given by  Lemma \ref{l:generatorW}.

\begin{lemma}\label{l:hypparell}
The tensor $J$ satisfies $\nabla_T^h J = 0$. 
\end{lemma}

\proof 
Since $D_1$ and $D_2$ are linearly independent, we may assume 
that $D_1=a_1I+b_1J$, with $b_1\neq 0$.
By part $(ii)$ of Lemma \ref{l:propofD} we have
$$
0=(\nabla_T^hD_1)=T(a_1)I+T(b_1)J+b_1\nabla^h_TJ
$$
for any $T\in\Gamma(\Delta)$. Hence
$$
T(a_1)J+\epsilon\, T(b_1)I+b_1(\nabla^h_TJ)J=0\;\;\mbox{and}\;\;
T(a_1)J+\epsilon\, T(b_1)I+b_1J(\nabla^h_TJ)=0.
$$
Adding the two equations yields $T(a_1)=T(b_1)=0$, 
and hence $\nabla^h_TJ=0$. \qed
\vspace{1ex}

A hypersurface $f:M^n\to \mathbb{R}^{n+1}$, $n \geq 3$, is said to be 
\emph{conformally ruled} if it carries an umbilical distribution $L$ 
of rank $n-1$ such that the restriction of $f$ to each leaf of $L$ 
is also umbilical. We now prove that the  parabolic case occurs 
precisely when $f$ is conformally ruled.

\begin{lemma}\label{l:ruledparabolic}
Let $f:M^n \to \mathbb{R}^{n+1}$ be an oriented  hypersurface with a nowhere vanishing 
principal curvature of constant multiplicity $n-2$. Assume that $f$ is not a  Cartan 
hypersurface  on any open subset of $M^n$ and that it admits a genuine conformal 
deformation $\tilde{f}:M^n \to \mathbb{R}^{n+2}$.  If  $f$ is parabolic with respect 
to  $J \in \Gamma(\text{End}(\Delta^\perp))$, then it is conformally ruled.  
\end{lemma}

\proof
Pick an orthonormal frame $\{X,Y\}$  of $\Gamma(\Delta^\perp)$ such that $JY = 0$ 
and $JX = \delta Y$ with $\delta \neq 0$. We will prove that the distribution
$$L(x) = \Delta(x) \oplus Y(x)$$
is umbilical, that is,  there exists $\rho\in C^{\infty}(M)$  such that 
$\left<\nabla_U V,X\right> = \rho\left<U,V\right>$
  for all $U$,$V \in \Gamma(L)$.
  From $C_T \in \text{span}\{I,J\}$ and $JY = 0$ 
  we get ${\left<C_TY,X\right> = 0}$, hence
\begin{equation}\label{e:parabolic1}
\left<\nabla_Y T,X\right> = - \left<C_TY,X\right> = 0 
\quad \text{for all}\,\, T \in \Gamma(\Delta).
\end{equation}

Since $J\nabla^h_TY=(\nabla_T^h J)Y = 0$ by Lemma \ref{l:hypparell}, and 
$\nabla_T^h Y$ is orthogonal to $Y$, it follows that 
$\nabla_T^h Y=0$,  or  equivalently,
\begin{equation}\label{e:parabolic2}
\left<\nabla_T Y,X\right> = 0.	
\end{equation}

Using that $(A-\lambda I)C_T = \nabla_T^h A$  is symmetric and 
$\text{span}\{I\} < C(\Delta) \leq \text{span}\{I,J\}$, 
we conclude that $(A-\lambda I)J$ is symmetric.  Therefore,
\begin{equation}\label{e:parabolic0}
\left<(A-\lambda I)Y,Y\right> = \delta^{-1}\left<(A-\lambda I)JX,Y\right> 
=  \delta^{-1}\left<X,(A-\lambda I)JY\right> = 0.
\end{equation}
It follows that in the orthonormal frame $\{X,Y\}$ of $\Delta^\perp$ we have
\begin{equation}\label{e:parabolicshapeA}
A - \lambda I=
\begin{pmatrix}
\beta&\mu\\
\mu&0	
\end{pmatrix}
\end{equation}
with $\mu \neq 0$, for $A - \lambda I$ restricted to $\Delta^\perp$ is an isomorphism.  
Since $D_i \in \text{span}\{I,J\}$, with $\det D_i = 1/2$, and $D_1$ and $D_2$ 
are linearly independent,  we can suppose that 
\begin{equation}\label{e:ruledD}
\sqrt{2}D_i = I + b_iJ,
\end{equation}
with $b_1 \neq 0$. Therefore, 
$$\sqrt{2}A_{\xi_i}Y = (A-\lambda I)\sqrt{2}D_iY = (A-\lambda I)Y = \mu X$$
and
$$ \sqrt{2}A_{\xi_i}X = (A-\lambda I)\sqrt{2}D_iX = (A-\lambda I)(X + b_i \delta Y) 
= (\beta + b_i\delta\mu)X + \mu Y.$$
Define $\theta = b_1\delta \mu \neq 0$ and $\tilde{\theta} = b_2\delta \mu$, 
so in the orthonormal frame $\{X,Y\}$ we have
\begin{equation}\label{e:parabolicshapeAxi}
\sqrt{2}A_{\xi_1}=
\begin{pmatrix}
\beta + \theta&\mu\\
\mu&0	
\end{pmatrix}
\quad \text{and} \quad
\sqrt{2}A_{\xi_2}=
\begin{pmatrix}
\beta + \tilde{\theta}&\mu\\
\mu&0	
\end{pmatrix}.
\end{equation}

Applying the Codazzi equation of $A$  to $T \in \Gamma(\Delta)$ of unit length 
and $Y \in \Gamma(\Delta^\perp)$,  and then taking the inner product with $T$,
we obtain  using  \eqref{e:parabolicshapeA} that
\begin{equation}\label{e:parabolic3}
\mu\left<\nabla_T T,X\right> = - Y(\lambda).	
\end{equation}

Now, applying the Codazzi equation for $A$  to $X$, $Y \in \Gamma(\Delta^\perp)$, 
and then taking the inner product with $Y$ yields
\begin{equation}\label{e:parabolic4}
	0 = 2\mu\left<\nabla_X X,Y\right> + X(\lambda) +\beta\left<\nabla_Y Y,X\right> - Y(\mu).
\end{equation}

Next, applying the Codazzi equation for $A_{\xi_i}$, $1\leq i\leq 2$,  
to $X$, $Y \in \Gamma(\Delta^\perp)$, and using   \eqref{e:codxi} 
and \eqref{e:parabolicshapeAxi}, give, respectively,
 \begin{equation}\label{e:parabolic5}
0 = 2\mu\left<\nabla_X X,Y\right> + (\beta + \theta)\left<\nabla_Y Y,X\right> 
- Y(\mu) -\sqrt{2}\lambda \omega_1(X) +  \mu \psi(Y).	
\end{equation}
and
\begin{equation}\label{e:parabolic6}
0 = 2\mu \left<\nabla_X X,Y\right> + (\beta + \tilde{\theta})\left<\nabla_Y Y,X\right> 
- Y(\mu) - \sqrt{2}\lambda\omega_2(X) - \mu\psi(Y).	
\end{equation}

Replacing  \eqref{e:parabolic4} into  \eqref{e:parabolic5} and \eqref{e:parabolic6} we obtain
$$\theta\left<\nabla_Y Y,X\right> - X(\lambda) - \sqrt{2}\lambda\omega_1(X) + \mu\psi(Y) = 0$$
and
$$\tilde{\theta}\left<\nabla_Y Y,X\right> - X(\lambda) 
-  \sqrt{2}\lambda\omega_2(X) - \mu\psi(Y) = 0.$$
Adding both equations yields
$$(\theta + \tilde{\theta})\left<\nabla_Y Y,X\right> - 2X(\lambda) 
- \sqrt{2}\lambda\left(\omega_1(X) + \omega_2(X)\right)=0.$$
Using  \eqref{e:omegaequivdef} and \eqref{e:ruledD}, and that 
$(\theta + \tilde\theta) = (b_1+b_2)\delta\mu$, we get
\begin{equation}\label{e:parabolic7}
(\theta + \tilde{\theta})\left(\mu\left<\nabla_Y Y,X\right> + Y(\lambda)\right)= 0.
\end{equation}

Suppose  that $\theta +\tilde{\theta}=0$. From  \eqref{e:parabolicshapeA} and  
\eqref{e:parabolicshapeAxi}, the vector fields
$$\xi = \frac{1}{\sqrt{2}}\left(\xi_1 + \xi_2\right) \quad \text{and} 
\quad \eta = \frac{1}{\sqrt{2}}\left(\xi_1 - \xi_2\right)$$
define an orthonormal frame $\{\xi,\eta\}$ of $\mathbb{P}$ satisfying
\begin{equation}\label{e:parabolicshapeAxieta}
A_\xi =
\begin{pmatrix}
\beta&\mu\\
\mu&0	
\end{pmatrix}
=(A-\lambda I)
\quad \text{and} \quad
A_\eta=
\begin{pmatrix}
\theta&0\\
0&0	
\end{pmatrix}.
\end{equation}
In particular, $D_\eta=(A-\lambda I)^{-1}A_\eta$ satisfies 
$$D_\eta X = \frac{\theta}{\mu}Y\;\;\;\mbox{and}\;\;\;D_\eta Y = 0.$$
From  \eqref{e:codmu2} for $\xi_1=\xi$ and $\xi_2=\eta$ we obtain
$$Y(\lambda)X - X(\lambda)Y = \lambda\tilde{\omega}_1(X)Y 
- \lambda\tilde{\omega}_1(Y)X -\frac{\lambda\theta}{\mu}\tilde{\omega}_2(Y)Y.$$ 
Hence,
\begin{equation}\label{e:parabolic8}
\tilde{\omega}_1(Y) + \lambda^{-1}Y(\lambda) 
= 0 \quad \text{and} \quad \tilde{\omega}_1(X)
 + \lambda^{-1}X(\lambda)  -\frac{\theta}{\mu}\tilde{\omega}_2(Y) = 0.
\end{equation}
Now,  the Codazzi equation of $A_\xi = A - \lambda I$ yields 
$$(Z\wedge W)\text{grad}\,\lambda = \lambda \tilde{\omega}_1(Z)W 
+ \tilde{\psi}(Z)A_\eta W - \lambda\tilde{\omega}_1(W)Z - \tilde{\psi}(W)A_\eta Z.$$
$Z$, $W \in \mathfrak{X}(M)$. For $Z = T \in \Delta$ and $W = X$, 
using  \eqref{e:parabolicshapeAxieta} and Lemma \ref{l:proptensorw} we obtain
\begin{equation}\label{e:paraboliccodazzi2}
X(\lambda) = - \lambda\tilde{\omega}_1(X) \quad \text{and} \quad \Delta \leq \ker \tilde{\psi}.
\end{equation}
Replacing now $Z=X$ and $W=Y$ and using  \eqref{e:parabolicshapeAxieta} we get
$$(X\wedge Y)\text{grad}\,\lambda =\lambda \tilde{\omega}_1(X)Y  
- \lambda\tilde{\omega}_1(Y)X - \theta\tilde{\psi}(Y)X, $$
hence
\begin{equation}\label{e:paraboliccodazzi3}
Y(\lambda) = - \theta\tilde{\psi}(Y)  - \lambda \tilde{\omega}_1(Y) 
\quad \text{and} \quad -X(\lambda) = \lambda \tilde{\omega}_1(X).
\end{equation}
It follows from   \eqref{e:parabolic8}, \eqref{e:paraboliccodazzi2} 
and \eqref{e:paraboliccodazzi3} that 
\begin{equation}\label{e:paraboliccodazziconclusions}
\Delta \oplus \text{span}\{Y\} \leq \ker \tilde{\psi} \cap \ker \tilde{\omega}_2 
\quad \text{and} \quad \lambda^{-1}Z(\lambda) + \tilde{\omega}_1(Z) = 0, 
\,\, \text{for}\,\,Z \in \mathfrak{X}(M).
\end{equation}

Now, the second fundamental form of $\tilde{F}$ is given by
\begin{align*}
\alpha^{\tilde{F}}(X,Y) &= \left<AX,Y\right>\mu + \left<(A-\lambda I)X,Y\right>\xi 
+ \left<A_\eta X,Y\right>\eta - \left<(A-\lambda I)X,Y\right>\zeta\\
&= \left<AX,Y\right>(\mu + \xi - \zeta) - \lambda\left<X,Y\right>(\xi - \zeta) 
+ \left<A_\eta X,Y\right>\eta.	
\end{align*}
From  \eqref{e:norcovderxi}, \eqref{e:norcovdermu}
 and \eqref{e:paraboliccodazziconclusions} we get
\begin{equation}\label{e:bundleisometryinduced1}
\nabla_X^\perp (\mu + \xi - \zeta) = \lambda^{-1}X(\lambda)(\mu - \zeta) 
+ \tilde{\omega}_1(X)(\mu -\zeta) + \tilde{\psi}(X)\eta = \tilde{\psi}(X)\eta,
\end{equation}
while using  \eqref{e:norcovderxi}
and \eqref{e:paraboliccodazziconclusions} we get
\begin{align}\label{e:bundleisometryinduced2}
\nabla_X^\perp \lambda(\xi-\zeta) &= X(\lambda)(\xi - \zeta) 
+ \lambda \nabla_X^\perp(\xi-\zeta)\\
&= \lambda\big(\tilde{\psi}(X)+\tilde{\omega}_2(X)\big)\eta,\nonumber
\end{align}
for all $X \in \mathfrak{X}(M)$. On the other hand, the second fundamental form 
of the isometric light-cone representative 
$F\colon M^n\to \mathbb{V}^{n+2} \subset \mathbb{L}^{n+3}$ of $f$ is given by
\begin{align*}
\alpha^{{F}}(X,Y) &= \left<AX,Y\right>\Psi_*N -\left<X,Y\right>w.	
\end{align*}
Define a vector-bundle isometry $\tau\colon N_FM\to L=\{\eta\}^\perp$ by setting
$$\tau \Psi_*N=\mu + \xi - \zeta,\;\;\;\tau w=\lambda (\xi - \zeta)
\;\;\;\mbox{and}\;\;\;\tau F=\tilde F.$$
From  \eqref{e:bundleisometryinduced1} and \eqref{e:bundleisometryinduced2}, 
the vector bundle isometry is parallel with respect to the induced connection on $L$.  
By Lemma \ref{le:blem}, there exists an isometric immersion
 $H\colon W\subset \mathbb{V}^{n+2} \to \mathbb{V}^{n+3}$, with $F(M^n) \subset W$, 
 such that $\tilde F = H \circ F$. It follows from Proposition \ref{p:conformalchar} 
 that there exists a conformal immersion $h\colon V\to \mathbb{R}^{n+p}$ of an open subset 
 $V \supset f(M^n)$ of $\mathbb{R}^{n+1}$ such that $\tilde f=h \circ f$, 
 contradicting the assumption that $\tilde f$ is a genuine conformal deformation of $f$. 

Thus $(\theta + \tilde{\theta}) \neq 0$, and from  \eqref{e:parabolic1}, \eqref{e:parabolic2}, 
\eqref{e:parabolic3} and \eqref{e:parabolic7} it follows that $L$ is an 
umbilical distribution with mean curvature vector $Z = -(Y(\lambda)/\mu) X$.

It remains to prove  that the restriction $g=f \circ i\colon \sigma\to \mathbb{R}^{n+1}$ 
of $f$ to each leaf $\sigma$ of $L$ is also umbilical.    From  \eqref{e:parabolic0} we get 
$$\alpha^{g}(Y,Y) = f_*\alpha^i(Y,Y) + \alpha^f(i_*Y,i_*Y) = f_*Z + \lambda N,$$
whereas for all $T$, $S \in \Gamma(\Delta)$ we have
$$\alpha^{g}(T,S) = f_*\alpha^i(T,S) + \alpha^f(i_*T,i_*S) = \left<T,S\right>f_*Z 
+ \lambda \left<T,S\right>N.$$
Thus $g$ is umbilical with $f_*Z + \lambda N$ as its mean curvature vector field.
\vspace{1ex}\qed

 Since conformally ruled hypersurfaces are Cartan hypersurfaces (see \cite{mD2000}), 
 in view of Lemma \ref{l:ruledparabolic} the parabolic case is ruled out by the assumption. 
 Therefore, to complete the proof of the direct statement it remains to prove condition (ix).

\begin{lemma}\label{l:rankgenuinedef}
The tensors $D_1$ and $D_2$ satisfy
$$\text{rank}\,(D_1^2 + D_2^2 - I) = 2.$$	
\end{lemma}

\proof
We will argue separately for the  elliptic and  hyperbolic cases.

\subsubsection{Elliptic Case}

This case is almost trivial.  Write $D_1 = aI + bJ$ and $D_2 = cI + dJ$.  
Since $\det D_i = 1/2$, we  have ${a^2 + b^2 = c^2 + d^2 = 1/2}$, hence
\begin{equation*}
D_1^2 + D_2^2 - I = 
\begin{pmatrix}
a^2 - b^2 + c^2 - d^2 -1 &	2(ab + cd)\\
-2(ab + cd)&a^2 - b^2 + c^2 - d^2 -1
\end{pmatrix}.
\end{equation*}
The conclusion follows, 
for otherwise $D_1^2 + D_2^2 - I=0$, hence $b=0=d$ 
from $a^2 + b^2 + c^2 + d^2 = 1= a^2 - b^2 + c^2 - d^2$, 
contradicting the linear independence of $D_1$ and $D_2$.

\subsubsection{Hyperbolic Case}

Suppose that $\text{rank}\,(D_1^2 + D_2^2 - I) < 2$ and let 
\begin{equation}
\sqrt{2}D_1 =
\begin{pmatrix}
\theta_1 & 0 \\
0 & \theta_1^{-1}	
\end{pmatrix}
\quad \text{and} \quad
\sqrt{2}D_2 =
\begin{pmatrix}
\theta_2 & 0 \\
0& \theta_2^{-1}	
\end{pmatrix}.
\end{equation}
Then,
\begin{equation}
2D_1^2 + 2D_2^2 - 2I =
\begin{pmatrix}
\theta_1^2 + \theta_2^2 -2 & 0\\
0 & \theta_1^{-2} + \theta_2^{-2} -2
\end{pmatrix},	
\end{equation}
and we may assume that
$\theta_1^2 + \theta_2^2 = 2$.  Thus, the orthonormal frame  $\{\xi,\eta\}$ of $\mathbb{P}$ given by
$$\sqrt{2}\xi = \theta_1 \xi_1 + \theta_2\xi_2 \quad \text{and} 
\quad \sqrt{2}\eta = -\theta_2\xi_1 + \theta_1\xi_2$$
satisfies $D_\xi = I$ and $\text{rank} \,D_\eta=1$. Let $\{X, Y\}$ be an orthogonal frame of 
$\Delta^\perp$ with $D_\eta X=0$.  From  \eqref{e:codmu2} for $\tilde\xi_1=\xi$ 
and $\tilde \xi_2=\eta$ we obtain
\begin{equation}\label{e:rankcodazimu}
 \left[\tilde{\omega}_1(X) + \lambda^{-1}X(\lambda)\right]Y + \tilde{\omega}_2(X)D_\eta Y 
 = \left[\tilde{\omega}_1(Y) + \lambda^{-1}Y(\lambda)\right]X.
\end{equation}
On the other hand, bearing in mind that $A_\xi = A - \lambda I$, Eq. \eqref{e:codxi} yields
\begin{equation}\label{e:rankcodazzixi}
(Z\wedge W)\text{grad}\,\lambda =\lambda \tilde{\omega}_1(Z)W + \tilde{\psi}(Z)A_\eta W 
- \lambda\tilde{\omega}_1(W)Z - \tilde{\psi}(W)A_\eta Z.
\end{equation}
For $Z=X$ and $W=T \in \Gamma(\Delta)$, using part $(ii)$ of Lemma \ref{l:proptensorw} 
and $A_\eta X = 0=T(\lambda) = 0$, the preceding equation gives
$-X(\lambda)T = \lambda \tilde{\omega}_1(X)T$, 
hence $X(\lambda) = -\lambda\tilde{\omega}_1(X)$.  Substituting in  \eqref{e:rankcodazimu} yields
\begin{equation}\label{e:rankcodazimu2}
\tilde{\omega}_2(X)D_\eta Y = \left[\tilde{\omega}_1(Y) + \lambda^{-1}Y(\lambda)\right]X.
\end{equation}
 Eq. \eqref{e:rankcodazzixi} for $Z=T$ and $W=Y$ gives
$Y(\lambda)T = \tilde{\psi}(T)A_\eta Y - \lambda\tilde{\omega}_1(Y)T$,  
so $\Delta \leq \ker \tilde{\psi}$ and $-Y(\lambda) = \lambda\tilde{\omega}_1(Y)$. 
Therefore, taking into account that $A_\eta Y \neq 0$, substituting in  \eqref{e:rankcodazimu2} 
we obtain $\tilde{\omega}_2(X) = 0$.  Lastly, for $Z=X$ and $W=Y$, 
$$(X\wedge Y)\text{grad}\,\lambda =\lambda \tilde{\omega}_1(X)Y 
+ \tilde{\psi}(X)A_\eta Y - \lambda\tilde{\omega}_1(Y)X,
 $$
thus $\tilde{\psi}(X) = 0$.  In summary, we have 
\begin{equation}\label{e:rankcondition1}
\Delta \oplus \text{span}\{X\} \leq \ker \tilde{\psi} \cap \ker \tilde{\omega}_2
\end{equation}
and 
\begin{equation}\label{e:rankcondition2}
\lambda^{-1}Z(\lambda) + \tilde{\omega}_1(Z) = 0,
\end{equation}
for $Z \in \mathfrak{X}(M)$.  Using  \eqref{e:norcovderxi}, \eqref{e:norcovdermu}
 and \eqref{e:rankcondition2} we obtain 
\begin{equation}\label{e:parallelinducedtau1}
\nabla_Z^\perp (\mu + \xi - \zeta)
= \tilde{\psi}(Z)\eta,
\end{equation}
for $Z \in \mathfrak{X}(M)$.  Similarly, using \eqref{e:norcovderxi}, \eqref{e:norcovdermu} 
and \eqref{e:rankcondition2} we get  
\begin{equation}\label{e:parallelinducedtau2}
\nabla_Z^\perp \lambda(\xi - \zeta)
= \lambda\left(\tilde{\psi}(Z) + \tilde{\omega}_2(Z)\right)\eta.
\end{equation}
The second fundamental form of $\tilde{F}$ can be rewritten as 
$$
\alpha^{\tilde{F}}(X,Y) 
 = \left<AX,Y\right>(\mu + \xi - \zeta) + \left<A_\eta X,Y\right>\eta
  - \lambda\left<X,Y\right>(\xi - \zeta).
$$
Let $L=\mbox{span}\{\eta\}^\perp$ and let $F$ be the isometric light-cone 
representative of $f$. 
Define a vector bundle isometry $\tau \colon N_FM\to L$
by setting
$$\tau \Psi_*N=\mu + \xi - \zeta,\;\;\;\tau w
=\lambda(\xi-\zeta)\;\;\;\mbox{and}\;\;\;\tau F=\tilde F.$$
From  \eqref{e:parallelinducedtau1} and \eqref{e:parallelinducedtau2}, 
the vector bundle isometry $\tau$ is parallel with respect to the induced connection on $L$, 
and all the conditions of Lemma \ref{le:blem} are satisfied.  
As in the proof of Lemma \ref{l:ruledparabolic}, it follows from Lemma \ref{le:blem} 
and Proposition \ref{p:conformalchar} that $\tilde f$ is not a genuine 
conformal deformation of $f$, a contradiction.
 \vspace{1ex}

We now prove the converse. Start by choosing an orthonormal frame 
$\{\mu, \xi_1, \xi_2, \zeta\}$ of the trivial bundle $E = M^n \times \mathbb{L}^4$, 
with  $\zeta$ time-like. Extend  the tensors $D_i$ to $\Delta$ by requiring that 
$\Delta \leq \ker D_i$.  Define a compatible connection $\hat{\nabla}$ on $E$ by declaring
\begin{align}\label{eq:hatnabla}
\hat{\nabla}_X \mu &=-\omega_1(X)\xi_1 - \omega_2(X)\xi_2 
-\lambda^{-1}X(\lambda)\zeta=\nabla_X^\perp \zeta -\lambda^{1}X(\lambda)(\zeta-\mu),\nonumber\\
\hat{\nabla}_X \xi_1 &= \omega_1(X) (\mu - \zeta)+ \psi(X)\xi_2,\\
\hat{\nabla}_X \xi_2 &= \omega_2(X)(\mu - \zeta) - \psi(X)\xi_1,\nonumber
\end{align}
where 
\begin{equation}\label{eq:omegai}
\omega_i(X) = -\frac{1}{\lambda}\left<D_i X,\text{grad}\, \lambda\right>.
\end{equation} 
In particular, since $T(\lambda)=0$ for all $T\in \Delta$, 
by condition $(i)$ the sections $\mu$, $\xi_1$, $\xi_2$ and $\zeta$ are parallel  
along $\Delta$ with respect to  $\hat{\nabla}$.  

Let $\hat{\alpha}\colon \mathfrak{X}(M) \times \mathfrak{X}(M) \to \Gamma(E)$ 
be the bilinear map defined by 
\begin{align*}
\hat{\alpha}(X,Y) &= \left<AX,Y\right>\mu + \left<(A-\lambda I)D_1X,Y\right>\xi_1 
+ \left<(A-\lambda I)D_2X,Y\right>\xi_2\\
&\quad - \left<(A-\lambda I)X,Y\right>\zeta.
\end{align*}
From the symmetry of $(A-\lambda I)C_T$ (see  \eqref{e:covderhorA}), 
and because $C(\Gamma(\Delta))\subset \text{span}\{I,J\}$ and 
$C(\Gamma(\Delta))\not\subset \text{span}\{I\}$, for $f$ is not conformally surface-like,  
$(A-\lambda I)J$ is  symmetric.  Since $D_i \in \text{span}\{I,J\}$, 
also $(A-\lambda I)D_i$ is symmetric.  Thus  $\hat{\alpha}$ is symmetric.  

  We shall prove that $\hat{\alpha}$ satisfies the Gauss, Codazzi and Ricci equations 
  for an isometric immersion $\tilde F\colon M^n\to \mathbb{L}^{n+4}$. 
  For the Gauss equation, in view of the Gauss equation for $f$, 
  it  is enough to show that the bilinear form 
  $$\gamma(X,Y)=\left<(A-\lambda I)D_1X,Y\right>\xi_1 
  + \left<(A-\lambda I)D_2X,Y\right>\xi_2- \left<(A-\lambda I)X,Y\right>\zeta$$
is flat. Since 
\begin{equation}\label{e:proofdelta}
\Delta = \ker (A-\lambda I)D_1 \cap \ker (A - \lambda I)D_2 \cap \ker (A-\lambda I)=\ker \gamma,
\end{equation}
this is equivalent to
$ \det (A-\lambda I)D_1 + \det (A-\lambda I)D_2 - \det (A-\lambda I)=0,$
which holds in view of condition $(ii)$.  

To show that $\hat \alpha$ satisfies the Codazzi equations, we must prove that 
$$A_\mu = A, \quad {A_{\xi_1} = (A-\lambda I)D_1}, \quad A_{\xi_2} 
= (A-\lambda I)D_2 \quad \text{and} \quad A_\zeta = A - \lambda I.$$ 
satisfy the Codazzi equations.  The Codazzi equation for $A_\mu = A$ is equivalent to
\begin{equation}\label{eq:codmu} 
A_{\hat{\nabla}_Z \mu}W - A_{\hat{\nabla}_W \mu}Z=0
\end{equation}
for all $Z$, $W \in \mathfrak{X}(M)$.   For $W=T\in \Gamma(\Delta))$ and $Z \in \mathfrak{X}(M)$, 
this follows from  \eqref{e:proofdelta}
and the fact that $T(\lambda)=0$.  On the other hand,  by \eqref{eq:hatnabla}, \eqref{eq:omegai} 
and  item (ii),  for $Z=X$ and  $W=Y \in \Gamma(\Delta^\perp)$ 
the left-hand-side of \eqref{eq:codmu} is 
$$ \begin{array}{l}
 \lambda^{-1}(A-\lambda I)\left(-(D_1X \wedge D_1Y)\text{grad}\,\lambda 
-(D_2X\wedge D_2Y)\text{grad}\,\lambda\right)\vspace{.5ex}\\
\quad +  \lambda^{-1}(A-\lambda I)(X \wedge Y)\text{grad}\,\lambda = 0.
\end{array}
$$

Let us prove  the Codazzi equation of $A_\zeta = A - \lambda I$.  Using \eqref{eq:hatnabla} 
and the Codazzi equation for $A$,  taking into account that $\zeta$ is parallel along 
$\Delta$ and  that $T(\lambda)=0$ for $T\in \Gamma(\Delta)$, we obtain
\begin{align*}
(\nabla_Z A_\zeta)T &- (\nabla_T A_\zeta)Z - A_{\hat{\nabla}_Z \zeta}T + A_{\hat{\nabla}_T \zeta}Z \\
&= -Z(\lambda)T + T(\lambda)Z + \lambda^{-1}Z(\lambda)AT + \omega_1(Z)A_{\xi_1}T + \omega_2(Z)A_{\xi_2}T \\
&= 0,
\end{align*}
for all  $Z \in \mathfrak{X}(M)$.  For  $X, Y \in \Gamma(\Delta^\perp)$, using item (ii), 
\eqref{eq:hatnabla} and \eqref{eq:omegai}  we obtain
\begin{align*}
(\nabla_X A_\zeta)Y &- (\nabla_Y A_\zeta)X - A_{\hat{\nabla}_X \zeta}Y + A_{\hat{\nabla}_Y \zeta}X \\
&= -X(\lambda)Y + Y(\lambda)X + \lambda^{-1}X(\lambda)AY + \omega_1(X)A_{\xi_1}Y + \omega_2(X)A_{\xi_2}Y \\
&\quad - \lambda^{-1}Y(\lambda)AX - \omega_1(Y)A_{\xi_1}X - \omega_2(Y)A_{\xi_2}X\\
&= \lambda^{-1}(A-\lambda I)\left( -(X\wedge Y)\text{grad}\,\lambda \right)\\
&\quad + \lambda^{-1}(A-\lambda I)\left((D_1X\wedge D_1Y)\text{grad}\,\lambda 
+ (D_2X\wedge D_2Y)\text{grad}\,\lambda \right)\\
&= 0.
\end{align*}

Now we prove the Codazzi equation 
\begin{equation} \label{eq:codaxii}
(\nabla_Z A_{\xi_i})W-(\nabla_W A_{\xi_i})Z = A_{\hat{\nabla}_Z\xi_i}W - A_{\hat{\nabla}_W\xi_i}Z
\end{equation}
for $A_{\xi_i}=(A-\lambda I)D_i$.  First, let us suppose $Z=T$, $W=S \in \Gamma(\Delta)$.  
Then, because $\xi_i$ is parallel along $\Delta$, the right hand side of the equation is zero.  
Since $\Delta \leq \ker A_{\xi_i}$, we must show that
$$A_{\xi_i}\nabla_S T - A_{\xi_i}\nabla_T S = 0,$$
which follows easily using that $\Delta$ is an umbilical distribution.

Now, suppose $Z=X \in \Gamma(\Delta^\perp)$ and $W=T \in \Gamma(\Delta)$. 
By \eqref{eq:hatnabla} and  the fact that $\xi_i$ is parallel along $\Delta$,
 we get
\begin{align*}
(\nabla_X & A_{\xi_i})T-(\nabla_T A_{\xi_i})X - A_{\hat{\nabla}_X\xi_i}T + A_{\hat{\nabla}_T\xi_i}X\\
&= -(A-\lambda I)D_i\nabla_X T - \nabla_T (A-\lambda I)D_i X 
+ (A-\lambda I)D_i\nabla_T X - \lambda\omega_i(X)T.
\end{align*}
Taking the inner product with $S \in \Gamma(\Delta)$, using \eqref{eq:omegai}
and the fact that $\Delta$ is an umbilical distribution whose mean curvature
vector field $\delta$ satisfies $(A-\lambda I)\delta=- \text{grad}\,\lambda$, we get
$$
\left<(A-\lambda I)D_iX,\nabla_T S\right> - \lambda \omega_i(X)\left<T,S\right>=0.
$$\\
Equality between  the horizontal components follows from
\begin{align*}
\nabla_T^h(A-\lambda I)D_i 
&= (\nabla_T^hA)D_i\\
&= (A-\lambda I)D_i C_T
\end{align*}
where we have used  \eqref{e:covderhorA} and item $(iii)$.
The last case is when $X$, $Y \in \Gamma(\Delta^\perp)$.  
We have that $A_{\hat{\nabla}_X\xi_i}Y - A_{\hat{\nabla}_Y\xi_i}X$ coincides with
$$\begin{array}{l}
\omega_i(X)(AY - A_\zeta Y) + (-1)^j\psi(X)A_{\xi_j}Y -\omega_i(Y)(AX -A_\zeta X) - (-1)^j\psi(Y)A_{\xi_j}X\\
= \lambda \omega_i(X)Y - \lambda \omega_i(Y)X  + (-1)^j(A-\lambda I)\left(\psi(X)D_jY - \psi(Y)D_jX\right)\\
=-D_iX(\lambda)Y + D_iY(\lambda)X + (-1)^j(A-\lambda I)\left(\psi(X)D_jY - \psi(Y)D_jX\right)\\
= (X\wedge Y)D_i^t \text{grad}\,\lambda + (-1)^j(A-\lambda I)\left(\psi(X)D_jY - \psi(Y)D_jX\right),
\end{array} $$
and this is equal to the left-hand-side of \eqref{eq:codaxii} by item (iv).

Now, let us move on to the Ricci equations.  It is easily checked using \eqref{eq:hatnabla} 
that $\big<\hat{R}(Z,W)\mu,\zeta\big>=0$, hence
the Ricci equation for $\mu$ and $\zeta$ is satisfied because $A_\mu = A$ and 
$A_\zeta = (A - \lambda I)$ commute. It is also easily seen that  the Ricci equation for 
$\zeta$ and ${\xi_i}$ is equivalent to that for $\mu$ and ${\xi_i}$. 

Let us prove the Ricci equation for $\mu$ and ${\xi_i}$.    
First, let us prove for $X$, $Y \in \Gamma(\Delta^\perp)$.  
On one hand, by the symmetry of $A$ and $(A-\lambda I)D_i$  we have 
$$\left<[A_{\xi_i},A_\mu]X,Y\right> = \left<AX,(A-\lambda I)D_iY\right> - \left<(A-\lambda I)D_iX,AY\right>.$$
On the other hand, a straightforward computation using \eqref{eq:hatnabla} and \eqref{eq:omegai} gives
\begin{align*}
\big<\hat{R}(X,Y)\xi_i,\mu\big> &= \lambda^{-1}\left(\left<(\nabla_Y D_i)X - 
(\nabla_X D_i)Y,\text{grad}\,\lambda\right>\right)\\
&\quad  + \lambda^{-1}\left(\text{Hess}\,\lambda(Y,D_iX)-\text{Hess}\,\lambda(X,D_iY)\right) \\
&\quad + \lambda^{-1}\left((-1)^j\psi(X)\left<D_jY,\text{grad}\,\lambda\right>- 
(-1)^j\psi(Y)\left<D_jX,\text{grad}\,\lambda\right> \right).
\end{align*}
Thus  the Ricci equation for $\xi_i$ and $\mu$ for $X$, 
$Y \in \Gamma(\Delta^\perp)$ follows from item (v).

Now for $X \in \Gamma(\Delta^\perp)$ and $T \in \Gamma(\Delta)$, we have, on one hand,
$$\left<[A_{\xi_i},A_\mu]X,T\right> = 0 $$
while, on the other hand, 
\begin{align*}
\big<\hat{R}(X,T)\xi_i,\mu\big> &=-T\omega_i(X) - \omega_i([X,T]) \\
&= T\left(\frac{1}{\lambda}\left<D_iX,\text{grad}\,\lambda\right>\right) 
+ \frac{1}{\lambda}\left<D_i[X,T],\text{grad}\,\lambda\right>\\
&= \frac{1}{\lambda}\left(D_iXT(\lambda) - \left<\nabla_{D_iX}T,\text{grad}\,\lambda\right> 
+ \left<D_i \nabla_X T,\text{grad}\,\lambda\right>\right)\\
&=\frac{1}{\lambda}\left( \left<C_T D_i X,\text{grad}\,\lambda\right> 
- \left<D_iC_TX,\text{grad}\,\lambda\right>\right) = 0,
\end{align*}
where we have used both equalities in item (iii).
Lastly, for $T$ and $S \in \Gamma(\Delta)$,  on one hand, $\left<[A_{\xi_i},A_\mu]T,S\right> = 0$
 because $\ker A_{\xi_i} = \Delta$.  On the other hand, $\big<\hat{R}(T,S)\xi_i,\mu\big>=0$
because $\xi_i$ is parallel along $\Delta$ and $[T,S] \in \Gamma(\Delta)$.

It remains to verify the Ricci equation for ${\xi_1}$ and ${\xi_2}$.  
From \eqref{eq:hatnabla} we obtain
$$
\big<\hat{R}(Z,W)\xi_1,\xi_2\big> = \text{d}\psi([Z,W]).
$$
Thus the Ricci equation for ${\xi_1}$ and ${\xi_2}$  follows from item (vi) 
if either $Z$ or $W$ belongs to $\Gamma(\Delta)$, 
and from item (vii) if  both $Z$ and $W$ belong to $\Gamma(\Delta^\perp)$.

By the Fundamental Theorem of Submanifolds, there exist an isometric immersion 
$\tilde{F}\colon M^n \to \mathbb{L}^{n+4}$ and a vector bundle isometry $\Phi\colon  E \to N_{\tilde{F}}M$ such that
$$\Phi \circ \hat{\alpha}=\alpha^{\tilde{F}} \quad \text{and} \quad \Phi \hat{\nabla} = \nabla^\perp \Phi.$$
  Moreover, the vector field $\rho = \lambda^{-1}\Phi(\zeta - \mu)$ satisfies
\begin{align*}
\lambda\tilde{\nabla}_X\rho	 &= \lambda X(\lambda^{-1})\Phi(\zeta - \mu) + \tilde{\nabla}_X\Phi(\zeta - \mu)\\
&= - \lambda^{-1}X(\lambda)\Phi(\zeta-\mu) -\tilde{F}_*A_{\Phi(\zeta-\mu)}X + \nabla^\perp_X\Phi(\zeta-\mu)\\
&= - X(\lambda)\rho - \tilde{F}_*A_{\zeta-\mu}X + \Phi\hat{\nabla}_X(\zeta-\mu)\\
&=- X(\lambda)\rho + \lambda \tilde{F}_*X + \lambda^{-1}X(\lambda)\Phi(\zeta-\mu)\\
&= \lambda \tilde F_*X
\end{align*}
for all $X \in\mathfrak{X}(M)$.  Therefore
$\tilde{F} - \rho$ is a constant vector $P_0 \in \mathbb{L}^{n+4}$, with
$$\big<\tilde{F}-P_0,\tilde{F}-P_0\big> = \big<\rho,\rho\big> 
= \lambda^{-2}\big<\zeta-\mu,\zeta-\mu\big>= 0,$$
that is, $\tilde{F}$ takes values in $P_0 + \mathbb{V}^{n+3}$.  
Without loss of generality, suppose $P_0=0$, otherwise redefine $\tilde{F}$ 
by $\tilde{F} - P_0$. Then, $\tilde{F}$ gives rise to a conformal immersion 
$\tilde{f}=\mathcal{C}(\tilde{F})\colon M^n \to \mathbb{R}^{n+2}$.

  We now prove that  $\tilde{f}$ is a genuine deformation of $f$. 
Assume otherwise.  By Proposition \ref{p:conformalchar}, there exist an open set 
$U \subset M^n$ and an isometric immersion $H\colon W \to \mathbb{V}^{n+3}$, 
with $W \supset F(U)$ open in $\mathbb{V}^{n+2}$, such that $\tilde{F}|_U = H \circ F|_U$.  
For simplicity, we will suppose $U = M^n$.  Because $f$ is an isometric immersion, 
its  isometric light-cone representative  is  $F = \Psi \circ f$.  
We conclude that $\tilde{F} = T \circ f$ for  
$T=H\circ \Psi\colon V\subset \mathbb{R}^{n+1} \to \mathbb{V}^{n+3} \subset \mathbb{L}^{n+4}$.  

Since $T$ is an isometric immersion into the light-cone, the position vector field  
$T$ is a section of its normal bundle $N_T\mathbb{R}^{n+1}$ such that
\begin{equation}\label{eq:sffpos}
\left<\alpha^T(Z,W),T\right> = - \left<Z,W\right>
\end{equation}
 for all $Z$, $W \in \mathfrak{X}(\mathbb{R}^{n+1})$.  
 Complete $T$ to a pseudo-orthonormal frame $\{\rho, T, \tilde{\zeta}\}$ 
 of $\Gamma(N_T\mathbb{R}^{n+1})$, where $\tilde{\zeta}$ is a light-like vector field 
 such that $\left<\tilde{\zeta},T\right>=1$.  
 We can associate to this frame the orthonormal frame given by 
 $\{ \rho, (T + \tilde{\zeta})/\sqrt{2}, (T - \tilde{\zeta})/\sqrt{2}\}$.

By the Gauss equation of $T$,  the bilinear form $\alpha^T$ is flat.  
It follows from \eqref{eq:sffpos} that $\mathcal{N}(\alpha^T) = \{0\}$,
hence    $\dim \Omega = \dim \left( \mathcal{S}(\alpha^T) \cap \mathcal{S}(\alpha^T)^\perp\right) = 1$
by Lemma \ref{l:mainlemma}. 
The projections $P_i$, $i=1$, $2$, of $N_T\mathbb{R}^{n+1}$ onto the subspaces 
$$W_1 = \text{span}\left\{ \rho, \frac{T + \tilde{\zeta}}{\sqrt{2}}\right\} 
\quad \text{and} \quad W_2=\text{span}\left\{\frac{T - \tilde{\zeta}}{\sqrt{2}}\right\}$$
map $\Omega$  isomorphically onto their images.  
By dimensional reasons, $P_2|_\Omega$ is an isomorphism.  
Let $\beta \in \Omega$ be such that $P_2(\beta) = (T-\tilde{\zeta})/\sqrt{2}.$
Then $\beta$ is a light-like vector field with $A_\beta^T = 0$, and we can write
$$\beta = \cos \theta \rho + \sin \theta \frac{T+\tilde{\zeta}}{\sqrt{2}} + \frac{T-\tilde{\zeta}}{\sqrt{2}},$$
where $\theta \in [0, 2\pi)$.  Define $\{\gamma, \delta, \tilde{\gamma}\}$ by
$$\gamma =  \cos \theta \rho + \sin \theta \frac{T+\tilde{\zeta}}{\sqrt{2}}, \quad \delta 
=  -\sin \theta \rho + \cos \theta \frac{T+\tilde{\zeta}}{\sqrt{2}} 
\quad \text{and} \quad \tilde{\gamma} = \frac{T-\tilde{\zeta}}{\sqrt{2}}.$$
Since $\beta = \gamma + \tilde{\gamma}$ and $A_\beta^T = 0$, then $A_\gamma^T = -A_{\tilde{\gamma}}^T$.  
Moreover, because
\begin{align*}
\alpha^T(Z,W) &= \left<A_\delta^T Z,W\right>\delta + \left<A_\gamma^T Z,W\right>\gamma 
- \left<A_{\tilde{\gamma}}^TZ,W\right>\tilde{\gamma}\\
&= 	\left<A_\delta^T Z,W\right>\delta + \left<A_\gamma^T Z,W\right>\beta
\end{align*}
for all $Z$, $W \in \mathfrak{X}(\mathbb{R}^{n+1})$, we conclude from the flatness of $\alpha^T$, 
and the fact that $\beta$ is light-like and orthogonal to $\delta$, 
that $\text{rank}\, A_\delta^T \leq 1$.  Therefore, 
\begin{equation}\label{e:propertiesinherited}
A^{\tilde{F}}_{T_*N} =  A, \quad A^{\tilde{F}}_{\gamma \circ f} = - A_{\tilde{\gamma} \circ f}^{\tilde{F}}X 
\quad \text{and} \quad \text{rank}\,A_{\delta \circ f}^{\tilde{F}} \leq 1.
\end{equation}
Notice that, since
$T = \frac{\sqrt{2}}{2}\left(\cos \theta \delta + \sin \theta \gamma + \tilde{\gamma}\right)$ 
and $\tilde{F} = T \circ f$, then
\begin{equation}\label{e:Fequation}
\tilde{F} = \frac{\sqrt{2}}{2}\left(\cos \theta (\delta \circ f) 
+ \sin \theta (\gamma \circ f) + (\tilde{\gamma} \circ f)\right).
\end{equation}
On the other hand, by \eqref{e:formatF} we have  
$\alpha_{\tilde F}(S,S)=\lambda \mu$
for all $S\in \Gamma(\Delta) = \ker (A-\lambda I)$. Comparing with
\begin{align*}
\alpha^{\tilde{F}}(S,S) &= \lambda T_*N + \left<A_\gamma^T f_*S, f_*S\right>\left((\gamma \circ f) 
+ (\tilde{\gamma}\circ f)\right)\\
&=\lambda T_*N - \frac{\sqrt{2}}{\sin\theta -1}	\left((\gamma \circ f) + (\tilde{\gamma}\circ f)\right).
\end{align*}
 we obtain
\begin{equation}\label{e:mueq}
\mu = T_*N - \frac{\sqrt{2}}{\lambda(\sin\theta -1)}	\left((\gamma \circ f) 
+ (\tilde{\gamma}\circ f)\right).
\end{equation} 
It is now straightforward to verify  that 
\begin{align*}
\xi_1 &= T_*N + \left( \frac{\lambda \cos^2 \theta}{\sqrt{2}(1-\sin \theta)} 
- \frac{\lambda}{\sqrt{2}} \right)(\gamma \circ f) 
+ \left( \frac{\lambda \cos^2 \theta}{\sqrt{2}(1-\sin \theta)} 
- \frac{\lambda \sin \theta}{\sqrt{2}} \right) (\tilde{\gamma} \circ f)\\
&\quad + \frac{\lambda \cos \theta}{\sqrt{2}}(\delta \circ f)\quad \text{and}\quad 
\xi_2 = \frac{\cos \theta}{1 - \sin \theta}\left( (\gamma \circ f)
+ (\tilde{\gamma}\circ f) \right) + (\delta \circ f).
\end{align*}
is an orthonormal frame for $\mathbb{P} = \{\mu, \zeta\}^\perp$.  From 
\eqref{e:propertiesinherited} we have 
\begin{equation}\label{eq:axi1tilF}	
A_{\xi_1}^{\tilde{F}} = A +\frac{\lambda}{\sqrt{2}}\left((\sin \theta - 1)A_{\gamma \circ f}^{\tilde{F}}
 + \cos \theta A_{\delta \circ f}^{\tilde{F}}\right)
\;\;\;\mbox{and}\;\;\;
A_{\xi_2}^{\tilde{F}} = A_{\delta \circ f}^{\tilde{F}}.
\end{equation}
The last relation in \eqref{e:propertiesinherited} implies that the rank of $D_{\xi_2} 
= (A-\lambda I)A_{\xi_2}$ is less than or equal to one.  
Now, by \eqref{e:Fequation} and the second  relation in \eqref{e:propertiesinherited} we have 
\begin{equation}
-I=A_{\tilde F}^{\tilde{F}}
=\frac{\sqrt{2}}{2}(\cos \theta A_{\delta \circ f}^{\tilde{F}}+(\sin \theta -1) A_{\gamma \circ f}^{\tilde{F}}).
\end{equation}
Substituting this expression in the first equation of \eqref{eq:axi1tilF} implies that
 $A^{\tilde{F}}_{\xi_1} = A - \lambda I$,  and hence $D_{\xi_1} = I$.   
Let $\theta \in [0,\pi/2]$ be such that
$$D_1 = \cos \theta D_{\xi_1} + \sin \theta D_{\xi_2}\;\;\;\mbox{and}\;\;\; D_2 
= - \sin \theta D_{\xi_1} + \cos \theta D_{\xi_2},$$
where $D_1$ and $D_2$ have  determinant $1/2$.  Then $D_1^2 + D_2^2 - I = D_{\xi_2}^2$,
and this means that $\text{rank }D_1^2 + D_2^2 - I < 2$, a contradiction with (ix).  

It remains to prove the last statement of Proposition \ref{p:equivalentcartanhyp}.  
First, suppose that the triples $(D_1,D_2,\psi)$ and $(\hat{D}_1,\hat{D}_2,\hat{\psi})$ 
give rise to congruent conformal immersions $\tilde{f}$ and $\tilde{g}$.  
Then, by Proposition \ref{p:conformalchar}, their isometric light-cone representatives 
$\tilde{F}$ and $\tilde{G}$ are congruent isometric immersions, that is, there exists 
$T \in O_1^+(m+4)$ such that $\tilde{G} = T \circ \tilde{F}$.  
Hence, $\alpha^{\tilde{G}} = T \circ \alpha^{\tilde{F}}$ and  
$\hat{\nabla}^\perp T = T\nabla^\perp$. 
From the equality regarding second fundamental forms applied to $(T,T) \in \Delta \times \Delta$ 
we conclude that $T(\mu) = \hat{\mu}$.  Taking into account the last fact, 
from the equality $\tilde{G} = T \circ \tilde{F}$ we get $T(\zeta) = \hat{\zeta}$.  
Now, from
$$\big<A_{T(\xi_i)}^{\tilde{G}}X,Y\big> = \big<\alpha^{\tilde{G}}(X,Y),T(\xi_i)\big> 
=\big<\alpha^{\tilde{F}}(X,Y),\xi_i\big> = \big<A_{\xi_i}^{\tilde{F}}X,Y\big>$$
and the uniqueness of the sections $\hat{\xi_i}$ such that $\det D_{\hat{\xi_i}} = 1/2$, 
we conclude that $T(\xi_i) = \hat{\xi_i}$ and $D_i = \hat{D}_i$.  
From $\hat{\nabla}^\perp T = T\nabla^\perp$ we obtain that $\psi=\hat{\psi}$.

For the converse, suppose the conformal immersions $\tilde{f}$ and $\tilde{g}$ 
have the same triples.  By the uniqueness of the frame $\{\xi_1, \xi_2\}$, 
we can define $T: N_{\tilde{F}}M \to N_{\tilde{G}}M$ by  
$T(\mu) = \hat{\mu}$, $T(\xi_i) = \hat{\xi}_i$ and $T(\zeta) = \hat{\zeta}$.  
Since  the triples are the same, we have 
$\hat{\nabla}^\perp T = T\nabla^\perp$ and $\alpha^{\tilde{G}} 
= T \circ \alpha^{\tilde{F}}$, hence  $\tilde F$ and $\tilde G$  are congruent.  

\section{The Reduction}
\label{ch:reduction}

In this section, for a hypersurface  $f:M^n \to \mathbb{R}^{n+1}$   
that is not conformally surface-like  and envelops a two-parameter 
congruence of hyperspheres $s:L^2 \to \mathbb{S}^{n+2}_{1,1}$, 
the problem of finding a pair of tensors $(D_1, D_2)$  and a 
one-form $\psi$ on $M^n$ satisfying all the conditions in 
Proposition \ref{p:equivalentcartanhyp} is reduced to a 
similar but easier one on the surface $s$.  
First we give a few definitions.

The surface $s\colon L^2 \to \mathbb{S}^{n+2}_{1,1}$ is said to be \emph{hyperbolic} 
(respectively, \emph{elliptic}) with respect to a tensor $\bar{J}$ on $L^2$ 
satisfying $\bar{J}^2 = \bar{I}$ (respectively, $\bar{J}^2 = -\bar{I}$) if
$$
\alpha'(\bar X, \bar J\bar Y)=\alpha'(\bar J\bar X, \bar Y)
$$
for all $\bar X, \bar Y\in \mathfrak{X}(L)$, where $\alpha'$ is the 
second fundamental form of $s$.

Now let $\pi:M \to L$ be a submersion.  A vector field $X \in \mathfrak{X}(M)$ is said to be 
\emph{projectable} if it is $\pi$-related to a vector field $\bar{X} \in \mathfrak{X}(L)$, 
that is, there exists $\bar{X} \in \mathfrak{X}(L)$ such that $\pi_*X=\bar X\circ \pi$.  
A tensor $D\in \Gamma(\text{End}\,(TM))$  is \emph{projectable}
if there exists $\bar{D}\in \Gamma(\text{End}\,(TL))$  such that 
$\bar{D} \circ \pi_* = \pi_* \circ D$. Similarly, a one-form $\omega$ on $M$ is 
\emph{projectable} if there exists a one-form $\bar{\omega}$ on $L$ such that 
$\bar{\omega}\circ \pi_* = \omega$.  

We will need the following result of \cite{mD2013}, which gives  conditions 
for  tensors and one-forms to be projectable.  


\begin{proposition}\label{c:projoneform}
Let $\Delta$ be an integrable distribution on a differentiable manifold $M$, 
let $L = M/\Delta$ be the (local) quotient space of leaves of $\Delta$ and 
let $\pi\colon M \to L$ be the quotient map.  Then the following assertions hold:
\begin{itemize}
\item[(i)] a one-form $\omega$ on 
$M$ is projectable if and only if $\omega(T) = 0$ and 
$\text{d}\omega(T,X) = 0$ for any $T \in \Gamma(\Delta)$ and 
$X \in \Gamma(\Delta^\perp)$; 	
\item[(ii)] if $M^n$ is a Riemannian manifold and 
$C\colon \Gamma(\Delta) \to \Gamma(\text{End}(\Delta^\perp))$
is the splitting tensor of $\Delta$, then 
$D\in \Gamma(\text{End}(\Delta^\perp))$ is projectable if and only if
$\nabla_T^h D = [D,C_T]$ 	for all $T \in \Gamma(\Delta)$.
	\end{itemize}
\end{proposition}

The reduction lemma is as follows. 

\begin{lemma}\label{l:reductionlemma}
Let $f:M^n \to \mathbb{R}^{n+1}$ be a hypersurface that is not conformally 
surface-like and envelops a two-parameter congruence of hyperspheres 
$s:L^2 \to \mathbb{S}^{n+2}_{1,1}$.
Let $\Delta$ be the eigenbundle of $f$ correspondent to its principal curvature 
$\lambda$ of multiplicity $n-2$.  If $f$ is hyperbolic (respectively, elliptic)
with respect to $J \in \Gamma(\text{End}(\Delta^\perp))$ and there exists a triple
$(D_1,D_2,\psi)$  satisfying conditions (i)-(ix) 
in Proposition~\ref{p:equivalentcartanhyp}, with $D_i \in \text{span}\{I,J\}$ 
for $i=1,2$ and $\psi$ a one-form on $M^n$, then $J$, $D_1$ and $D_2$ are the horizontal 
lifts of tensors $\bar{J}$, $\bar{D}_1$, $\bar{D}_2 \in \text{span}\{\bar{I},\bar{J}\}$ 
on $L^2$, with $\bar{J}^2 = I$ (respectively, $\bar{J}^2 = -I$) and $\psi$ is the horizontal 
lift of a one-form $\bar{\psi}$ on $L^2$ such that $s$ is hyperbolic (respectively, elliptic) 
with respect to $\bar{J}$ and the triple $(\bar{D}_1,\bar{D}_2,\bar{\psi})$ satisfies:
\begin{enumerate}[(a)]
\item $\det \bar{D}_i = 1/2$,
\item $(\nabla_X'\bar{D}_i)Y - (\nabla_Y'\bar{D}_i)X 
= (-1)^j\left((\bar{\psi}(X)\bar{D}_jY - \bar{\psi}(Y)\bar{D}_j(X)\right)$,
\item $\text{d}\bar{\psi}(X,Y) 
= \left<\bar{D}_2X,\bar{D}_1 Y\right>'-\left<\bar{D}_1X,\bar{D}_2 Y\right>'$,
\item $\bar{D}_2^2 \neq \pm \bar{D}_1^2$,
\item $\text{rank}\, (\bar{D}_1^2+\bar{D}_2^2-\bar{I}) =2$.
\end{enumerate}
Conversely, if $s\colon L^2 \to \mathbb{S}^{n+2}_{1,1}$ is hyperbolic (respectively, elliptic) 
with respect to a tensor $\bar{J}$ on $L^2$ satisfying $\bar{J}^2 = \bar{I}$ 
(respectively, $\bar{J}^2 = -\bar{I}$), then the hypersurface $f$ is hyperbolic 
(respectively, elliptic) with respect to the horizontal lift $J$ of $\bar{J}$, 
and the horizontal lifts $D_1$, $D_2$ and $\psi$ of 
 $\bar{D}_1$, $\bar{D}_2 \in \text{span}\{\bar{I},\bar{J}\}$ and the one-form 
 $\bar{\psi}$ satisfying items (a) to (e) have all the properties (i) to (ix) 
 in Proposition \ref{p:equivalentcartanhyp}.
\end{lemma}

\proof
Conditions (i) and (vi) of Proposition \ref{p:equivalentcartanhyp}, 
together with part $(i)$ of Proposition  \ref{c:projoneform}, 
assure us that the one-form $\psi$ is projectable with respect to the canonical 
projection $\pi: M \to L^2$  onto the (local) quotient of leaves of the distribution 
$\Delta$, that is, there exists a one-form $\bar{\psi}$ on $L^2$ such that
	$\bar{\psi}\circ \pi_* = \psi.$
	
	The tensors $D_1$ and $D_2$ are also projectable, because of item (iii) 
	of Proposition~\ref{p:equivalentcartanhyp} and part $(ii)$ of 
	Proposition \ref{c:projoneform}, that is, there exist tensors 
	 $\bar{D}_1$ and $\bar{D}_2$ on $L^2$ such that  
	\begin{equation}\label{eq:d1d2}
	\bar{D}_1 \circ \pi_* = \pi_* \circ D_1 \quad \text{and} \quad \bar{D}_2 \circ \pi_* 
	= \pi_* \circ D_2.
	\end{equation}
	
	From item (iii) we have that  $D_i$, $i=1,2$, commute with  $C_T$ for all  
	$T\in \Gamma(\Delta)$. Since  $D_i\in \text{span}\,\{I,J\}$, and taking into account 
	item (viii), at least one $D_i$ is of the form $D_i = a_i I + b_i J$ with $b_i\neq 0$. 
	It follows that  $[C_T,J]=0$.  The fact that ${f\colon M^n \to \mathbb{R}^{n+1}}$ 
	is hyperbolic or elliptic gives us that $\nabla_T^h J= 0$.  
	Therefore $J$ is  projectable, that is, there is $\bar{J} \in \text{End}(TL)$ such that
	$\pi_* \circ J = \bar{J} \circ \pi_*.$
		Since $D_i \in \text{span}\{I,J\}$, we get that 
		$\bar{D}_i \in \text{span}\{\bar{I},\bar{J}\}$ from \eqref{eq:d1d2}.  
		From $J^2 = \epsilon I$, where $\epsilon = 1$ or $\epsilon=-1$ according to whether 
		$f$ is hyperbolic or elliptic, it follows that $\bar{J}^{2} = \epsilon \bar{I}$ 
		and that $(a)$, $(d)$ and $(e)$ hold. 
	
	Let $S\colon M^n\to \mathbb{S}_{1,1}^{n+2}\subset \mathbb{L}^{n+3}$ be the two-parameter congruence 
	of hyperspheres enveloped by $f$, so that $S=s \circ \pi$. We have
	\begin{equation}\label{e:envelopeS}
	S(x) = \Psi_{*}(f(x))N(x) + \lambda(x)\Psi(f(x))
	\end{equation}
for all $x\in M^n$. 	Differentiating  \eqref{e:envelopeS} with respect to 
$Y \in \mathfrak{X}(M)$ gives
\begin{equation}\label{e:envelopeSdif1}
	S_*Y = -\Psi_*f_*(A - \lambda I)Y + Y(\lambda)\Psi \circ f.
	\end{equation}
	In particular, 
	\begin{equation}\label{eq:indmetric}
	\langle S_*X, S_*Y\rangle = \langle (A - \lambda I)X, (A - \lambda I)Y\rangle
	\end{equation}
	for all $X, Y\in \mathfrak{X}(M)$.		
	Replacing $Y$ by $D_iY$ in  	\eqref{e:envelopeSdif1} we get
	\begin{equation}\label{e:envelopeSdif2}
     \Psi_*f_*(A - \lambda I)D_iY=  \left<D_iY,\text{grad}\,\lambda\right>\Psi \circ f - S_*D_iY.
	\end{equation}
	Differentiating one more time  \eqref{e:envelopeSdif2} with respect to 
	$X \in \Gamma(\Delta^\perp)$ yields
	\begin{equation}\label{e:envelopeSdif3}
	\begin{array}{l}
	\tilde{\nabla}_X \Psi_*f_*(A - \lambda I)D_iY = 	 
	\left<\nabla_X D_iY,\text{grad}\,\lambda\right>\Psi \circ f 
	 + \text{Hess}\,\lambda (X,D_iY)\Psi \circ f \vspace{1ex}\\
	 \hspace*{25ex}+  \left<D_iY,\text{grad}\,\lambda\right>\Psi_*f_*X - \tilde{\nabla}_X S_{*}D_iY.
	 \end{array}
	 \end{equation}
	Let $\hat \nabla$ be the connection of $\mathbb{S}^{n+2}_{1,1}$, 
	 $\left<\cdot, \cdot\right>'$  be the metric on $L^2$ induced by $s$ 
	 and $\nabla'$ its Levi-Civita connection. Then 
	\begin{align}\label{e:envelopeSdif2exp1}
	\tilde{\nabla}_X S_{*}D_iY 
	&=\tilde{\nabla}_{\pi_*X}s_*\bar{D}_i\pi_*Y\\
	&= \hat{\nabla}_{\pi_*X}s_*\bar{D}_i\pi_*Y 
	- \left<\pi_*X, \bar{D}_i\pi_*Y\right>'s\circ\pi\nonumber\\
	&= s_* \nabla'_{\pi_*X}\bar{D}_i\pi_*Y + \alpha'(\pi_*X,\bar{D}_i\pi_*Y) 
	- \left<\pi_*X, \bar{D}_i\pi_*Y\right>'s\circ\pi\nonumber,
	\end{align}
	for all projectable vector fields $X$, $Y \in \Gamma(\Delta^\perp)$. 
	By  \eqref{eq:indmetric} we have
\begin{align}\label{e:envelopeSdif2exp2}
	\left<\pi_*X, \bar{D}_i\pi_*Y\right>' &= \left<s_*\pi_*X, s_*\bar{D}_i\pi_*Y\right>	\\
	&=\left<(A-\lambda I)X,(A-\lambda I)D_iY\right> \nonumber.
	\end{align}
Therefore, substituting  \eqref{e:envelopeSdif2exp1} and \eqref{e:envelopeSdif2exp2} 
in  \eqref{e:envelopeSdif3} we obtain
\begin{equation} \label{eq:firstdif}
\begin{array}{l}
\tilde{\nabla}_X \Psi_*f_*(A - \lambda I)D_iY= 
\left<\nabla_X D_iY,\text{grad}\,\lambda\right>\Psi \circ f 
+ \text{Hess}\,\lambda (X,D_iY)\Psi \circ f +\vspace{1ex}\\
\hspace*{5ex}+ \left<D_iY,\text{grad}\,\lambda\right>\Psi_*f_*X 
	 - s_* \nabla'_{\pi_*X}\bar{D}_i\pi_*Y - \alpha'(\pi_*X,\bar{D}_i\pi_*Y) + \vspace{1ex}\\
\hspace*{5ex}+ \left<(A-\lambda I)X,(A-\lambda I)D_iY\right>(\Psi_*N + \lambda(\Psi \circ f)).
\end{array}
\end{equation}

On the other hand, from  \eqref{eq:sffpsi} and  \eqref{e:envelopeSdif2} we get
\begin{align}\label{e:envelopevital2}
&\tilde{\nabla}_X \Psi_*f_*(A - \lambda I)D_iY= \\
&=	 \Psi_* \bar{\nabla}_X f_*(A-\lambda I)D_iY  + 
\alpha^\Psi(f_*X, f_*(A-\lambda I)D_iY) \nonumber\\
&=\Psi_* f_*\nabla_X (A-\lambda I)D_iY + \left<AX,(A-\lambda I)D_iY\right>\Psi_*N 
- \left<X,(A-\lambda I)D_iY\right>w\nonumber\\
&= \Psi_* f_*(\nabla_X (A-\lambda I)D_i)Y + \Psi_* f_* (A-\lambda I)D_i\nabla_X Y 
+ \left<AX,(A-\lambda I)D_iY\right>\Psi_*N\nonumber\\
&\quad   - \left<X,(A-\lambda I)D_iY\right>w\nonumber\\
&= \Psi_* f_*(\nabla_X (A-\lambda I)D_i)Y 
+ \left<D_i\nabla_XY,\text{grad}\,\lambda\right>\Psi \circ f 
- s_{*}\bar{D}_i\pi_*\nabla_{X}Y \nonumber\\
&\quad + \left<AX,(A-\lambda I)D_iY\right>\Psi_*N  
- \left<X,(A-\lambda I)D_iY\right>w.\nonumber
\end{align}
Computing $\tilde{\nabla}_X \Psi_*f_*(A - \lambda I)D_iY-\tilde{\nabla}_Y \Psi_*f_*(A - \lambda I)D_iX$, 
first using \eqref{eq:firstdif} and then
\eqref{e:envelopevital2}, and comparing both expressions give
\begin{align}\label{e:equivalent}
\Psi_*&f_*B(X,Y) + \theta(X,Y)\Psi_*N + \varphi(X,Y)\Psi \circ f -\lambda^{-1}\theta(X,Y)w	\\
&\quad = s_*( (\nabla_{\pi_*Y}'\bar{D}_i)\pi_*X 
- (\nabla_{\pi_*X}'\bar{D}_i)\pi_*Y ) 
+ \alpha'(\pi_*Y,\bar{D}_i\pi_*X) - \alpha'(\pi_*X,\bar{D}_i\pi_*Y)\nonumber
\end{align}
where
$$B(X,Y) = (\nabla_X(A-\lambda I)D_i)Y - (\nabla_Y (A-\lambda I)D_i)X 
- X\wedge Y(D_i^t\text{grad}\,\lambda),$$
$$\theta(X,Y) = \lambda( \left<X,(A-\lambda I)D_iY\right> 
- \left<Y,(A-\lambda I)D_iX\right>),$$
\begin{align*}
\varphi(X,Y) &= \left<(\nabla_YD_i)X - (\nabla_XD_i)Y	,\text{grad}\,\lambda\right> 
+ \text{Hess}\,\lambda(D_iX,Y) - \text{Hess}\,\lambda(X,D_iY)\\
&\quad - \lambda( \left<(A-\lambda I)X,(A-\lambda I)D_iY\right> 
- \left<(A-\lambda I)D_iX,(A-\lambda I)Y\right>),
\end{align*}
for all projectable $X$, $Y \in \Gamma(\Delta^\perp)$.  
Since $(A-\lambda I)C_T$ is symmetric by  \eqref{e:covderhorA}, 
and $\rm{span}\{I\} < C(\Gamma(\Delta)) \leq \rm{span}\{I,J\} $ because $f$ 
is either hyperbolic or elliptic and not surface-like, we have that 
$(A-\lambda I)J$ is  symmetric.  Thus $(A-\lambda I)D_i$ is symmetric 
for $i=1, 2$, for $D_i \in \text{span}\,\{I,J\}$.   
Using this,  \eqref{e:equivalent} and  items (iv) and (v) 
of Proposition~\ref{p:equivalentcartanhyp}  we obtain
\begin{align}\label{e:equivalent2}
(-1)^j&\Psi_*f_*(A-\lambda I)(\psi(X)D_jY - \psi(Y)D_jX) \\
&\quad + ((-1)^j\psi(Y)\left<D_jX,\text{grad}\,\lambda\right> 
- (-1)^j\psi(X)\left<D_jY,\text{grad}\,\lambda\right>)\Psi \circ f\nonumber\\
&= s_*( (\nabla_{\pi_*Y}'\bar{D}_i)\pi_*X - (\nabla_{\pi_*X}'\bar{D}_i)\pi_*Y ) 
+ \alpha'(\pi_*Y,\bar{D}_i\pi_*X) - \alpha'(\pi_*X,\bar{D}_i\pi_*Y)\nonumber.
\end{align}
Using   \eqref{e:envelopeSdif1} we get 
\begin{align*}
(-1)^j&\bar{\psi}(\pi_*Y)s_*\bar{D}_j\pi_*X -(-1)^j\bar{\psi}(\pi_*X)s_*\bar{D}_j\pi_*Y	\\
&= s_*( (\nabla_{\pi_*Y}'\bar{D}_i)\pi_*X - (\nabla_{\pi_*X}'\bar{D}_i)\pi_*Y ) 
+ \alpha'(\pi_*Y,\bar{D}_i\pi_*X) - \alpha'(\pi_*X,\bar{D}_i\pi_*Y).
\end{align*}
Comparing the tangent and normal components we get the identities
$$(\nabla_{\pi_*Y}'\bar{D}_i)\pi_*X - (\nabla_{\pi_*X}'\bar{D}_i)\pi_*Y 
= (-1)^j\bar{\psi}(\pi_*Y)\bar{D}_j\pi_*X -(-1)^j\bar{\psi}(\pi_*X)\bar{D}_j\pi_*Y$$
and 
$$\alpha'(\pi_*Y,\bar{D}_i\pi_*X) = \alpha'(\pi_*X,\bar{D}_i\pi_*Y).$$
The first equation above gives us (b), while the second one means that $s$ 
is hyperbolic or elliptic with respect to $\bar{J}$, 
because $\bar{D}_i\in \text{span}\{\bar{I},\bar{J}\}$ for $i=1, 2$ 
and $\bar{D}_i$ is not a multiple of $\bar I$ for some $i$.  

The only thing left to prove in the direct statement is condition (c).  
Using that $\psi$ is projectable onto $\bar{\psi}$, 
item (vii) of Proposition \ref{p:equivalentcartanhyp} 
and  \eqref{e:envelopeSdif2} we obtain
\begin{align}\label{e:dpsi}
\text{d}\bar{\psi}(\bar{X},\bar{Y}) &= \text{d}\psi(X,Y)\\
&= \left<(A-\lambda I)D_2X,(A-\lambda I)D_1Y\right>
-\left<(A-\lambda I)D_1X,(A-\lambda I)D_2Y\right>\nonumber\\
&= \left< \left<D_2X,\text{grad}\,\lambda\right>\Psi \circ f 
- S_*D_2X, \left<D_1Y,\text{grad}\,\lambda\right>\Psi \circ f - S_*D_1Y\right>\nonumber\\
&\quad - \left< \left<D_1X,\text{grad}\,\lambda\right>\Psi \circ f 
- S_*D_1X, \left<D_2Y,\text{grad}\,\lambda\right>\Psi \circ f - S_*D_2Y\right>\nonumber\\
&= \left<S_*D_2X,S_*D_1Y\right> - \left<S_*D_1X,S_*D_2Y\right>\nonumber\\
&= \left<\bar{D}_2\bar{X},\bar{D}_1\bar{Y}\right>' 
- \left<\bar{D}_1\bar{X},\bar{D}_2\bar{Y}\right>'\nonumber.
\end{align}

Let us now prove the converse. 
Using \eqref{e:equivalent}, and taking into account condition (b) 
and the fact that $s$ is hyperbolic or elliptic, we have 
\begin{align}\label{e:equivalent3}
\Psi_*&f_*B(X,Y) + \theta(X,Y)\Psi_*N + \varphi(X,Y)\Psi \circ f  -\lambda^{-1}\theta(X,Y)w	\\
&\quad = s_*( (\nabla_{\pi_*Y}'\bar{D}_i)\pi_*X 
- (\nabla_{\pi_*X}'\bar{D}_i)\pi_*Y ) + \alpha'(\pi_*Y,\bar{D}_i\pi_*X) 
- \alpha'(\pi_*X,\bar{D}_i\pi_*Y)\nonumber\\
&\quad = (-1)^js_*\left(\bar{\psi}(\pi_*Y)\bar{D}_j(\pi_*X) 
- \bar{\psi}(\pi_*X)\bar{D}_j\pi_*Y\right) \nonumber \\
&\quad = (-1)^j\left(\psi(Y)S_*D_jX - \psi(X)S_*D_jY\right). \nonumber
\end{align}
From  \eqref{e:envelopeSdif2} we have 
\begin{align*}
(-1)^j&\left(\psi(Y)S_*D_jX - \psi(X)S_*D_jY\right)\\
&= (-1)^j\psi(Y)\left( \left<D_jX,\text{grad}\,\lambda\right>\Psi \circ f
-\Psi_*f_*(A-\lambda I)D_jX\right>\\
&\quad - (-1)^j\psi(X)\left(\left<D_jY,\text{grad}\,\lambda\right>\Psi \circ f
-\Psi_*f_*(A-\lambda I)D_jY\right).
\end{align*}
Therefore, if we arrange equation \eqref{e:equivalent3} with this new information, we end up with

$$\Psi_*f_*\tilde{B}(X,Y) + \theta(X,Y)\Psi_*N + \tilde{\varphi}(X,Y)\Psi \circ f  
-\lambda^{-1}\theta(X,Y)w	 = 0$$
where $\tilde{B}$ and $\tilde{\varphi}$ are proper modifications of $B$ and $\varphi$.  In particular,
$$0 = \theta(X,Y) = \lambda( \left<X,(A-\lambda I)D_iY\right> - \left<Y,(A-\lambda I)D_iX\right>),$$ 
for all projectable  $X$,$Y \in \Gamma(\Delta^\perp)$.  Thus $(A-\lambda I)D_i$ is symmetric.   

 Let $J \in \Gamma(\text{End}(\Delta^\perp))$ (respectively, $D_i \in \Gamma(\text{End}(\Delta^\perp))$
  be the horizontal lift of $\bar{J}$ (respectively, $\bar{D}_i$) and 
  $\psi$ the horizontal lift of $\bar{\psi}$.  Since 
  $\bar{D}_1$, $\bar{D}_2 \in \text{span}\,\{\bar{I},\bar{J}\}$ and $\pi_*|_{\Delta^\perp}$ 
  is an isomorphism, we have that $D_1$, $D_2 \in \text{span}\,\{I,J\}$ and $J^2 = \epsilon I$, 
  depending on whether $s$ is hyperbolic or elliptic.  Let us prove that $D_i$ and $\psi$ 
  satisfy (i) to (ix), and that $f$ is hyperbolic (respectively, elliptic) with respect to $J$.   
  Items (i) and (vii) are clear because $\psi$  projects to $\bar{\psi}$.  
  From item (a) we get item (ii), item (e) gives item (ix) and from item (d) we get item (viii). 
 
  To prove item (iii), first notice that $\nabla_T^h D_i = [D_i,C_T]$ for all $T \in \Gamma(\Delta)$, 
  because  $D_i$ is  projectable.  	 
  On the other hand,  $\nabla_T^h A = \nabla_T^h(A-\lambda I)$ and \eqref{e:covderhorA} give
\begin{align*}
\nabla_T^h(A-\lambda I)D_i &- (A-\lambda I)D_iC_T \\
&= (\nabla_T^h(A-\lambda I) - (A-\lambda I)C_T)D_i + (A-\lambda I)(\nabla_T^hD_i - [D_i,C_T])\\
&= 0.
\end{align*}
Hence $\nabla_T^h(A-\lambda I)D_i = (A-\lambda I)D_iC_T$. 
In particular, this implies that $(A-\lambda I)D_iC_T$ is symmetric.  Therefore
$(A-\lambda I)D_iC_T=(A-\lambda I)C_T D_i$, 
and item (iii) follows.
Since there is $i\in \{1, 2\}$ such that  $D_i = a_iI + b_iJ$ with $b_i$ not null,
it follows that $\nabla_T^h J = [J,C_T] = 0$.  
This easily implies that
$C(\Gamma(\Delta)) \leq \text{span}\{I,J\}$, hence $f$ is hyperbolic 
(respectively, elliptic) with respect to $J$.

Since $s$ is either hyperbolic or elliptic with respect to $\bar{J}$ 
and $\bar{D}_i \in \text{span}\{\bar{I},\bar{J}\}$, 
$$\alpha'(\bar{D}_i\pi_*X,\pi_*Y) = \alpha'(\pi_*X,\bar{D}_i\pi_*Y).$$
From  \eqref{e:equivalent},  the symmetry of $(A-\lambda I)D_i$ and the fact that $\theta=0$ we get
\begin{align}\label{e:converseeq1}
\Psi_*f_*&\left( \left(\nabla_X(A-\lambda I)D_i\right)Y 
- \left(\nabla_Y(A-\lambda I)D_i\right)X - X\wedge Y(D_i^t\text{grad}\,\lambda) \right)\\
&+ \left(\left<(\nabla_Y D_i)X - (\nabla_X D_i)Y,\text{grad}\,\lambda\right> 
+ \text{Hess}\,\lambda(D_iX,Y) - \text{Hess}\,\lambda(X,D_iY) \right)\Psi\!\circ\!f\nonumber\\
& -\lambda \left(\left<(A-\lambda I)X,(A-\lambda I)D_iY\right> 
- \left<(A-\lambda I)D_iX,(A-\lambda I)Y\right> \right)\Psi\circ f \nonumber\\
&= s_*\left( \left(\nabla_{\pi_*Y}'\bar{D}_i\right)\pi_*X 
- \left(\nabla_{\pi_*X}'\bar{D}_i\right)\pi_*Y	 \right)\nonumber.
\end{align}
Using item (b),  \eqref{e:envelopeSdif2} and the fact that $D_1$, $D_2$ and $\psi$  
project to $\bar{D}_1$, $\bar{D}_2$ and $\bar{\psi}$, respectively, we obtain
\begin{align}\label{e:converseeq2}
s_*&\left( \left(\nabla_{\pi_*Y}'\bar{D}_i\right)\pi_*X 
- \left(\nabla_{\pi_*X}'\bar{D}_i\right)\pi_*Y	 \right) \\
&\qquad\qquad= (-1)^j s_*\left( \bar{\psi}(\pi_*Y)\bar{D}_j\pi_*X 
- \bar{\psi}(\pi_*X)\bar{D}_j\pi_*Y\right)\nonumber\\
&\qquad\qquad= (-1)^j\psi(Y)S_*D_jX - (-1)^j\psi(X)S_*D_jY\nonumber	\\
&\qquad\qquad= (-1)^j\psi(Y)\left( \left<D_jX,\text{grad}\,\lambda\right>\Psi \circ f 
- \Psi_*f_*(A-\lambda I)D_jX\right) \nonumber\\
&\qquad\qquad\quad- (-1)^j\psi(X)\left( \left<D_jY,\text{grad}\,\lambda\right>\Psi \circ f 
- \Psi_*f_*(A-\lambda I)D_jY\right)\nonumber.
\end{align}
Combining equations \eqref{e:converseeq1} and \eqref{e:converseeq2}, we get
\begin{align*}
0 &= \Psi_*f_*\left( \left(\nabla_X(A-\lambda I)D_i\right)Y 
- \left(\nabla_Y(A-\lambda I)D_i\right)X - X\wedge Y(D_i^t\text{grad}\,\lambda) \right)\\
&\quad + (-1)^j\Psi_*f_*(A-\lambda I)\left(\psi(Y)D_jX - \psi(X)D_jY\right)\\
&\quad + (-1)^j\left(\psi(X)\left<D_jY,\text{grad}\,\lambda\right> 
- \psi(Y)\left<D_jX,\text{grad}\,\lambda\right>\right)\Psi \circ f\\
&\quad+ \left(\left<(\nabla_Y D_i)X - (\nabla_X D_i)Y,\text{grad}\,\lambda\right> 
+ \text{Hess}\,\lambda(D_iX,Y) - \text{Hess}\,\lambda(X,D_iY) \right)\Psi \circ f\nonumber\\
&\quad -\lambda \left(\left<(A-\lambda I)X,(A-\lambda I)D_iY\right> 
- \left<(A-\lambda I)D_iX,(A-\lambda I)Y\right> \right)\Psi\circ f. \nonumber
\end{align*}
Taking into account the symmetry of  $(A-\lambda I)D_i$,  items (iv) and (v) 
of Proposition \ref{p:equivalentcartanhyp} follow.  Going the other way around in 
\eqref{e:dpsi} gives us (vii).
\qed

\section{The Subset $\mathcal{C}_s$}
\label{ch:subset}

This section is devoted to  characterizing hyperbolic and elliptic surfaces 
$s\colon L^2 \to \mathbb{S}^{n+2}_{1,1}$ that admit a triple $(\bar{D}_1,\bar{D}_2,\bar{\psi})$ 
satisfying items (a) to (e) of Lemma \ref{l:reductionlemma}.
We follow closely the proof of Proposition 9 in \cite{mD2013}.

Let us start with the case in which $s\colon L^2 \to \mathbb{S}^{n+2}_{1,1}$ is an hyperbolic 
surface with respect to the tensor $\bar{J}$. 
Let $(u,v)$ be local coordinates  whose coordinate vector fields $\{\partial_u,\partial_v\}$ 
are eigenvectors of $\bar{J}$ with eigenvalues $1$ and $-1$, respectively.  Then 
$$\alpha'(\partial_u,\partial_v) = \alpha'(J\partial_u,\partial_v) 
= \alpha'(\partial_u,J \partial_v) = - \alpha'(\partial_u,\partial_v),$$ 
 hence
$\alpha'(\partial_u,\partial_v) = 0$.
The coordinates $(u,v)$ are called \emph{real-conjugate coordinates}\index{real-conjugate coordinates}.  
Define the Christoffel symbols $\Gamma^1$ and $\Gamma^2$ with respect to the frame 
$\{\partial_u,\partial_v\}$  by
\begin{equation}\label{e:definitionchristsymbolsreal}
\nabla_{\partial_u}\partial_v = \Gamma^1\partial_u + \Gamma^2\partial_v.
\end{equation}
  Denote  $F = \left<\partial_u,\partial_v\right>$
and define the differential operator  
\begin{equation}\label{e:definitionqreal}
Q(\theta) = \text{Hess}\,\theta(\partial_u,\partial_v) + F\theta = \theta_{uv} 
- \Gamma^1\theta_u - \Gamma^2\theta_v + F\theta.
\end{equation}

For each pair of smooth functions $U = U(u)$ and $V = V(v)$, define 
\begin{equation}\label{e:definitionvarphireal}
\varphi^U(u,v) = U(u)e^{2\int_0^v\Gamma^1(u,s)\text{d}s} \quad \text{and} \quad \phi^V(u,v)
 = V(v)e^{2\int_0^u\Gamma^2(s,v)\text{d}s}.
\end{equation}
These functions satisfy
\begin{equation}\label{e:diffequationprojected}
\varphi^U_v  = 2\Gamma^1\varphi^U \quad \text{and} \quad \phi^V_u =  2\Gamma^2\psi^V
\end{equation}
with initial conditions $\varphi^U(u,0) = U(u)$ and $\phi^V(0,v) = V(v)$.  
 Assume, in addition, that one of the following conditions holds:
\begin{equation}\label{e:3}
U , V > 0 \quad \text{or} \quad 0 < 2\varphi^U 
< -(2\phi^V + 1) \quad \text{or} \quad 0 < 2\phi^V < - (2\varphi^U + 1).
\end{equation}
Under one of these conditions, one can define
\begin{equation}\label{e:definitionrhoreal}
\rho^{UV} = \sqrt{|2(\varphi^U + \phi^V)+1|}
\end{equation}
and
$$\mathcal{C}_s=\left\{(U,V): \eqref{e:3}\, \text{holds and}\, Q\left(\rho^{UV}\right) = 0\right\}.$$

Let us now suppose that $s:L^2 \to \mathbb{S}^{n+2}_{1,1}$ is an elliptic  surface with respect to a tensor ${J}$. 
 Let $(u,v)$ be local coordinates whose coordinate vector fields satisfy $J\partial_u = \partial_v$ and 
 $J\partial_v = -\partial_u$.  Extend  $J$, $\nabla$ and $\alpha^s$ $\mathbb{C}$-linearly. 
Denoting $\partial_z = (\partial_u - i\partial_v)/2$ and $\partial_{\bar{z}} = (\partial_u + i\partial_v)/2$,  
we have $J\partial_z= i\partial_z\;\;\;\mbox{and}\;\;\;J\partial_{\bar{z}} =-i\partial_{\bar{z}}$. Then
$$i\alpha^s(\partial_z,\partial_{\bar{z}})=\alpha^s(J\partial_z,\partial_{\bar{z}}) = \alpha^s(\partial_z,J\partial_{\bar{z}}) =-i\alpha^s(\partial_z,\partial_{\bar{z}}),$$
so, $\alpha^s(\partial_z,\partial_{\bar{z}})=0$.  
The coordinates $(u,v)$ are now called \emph{complex-conjugate}.

We can define a complex-valued Christoffel symbol $\Gamma\colon W\subset L^2 \to \mathbb{C}$ by
$$\nabla_{\partial_z}\partial_{\bar{z}} = \Gamma \partial_z + \bar{\Gamma}\partial_{\bar{z}}.$$
Set $F = \left<\partial_z,\partial_{\bar{z}}\right>$, where $\left<\,,\right>$ is the 
$\mathbb{C}$-bilinear extension of the metric induced by $s$, and define the differential operator
$$Q(\theta) = \text{Hess}\,\theta(\partial_z,\partial_{\bar{z}}) + F\theta 
= \theta_{z\bar{z}} - \Gamma\theta_z - \bar{\Gamma}\theta_{\bar{z}} + F\theta,$$
where $\theta:W \subset L^2 \to \mathbb{C}$ is a smooth function.
For each holomorphic function $\zeta$, let $\varphi^\zeta(z,\bar{z})$ 
be the unique complex valued function defined by
$$\varphi_{\bar{z}}^\zeta = 2\Gamma \varphi^\zeta \quad \text{and} 
\quad \varphi^\zeta(z,0) = \zeta(z).$$
Assume further that 
\begin{equation}\label{e:4}
\varphi^\zeta \neq -\frac{1}{2} \quad \text{and} \quad 4\text{Re}\,(\varphi^\zeta) + 1 < 0
\end{equation}
and define 
$$\rho^\zeta = \sqrt{-(4\text{Re}\,(\varphi^\zeta) + 1)}$$
and
$$\mathcal{C}_s = \left\{ \zeta \,\, \text{holomorphic}:\, 
\text{equation} \, \eqref{e:4}\, \text{holds and}\,\, Q(\rho^\zeta)=0 \right\}.$$

We are now ready to state and prove the main result of the section.  

\begin{proposition}\label{p:paramet}
If $s\colon L^2 \to \mathbb{S}^{n+2}_{1,1}$ is an elliptic or hyperbolic surface, 
then there exists a triple $(\bar{D}_1,\bar{D}_2,\bar{\psi})$ satisfying 
all conditions in Lemma \ref{l:reductionlemma} if and only if $\mathcal{C}_s$ is nonempty.  
Distinct triples (up to signs and permutation) give rise to distinct elements 
of $\mathcal{C}_s$, and conversely.	
\end{proposition}

\proof
The proof will be divided into cases, depending on whether $s$ is hyperbolic or elliptic.  

\subsubsection{Hyperbolic case}

Assume that $s$ is hyperbolic with respect to $\bar J$, and let $(\bar{D}_1,\bar{D}_2,\bar{\psi})$ 
satisfy all conditions in Lemma \ref{l:reductionlemma}. 
Let $(u,v)$ be  real-conjugate coordinates whose  coordinate vector fields are eigenvectors of $\bar J$, 
and hence of $\bar{D}_i$, $1\leq i\leq 2$, for $\bar{D}_1$, ${\bar{D}_2 \in \text{span}\{\bar{I},\bar{J}\}}$.  
From condition (a), we can suppose that the endomorphisms $\bar{D}_i$ are represented in this basis by 
\begin{equation}\label{e:21}
\sqrt{2}\bar{D}_1 = 
\begin{pmatrix}
\theta_1 & 0 \\
0 & 1/\theta_1	
\end{pmatrix}
\quad \text{and} \quad
\sqrt{2}\bar{D}_2 =
\begin{pmatrix}
\theta_2 & 0\\
0 & 1/\theta_2	
\end{pmatrix}.
\end{equation}
From item (e), that is, the assumption that $\text{rank}\,\bar{D}_1^2 + \bar{D}_2^2 - \bar{I} = 2$,  and
$$(\sqrt{2}\bar{D}_1)^2 + (\sqrt{2}\bar{D}_2)^2 - 2\bar{I} = 
\begin{pmatrix}
\theta_1^2 + \theta_2^2 - 2 & 0\\
0 & 1/\theta_1^2 + 1/\theta_2^2 - 2
\end{pmatrix}
$$
we infer that $\theta_1^2 + \theta_2^2 \neq 2$ and $1/\theta_1^2 + 1/\theta_2^2 \neq 2$.  
Also, from item (d), we get $\theta_1 \neq \pm \theta_2$.    
The equation of item (b) can be written as
$$\nabla_{\partial_u}' \bar{D}_i\partial_v - \nabla_{\partial_v}'\bar{D}_i\partial_u 
= (-1)^j\left(\bar{\psi}^u\bar{D}_j\partial_v - \bar{\psi}^v\bar{D}_j\partial_u\right),\, i \neq j,$$
where $\bar{\psi}^u = \bar{\psi}(\partial_u)$ and $\bar{\psi}^v = \bar{\psi}(\partial_v)$.  Therefore
$$\nabla_{\partial_u}'\theta_i^{-1}\partial_v - \nabla_{\partial_v}'\theta_i\partial_u 
= (-1)^j\left(\bar{\psi}^u\theta_j^{-1}\partial_v - \bar{\psi}^v\theta_j\partial_u\right),\, i \neq j,$$
and hence
\begin{align*}
-\frac{(\theta_i)_u}{\theta_i^2}&\partial_v + \theta_i^{-1}(\Gamma^1\partial_u 
+ \Gamma^2\partial_v) - (\theta_i)_v\partial_u - \theta_i(\Gamma^1\partial_u + \Gamma^2\partial_v) \\
&= (-1)^j\left(\bar{\psi}^u\theta_j^{-1}\partial_v - \bar{\psi}^v\theta_j\partial_u\right),\, i \neq j.
\end{align*}
From the equality of the components of both sides of the preceding equation with respect 
to the coordinate vector fields, we get that item (b) is equivalent to 
the  system of partial differential equations
\begin{equation}
\frac{(\theta_i)_u}{\theta_i^2}	+ \left(\theta_i-\frac{1}{\theta_i}\right)\Gamma^2 
= - (-1)^j\frac{\bar{\psi}^u}{\theta_j},
\end{equation}
\begin{equation}
(\theta_i)_v + \left(\theta_i - \frac{1}{\theta_i}\right)\Gamma^1 = (-1)^j\bar{\psi}^v\theta_j,
\end{equation}
with $i \neq j$.  Defining $\tau_i = \theta_i^2$, and multiplying the first equation 
by $-2/\theta_i$ and the second equation by $2\theta_i$, the preceding system becomes
\begin{equation}\label{e:24}
\left(\frac{1}{\tau_i}\right)_u + 2\left(\frac{1}{\tau_i}-1\right)\Gamma^2 
= 2(-1)^j\frac{\bar{\psi}^u}{\theta_1\theta_2}	,
\end{equation}
\begin{equation}\label{e:25}
(\tau_i)_v + 2(\tau_i - 1)\Gamma^1 = 2(-1)^j\bar{\psi}^v\theta_1\theta_2, \quad 1\leq i \neq j \leq 2.	
\end{equation}
Considering  \eqref{e:24} for the cases $i=1$ and $i=2$ and summing them up yields
\begin{equation}
\left(\frac{1}{\tau_1} + \frac{1}{\tau_2}\right)_u 
+ 2\left(\frac{1}{\tau_1}+\frac{1}{\tau_2}-2\right)\Gamma^2 = 0.	
\end{equation}
With the same procedure, but using instead  \eqref{e:25}, we get
\begin{equation}
(\tau_1 + \tau_2)_v + 2(\tau_1 + \tau_2 -2)\Gamma^1 = 0.
\end{equation}
Defining $\alpha = \tau_1 + \tau_2$ and $\beta = 1/\tau_1 + 1/\tau_2$,  
one can write the preceding equations  as
\begin{equation}\label{e:31}
\beta_u + 2(\beta -2)\Gamma^2 = 0 \quad \text{and} \quad \alpha_v + 2(\alpha -2)\Gamma^1 = 0.
\end{equation}
From the definition of $\tau_i$ we have that $\alpha$, $\beta > 0$.  
Moreover, since $\theta_1^2 \neq \theta_2^2$,
we have that $\tau_1$ and $\tau_2$ are distinct real roots of 
 $$\tau^2 - \alpha\tau +(\alpha/\beta)=0.$$
Thus  $\alpha\beta > 4$ and
\begin{equation}\label{e:32}
2\tau_i = \alpha - (-1)^i\sqrt{\frac{\alpha}{\beta}\left(\alpha\beta -4\right)}, \quad 1\leq i \leq 2.	
\end{equation}
Since $\theta_1^2 + \theta_2^2 \neq 2$ and $1/\theta_1^2 + 1/\theta_2^2 \neq 2$, 
we have that $\alpha \neq 2$ and $\beta \neq 2$.  Then, we can define
\begin{equation}
\varphi = \frac{1}{\alpha -2} \quad \text{and}	\quad \phi = \frac{1}{\beta-2}.
\end{equation}
From $\alpha > 0$, $\beta > 0$, $\alpha\beta - 4 > 0$,
$$\alpha = 2 + \frac{1}{\varphi} \quad \text{and} \quad \beta = 2 + \frac{1}{\phi},$$
and noticing that $\varphi$ and $\phi$ cannot be both negative, we get 
$$0 < \frac{2}{\varphi} + \frac{2}{\phi} + \frac{1}{\varphi\phi} 
= \frac{1}{\varphi\phi}\left(2\phi + 2\varphi + 1\right),$$
and hence $(\varphi,\phi)$ satisfies \eqref{e:3}. Moreover, 
$$\frac{\varphi_v}{\varphi} = -\frac{\alpha_v}{\alpha-2} 
\quad \text{and} \quad \frac{\phi_u}{\phi} = -\frac{\beta_u}{\beta -2},$$
so, from  \eqref{e:31} we get
$$\frac{\varphi_v}{\varphi} = 2\Gamma^1 \quad \text{and} \quad \frac{\phi_u}{\phi} 
= 2\Gamma^2.$$
Now, differentiating
$\bar{\psi} = \bar{\psi}^u\text{d}u + \bar{\psi}^v\text{d}v$  we get
$$2\text{d}\bar{\psi}(\partial_u,\partial_v) 
= 2(\bar{\psi}^v_u-\bar{\psi}^u_v)\text{d}u\wedge\text{d}v(\partial_u,\partial_v) 
= 2(\bar{\psi}^v_u-\bar{\psi}^u_v).$$
On the other hand,
$$\big<\sqrt{2}\bar{D}_2\partial_u,\sqrt{2}\bar{D}_1\partial_v\big>
-\big<\sqrt{2}\bar{D}_1\partial_u,\sqrt{2}\bar{D}_2\partial_v\big> 
= \left(\frac{\theta_2}{\theta_1}-\frac{\theta_1}{\theta_2}\right)F 
= \frac{\tau_2 - \tau_1}{\theta_1\theta_2}F.$$
Therefore,  item (c) is equivalent to
\begin{equation}\label{e:28}
2(\bar{\psi}^v_u-\bar{\psi}^u_v)	= \frac{\tau_2 - \tau_1}{\theta_1\theta_2}F.
\end{equation}
Set
\begin{equation}\label{e:csrho}
\rho = \sqrt{|2(\varphi + \phi)+1|} = \sqrt{\left|\frac{2}{\alpha-2} 
+ \frac{2}{\beta-2} + 1\right|} = \frac{\sqrt{\alpha\beta - 4}}{\sqrt{|(\alpha-2)(\beta-2)}|}.
\end{equation}
We want to show now that 
\begin{equation}\label{e:qrho}
Q(\rho) = \rho_{uv} - \Gamma^1\rho_u - \Gamma^2\rho_v + F\rho
 = 0.
 \end{equation}
In order to do so, we  express the functions $\rho$, $\Gamma^1$ and 
$\Gamma^2$ in terms of $\theta_i$. 
Using  \eqref{e:24} and \eqref{e:25} we get
\begin{equation}\label{e:Csgamma1}
\Gamma^1=-\frac{\theta_1(\theta_1)_v + \theta_2(\theta_2)_v}{\theta_1^2 + \theta_2^2 -2},
\end{equation}
\begin{equation}\label{e:Csgamma2}
\Gamma^2 = -\frac{\theta_1^3(\theta_2)_u 
+ \theta_2^3(\theta_1)_u}{\theta_1\theta_2(2\theta_2^2\theta_1^2-\theta_2^2-\theta_1^2)},
\end{equation}
\begin{equation}\label{e:Cspsiu}
\bar{\psi}^u = \frac{(\theta_2)_u\theta_1^3 - (\theta_1)_u\theta_2^3 
- (\theta_2)_u\theta_1 + (\theta_1)_u\theta_2}{2\theta_2^2\theta_1^2 - \theta^2_2 - \theta_1^2}
\end{equation}
and
\begin{equation}\label{e:Cspsiv}
\bar{\psi}^v=-\frac{(\theta_2)_v\theta_2\theta_1^2
 - (\theta_1)_v\theta_2^2\theta_1- \theta_2(\theta_2)_v 
 + \theta_1(\theta_1)_v }{\theta_1\theta_2(\theta_1^2 +\theta_2^2 - 2)}.
\end{equation}
From \eqref{e:28}  we obtain
\begin{equation}\label{e:CsF} 
F = \frac{2\theta_1\theta_2(\bar{\psi}_u^v -\bar{\psi}_v^u)}{\theta_2^2 - \theta_1^2}.
\end{equation}
Lastly, using  \eqref{e:csrho} we get
\begin{equation}\label{e:Csrho}
\rho = \sqrt{\frac{(\theta_1^2 + \theta_2^2)^2/\theta_1^2\theta_2^2-4}{\left|(\theta_1^2+\theta_2^2 - 2)
(\frac{1}{\theta_1^2}+\frac{1}{\theta_2^2}-2)\right|}}.
\end{equation}
Using the preceding identities, a long but straightforward computaion 
shows that \eqref{e:qrho} is satisfied. Thus, the set $\mathcal{C}_s$ is non-empty.

Now we prove the converse statement.  Since $s\colon L^2 \to \mathbb{S}^{n+2}_{1,1}$
 is hyperbolic, there exist real conjugate coordinates $(u,v)$. If $(U,V) \in \mathcal{C}_s$, then 
$$\varphi^U(u,v) = U(u)e^{2\int_0^v\Gamma^1(u,s)\text{d}s} \quad \text{and} \quad \phi^V(u,v) 
= V(v)e^{2\int_0^u\Gamma^2(s,v)\text{d}s}$$
must satisfy  \eqref{e:diffequationprojected} and, together with the functions $U$ and $V$, 
also satisfy  \eqref{e:3}.  From the definition of the set $\mathcal{C}_s$, we must have 
$Q(\rho) = 0$, where $\rho = \sqrt{|2(\varphi^U + \phi^V) +1|}$.  Set
$\alpha = 2 +{1}/{\varphi^U}$ and $\beta= 2 + {1}/{\phi^V}$,
which are well defined because $U$, $V$, $\varphi^U$ and $\phi^V$ satisfy one of the equations 
in \eqref{e:3}, and therefore, $\varphi^U$ and $\phi^V$ cannot vanish at any point.

Since $(\varphi^U, \phi^V)$ satisfies  \eqref{e:3}, we claim that $\alpha > 0$, 
$\beta > 0$ and $\alpha \beta - 4 > 0$.  In the first possiblity, namely, 
if $U$, $V > 0$, then $\varphi^U > 0$ and $\phi^V > 0$, and our claim follows 
from the definition of $\alpha$ and $\beta$.  If 
$0 < 2\varphi^U < -(2\phi^V + 1)$,
then we immediately see that $\alpha > 0$.  We also have $\psi^V < -1/2$, so $\beta > 0$.  
Lastly,
$$\alpha \beta - 4 =  \frac{2}{\varphi^U} +\frac{2}{\phi^V}  +\frac{1}{\varphi^U \phi^V}
= \frac{1}{\varphi^U \phi^V}\left(2\varphi^U + 2\phi^V + 1\right).$$
Since, $\varphi^U >0$, $\phi^V < 0$ and $2\varphi^U + 2\phi^V + 1 < 0$, 
we conclude that $\alpha \beta - 4 > 0$.  The other case is symmetric, so our claim is proved.

With the information that $\alpha > 0$, $\beta > 0$ and $\alpha \beta - 4 > 0$, we can define 
the functions $\tau_i$ by  \eqref{e:32}, that is, $\tau_1$ and $\tau_2$ are the  
roots of   $\tau^2 - \alpha\tau + \alpha/\beta = 0$. We  conclude that 
$\tau_1 + \tau_2 = \alpha$ and $\tau_1\tau_2 = \alpha/\beta$.

As before, write $\tau_i=(\theta_i)^2$ and let $\bar{\psi}^u$ and $\bar{\psi}^v$ be given 
by \eqref{e:24} and \eqref{e:25}, respectively.  Substituting $\tau_i$ by $\theta_i^2$ 
in those equations, we arrive at the same equations as in the direct statement, 
so we can express $\Gamma^1$, $\Gamma^2$, $\bar{\psi}^u$ and $\bar{\psi^v}$ in terms of the 
$\theta_i$ by the identities \eqref{e:Csgamma1}, \eqref{e:Csgamma2}, \eqref{e:Cspsiu} and \eqref{e:Cspsiv}.  
From the fact that $\tau_1 + \tau_2 = \alpha$ and $\tau_1\tau_2 = \alpha/\beta$, we get 
$\alpha=\theta_1^2 + \theta_2^2$ and $\beta= 1/\theta_1^2 + 1/\theta_2^2$.
From the definition of $\rho$, we have that  \eqref{e:csrho} is valid, and so, 
replacing $\alpha$ and $\beta$ is terms of the $\theta_i$, we also obtain  \eqref{e:Csrho}.  
Since $\rho\neq 0$ at any point,  from $Q(\rho)=0$ we obtain
\begin{equation}\label{e:CsF2}
F = -\frac{\rho_{uv}-\Gamma^1\rho_u - \Gamma^2\rho_v}{\rho}
\end{equation}
which can  be written in terms of the $\theta_i$ using  \eqref{e:Csrho}, \eqref{e:Csgamma1} and \eqref{e:Csgamma2}.
Using those identities,
a long but straightforward computation shows that \eqref{e:28} is satisfied.   

Let $\bar{D}_1$ and $\bar{D}_2$ be defined by  \eqref{e:21} with respect to the frame 
$\{\partial_u,\partial_v\}$, and set $\bar{\psi} = \bar{\psi}^u\text{d}u + \bar{\psi}^v\text{d}v$.  
Then condition (a) is clear from the definition of $\bar{D}_i$, whereas condition (b) 
follows from   \eqref{e:24} and \eqref{e:25}.  Condition (c) is a consequence of  \eqref{e:28}.  
Since $\alpha>0$, we have $\tau_1 \neq - \tau_2$, so $\bar{D}_1^1 \neq - \bar{D}_2^2$.  
Because the discriminant is $\alpha\beta - 4> 0$, $\tau_1$ and $\tau_2$ are not equal, 
so $\bar{D}_1^1 \neq \bar{D}_2^2$, and item (d) is proved.  From the definition of $\alpha$ 
and $\beta$ we cannot have $\alpha = 2$ or $\beta = 2$, so item (e) follows.  
Distinct pairs ($\varphi,\phi$) give rise to distinct 4-tuples $(\tau^1,\tau^2,\bar{\psi}^u,\bar{\psi}^v)$, 
and hence to distinct triples  $(\bar{D}_1,\bar{D}_2,\bar{\psi})$.  This completes the proof for the hyperbolic case.

\subsubsection{Elliptic case}

Suppose $s:L^2 \to \mathbb{S}^{n+2}_{1,1}$ is an elliptic surface, and that there exists a 
triple $(\bar{D}_1,\bar{D}_2,\bar{\psi})$ satisfying all conditions in Lemma \ref{l:reductionlemma}.  
Since we will use complex conjugate operation, let us omit the bar notation  just for now.

Let $(u,v)$  be complex-conjugate coordinates on $L^2$. Then  
$\partial_z =(1/2)(\partial_u - i\partial_v)$ and $\partial_{\bar{z}} = (1/2)(\partial_u + i\partial_v)$ 
are eigenvectors of the complex linear extension of the tensor $J$ with eigenvalues $i$ and $-i$, respectively.  
From item (a) of Lemma \ref{l:reductionlemma} we can assume that 
$\sqrt{2}D_i = a_iI + b_iJ,$
 where $a_i^2 + b_i^2 =1$.  Then the complex-linear extensions  of $D_1$ and $D_2$, 
 which we  denote by the same symbols, are given with respect to the frame $\{\partial_z, \partial_{\bar{z}}\}$ by
\begin{equation}\label{e:34}
\sqrt{2}D_1=
\begin{pmatrix}
\theta_1 & 0\\
0& \bar{\theta}_1	
\end{pmatrix}
\quad \text{and} \quad	
\sqrt{2}D_2=
\begin{pmatrix}
\theta_2 & 0\\
0& \bar{\theta}_2
\end{pmatrix},
\end{equation}
where $\theta_i:L^2 \to \mathbb{S}^1$. Moreover,  from item (d) of Lemma \ref{l:reductionlemma}, 
we must have $\theta_1 \neq \pm \theta_2$.  

Set $\psi^z = \psi(\partial_z)$, $\psi^{\bar{z}} = \psi(\partial_{\bar{z}})= \psi^{\bar{z}}$
and define a complex-valued Christoffel symbol $\Gamma$ by
$$\nabla_{\partial_z}\partial_{\bar{z}} = \Gamma \partial_z + \bar{\Gamma}\partial_{\bar{z}}.$$  
Define $\tau^i =\theta_i^2$, $1 \leq i \leq 2$.  Then, from item (b) of Lemma \ref{l:reductionlemma}
we get
$$\nabla_{\partial_z}\bar{\theta}_i\partial_{\bar{z}} - \nabla_{\partial_{\bar{z}}}\theta_i\partial_z
 = (-1)^j\left(\psi^z \bar{\theta}_j\partial_{\bar{z}} - \psi^{\bar{z}}\theta_j\partial_z\right),$$
which is equivalent to
$$(\bar{\theta}_i)_z\partial_{\bar{z}} + \bar{\theta}_i\left(\Gamma\partial_z 
+ \bar{\Gamma}\partial_{\bar{z}}\right) - (\theta_i)_{\bar{z}}\partial_z 
- \theta_i\left(\Gamma\partial_z + \bar{\Gamma}\partial_{\bar{z}}\right) 
= (-1)^j\left(\psi^z\bar{\theta}_j\partial_{\bar{z}} - \psi^{\bar{z}}\theta_j\partial_z\right).$$
We obtain that
\begin{equation}\label{e:35pre}
(\theta_i)_{\bar{z}} - \bar{\theta}_i\Gamma + \theta_i\Gamma=(-1)^j\psi^{\bar{z}}\theta_j.
\end{equation} 
Multiplying both sides of  \eqref{e:35pre} by $2\theta_i$
we get
\begin{equation}\label{e:35}
(\tau_i)_{\bar{z}} + 2\left(\tau_i - 1\right)\Gamma=2(-1)^j\psi^{\bar{z}}\theta_1\theta_2.
\end{equation}
Now we  use item (c) of Lemma \ref{l:reductionlemma}.  On one hand, since 
$\text{d}\psi = (\psi^v_u - \psi^u_v)\text{d}u\wedge\text{d}v$, we obtain that
$ 2\text{d}\psi(\partial_z,\partial_{\bar{z}})= -4i\text{Im}\,\psi^z_{\bar{z}}$. 
  On the other hand,
$$\left<\sqrt{2}D_2\partial_z,\sqrt{2}D_1\partial_{\bar{z}}\right> 
- \left<\sqrt{2}D_1\partial_z,\sqrt{2}D_2\partial_{\bar{z}}\right> 
= \left(\bar{\theta}_1\theta_2 - \theta_1\bar{\theta}_2\right)F = \frac{\tau_2 -\tau_1}{\theta_1\theta_2}F.$$
Using item (c) of Lemma \eqref{l:reductionlemma} and multiplying both sides by $i$, we get
\begin{equation}\label{e:36}
4\text{Im}\,\psi^z_{\bar{z}} = i\frac{\tau_2 -\tau_1}{\theta_1\theta_2}F.
\end{equation}
Defining $\alpha = \tau_1 + \tau_2$, and  summing up cases $i=1$ and $i=2$ in  \eqref{e:35} yield
\begin{equation}\label{e:Csdiffalpha}
\alpha_{\bar{z}} + 2(\alpha -2)\Gamma = 0.
\end{equation}
Because $\theta_i \in S^1$, also $\tau_i \in S^1$.  From condition (d) in Lemma \ref{l:reductionlemma}, 
we have $\tau_i \neq \pm \tau_2$. Hence, $0 < |\alpha| = |\tau_1 + \tau_2| < 2$.  Thus $\varphi = {1}/{(\alpha -2)}$
is well defined and satisfies
$$\frac{\varphi_{\bar{z}}}{\varphi} =  - \frac{\alpha_{\bar{z}}}{\alpha-2} = 2\Gamma.$$
Since
$$4\text{Re}\,\varphi +1 = 2\frac{\alpha + \bar{\alpha} - 4 }{|\alpha-2|^2} + 1 
= \frac{|\alpha|^2-4}{|\alpha-2|^2}$$
and $|\alpha| < 2$, we conclude that $4\text{Re}\,\varphi +1 <0$.  
Since $\alpha \neq 0$, we have $\varphi \neq -1/2$, and the conditions in  \eqref{e:4} follow.  
From $\tau_1 + \tau_2 = \alpha$, $\tau_i \in \mathbb{S}^1$ and 
$$\tau_1\tau_2 = \frac{\tau_1 + \tau_2}{1/\tau_1+1/\tau_2} 
= \frac{\tau_1 + \tau_2}{\bar{\tau}_1 + \bar{\tau}_2} = \frac{\alpha}{\bar{\alpha}},$$
we obtain that
\begin{equation}\label{e:37}
\tau_j = \frac{\alpha}{2}\left(1 - (-1)^ji\frac{\sqrt{4-|\alpha|^2}}{|\alpha|}\right).
\end{equation}
In order to show that $\mathcal{C}_s$ is non-empty, we must prove that 
\begin{equation}\label{e:rhoellip}
\rho = \sqrt{-\left(4\text{Re}\,\varphi + 1\right)} = \frac{\sqrt{4-|\alpha|^2}}{|\alpha -2|}
\end{equation}
satisfies $Q(\rho) = 0$.  For that, as in the hyperbolic case we express $\Gamma$, $\psi^{\bar{z}}$, 
$F$ and $\rho$ in terms of the functions $\theta_i$.  First, notice that 
$\alpha = \theta_1^2 + \theta_2^2$ and $\bar{\alpha} = 1/\theta_1^2 + 1/\theta_2^2$.  
From  \eqref{e:Csdiffalpha}, and replacing $\alpha$ in terms of $\theta_i$, we get
\begin{equation}\label{e:Csgammaellip}
\Gamma = - \frac{(\theta_1^2 + \theta_2^2)_{\bar{z}}}{2\left(\theta_1^2 + \theta_2^2 -2\right)} 
= -\frac{\theta_1(\theta_1)_{\bar{z}} + \theta_2(\theta_2)_{\bar{z}}}{\theta_1^2 + \theta_2^2 -2}.	
\end{equation}
Using this and  \eqref{e:35} with $i=1$ we obtain
\begin{equation}\label{e:Cspsibarz}
\psi^{\bar{z}} = \frac{\theta_1\theta_2^2(\theta_1)_{\bar{z}} - 
\theta_1^2\theta_2(\theta_2)_{\bar{z}} - \theta_1(\theta_1)_{\bar{z}} 
+ \theta_2(\theta_2)_{\bar{z}}}{\theta_1\theta_2\left(\theta_1^2 + \theta_2^2 -2\right)}.
\end{equation}
Observing that
$$(\psi^{\bar{z}})_z - (\psi^z)_{\bar{z}} = \overline{(\psi^z)_{\bar{z}}} 
- (\psi^z)_{\bar{z}} = -2i\text{Im}\,(\psi^z)_{\bar{z}}$$
and using  \eqref{e:36} we get
$$2\left((\psi^{\bar{z}})_z - (\psi^z)_{\bar{z}}\right) = -4i\text{Im}\,(\psi^z)_{\bar{z}} 
= \frac{\theta_2^2 - \theta_1^2}{\theta_1\theta_2.}F$$
Solving for $F$ yields
\begin{equation}\label{e:CsFellip}
F=\frac{2\theta_1\theta_2\left((\psi^{\bar{z}})_z - (\psi^z)_{\bar{z}}\right)}{\theta_2^2 - \theta_1^2}.
\end{equation}
From  \eqref{e:rhoellip} and the expression of $\alpha$ and $\bar{\alpha}$ in terms of $\theta_i$ we have
\begin{equation}\label{e:Csrhoellip}
\rho =  \sqrt{  \frac{4-(\theta_1^2+\theta_2^2)/\theta_1^2\theta_2^2}{\left(\theta_1^2+\theta_2^2-2\right)\left(1/\theta_1^2 + 1/\theta_2^2 - 2\right)}  } = i\sqrt{  \frac{(\theta_1^2+\theta_2^2)/\theta_1^2\theta_2^2-4}{|\left(\theta_1^2+\theta_2^2-2\right)\left(1/\theta_1^2 + 1/\theta_2^2 - 2\right)|}  } .
\end{equation}
If we compare the expressions we got for $\Gamma$, $\bar{\Gamma}$, $\psi^{\bar{z}}$, 
$\psi^z$, $F$ and $\rho$, except for constant multiple $i$ in the $\rho$, 
they are the same equations as \eqref{e:Csgamma1}, \eqref{e:Csgamma2}, \eqref{e:Cspsiu}, 
\eqref{e:Cspsiv}, \eqref{e:CsF} and \eqref{e:Csrho}  we have found in the hyperbolic case, 
when we replace $(z,\bar{z})$, $(\Gamma,\bar{\Gamma})$, $(\psi^z,\psi^{\bar{z}})$ 
for $(u,v)$, $(\Gamma^1,\Gamma^2)$ and $(\psi^u,\psi^v)$, respectively.  Therefore $Q(\rho) =0$,
as one can confirm by direct  computation. This shows that $\mathcal{C}_s$ is non-empty.

We now prove the converse.  Let $(u,v)$ be  complex-conjugate coordinates for 
$s\colon L^2 \to \mathbb{S}^{n+2}_{1,1}$.  If $\zeta \in \mathcal{C}_s$ 
is an holomorphic function, then  \eqref{e:4} holds for the complex-valued function 
$\varphi^\zeta(z,\bar{z})$ defined by
$\varphi^z_{\bar{z}} = 2\Gamma \varphi^\zeta$ and $\varphi^\zeta(z,0) = \zeta$. 
Moreover,  $\rho^\zeta=\sqrt{-(4\text{Re}\, \varphi^\zeta + 1)}$ satisfies $Q(\rho^\zeta) = 0$. 

Define 
$\alpha = 2 + {1}/{\varphi^\zeta}.$
From the first condition of  \eqref{e:4} we have that $\alpha$ is not null.  Since
$$|\alpha|^2=\alpha\bar{\alpha} 
= \left(2+ \frac{\overline{\varphi^\zeta}}{|\varphi^\zeta|^2}\right)\left(2+ \frac{\varphi^\zeta}{|\varphi^\zeta|^2}\right) 
= 4 + \frac{4\text{Re}\,\varphi^\zeta + 1}{|\varphi^\zeta|^2},$$ 
from the second condition of  \eqref{e:4} we get $|\alpha| < 2$.  

Let $\tau_1$ and $\tau_2$ be the roots of
$x^2 - \alpha x + \frac{\alpha}{\bar{\alpha}} = 0$. In particular, 
$\alpha =\tau_1 + \tau_2$.  From the definition of $\tau_j$, we have 
$$|\tau_j| = \frac{|\alpha|}{2}\sqrt{\left(1 + \frac{4-|\alpha|^2}{|\alpha|^2}\right)} = 1,$$
for $j=1, 2$.  Also, since
$|\alpha|<2$ we have $\tau_1 \neq \pm\tau_2$.
Write $\tau_j = \theta_j^2$,  define $\psi^{\bar{z}}$ by   \eqref{e:35} and then $\psi^u$ and $\psi^v$ by  
$\psi^u = 2\text{Re}\,\psi^{\bar{z}}$ and $\psi^v = 2\text{Im}\,\psi^{\bar{z}}$.  
Define the complex-linear extensions $\sqrt{2}D_j$ by  \eqref{e:34}.  
To recover the original $\sqrt{2}D_j$ just remember that $\sqrt{2}D_j = a_jI + b_jJ$ for $\theta_j = a_j + ib_j$.  
So, we get a triple $(D_1,D_2,\psi)$.  We have to show that this triple satisfies 
conditions (a) to (e) of Lemma \eqref{l:reductionlemma}.

Since $|\tau_j| = 1$, then $|\theta_j|=1$, and so $\det \sqrt{2} D_j = 1$. 
This gives (a).  Because  \eqref{e:35} is satisfied, item (b) follows.  
From the fact that $\tau_1 \neq \pm \tau_2$ and how $\tau_j$ is defined we get item (d).  
Now, it is easily seen that one can have 
$\text{rank}\, (\sqrt{2}D_1)^2 + (\sqrt{2}D_2)^2 - 2I <2$ only if $(\sqrt{2}D_1)^2 + (\sqrt{2}D_2)^2 - 2I=0$. 
Since $\theta_j = a_j + ib_j$ satisfies 
$|\theta_j|=1$, this easily implies that $b_1 = 0 = b_2$ and $a_j = \pm  1$.  
Therefore, $\theta_1 = \pm \theta_2$, a contradiction because $\tau_1 \neq \pm\tau_2$, which proves (e).

Let us prove item (c).  Since $\varphi^\zeta = 1/(\alpha-2)$, and from the definition 
of $\rho^\zeta$, we get  \eqref{e:rhoellip}.  Eq. \eqref{e:Csrhoellip} then follows 
from  $\alpha = \theta_1^2 + \theta_2^2$.  Since $\psi^{\bar{z}}$ and $\Gamma$ 
satisfy  \eqref{e:35}, we have the validity of  \eqref{e:Csgammaellip}
and \eqref{e:Cspsibarz}.
 From  $Q(\rho) = 0$ we get
\begin{equation}\label{CsFellipconv}
F = -\frac{-\rho_{z\bar{z}}- \Gamma\rho_z -\bar{\Gamma}\rho_{\bar{z}}}{\rho},
\end{equation}
so we can express $F$ in terms of $\theta_i$ using Eqs \eqref{e:Csrhoellip} and \eqref{e:Csgammaellip}.
  Notice that the $\rho$ used in the hyperbolic case differs from this $\rho$ by a multiple of $i$.  
  We arrive at the same equations as in proof of the converse statement  of the hyperbolic case, 
  with  $(z,\bar{z})$, $(\Gamma,\bar{\Gamma})$, $(\psi^z,\psi^{\bar{z}})$ instead of 
  $(u,v)$, $(\Gamma^1,\Gamma^2)$ and $(\psi^u,\psi^v)$, respectively.  
  Thus, equation \eqref{e:36} is valid,  and so is  item (c).  

Finally, notice that distinct $\zeta's$ give rise to distinct ${\varphi^{\zeta}}'s$, 
and so distinct $\alpha's$.  Since the $\tau_i$ are the roots of 
$x^2 - \alpha x + \frac{\alpha}{\bar{\alpha}} = 0$, we get distinct $\tau_i's$, 
hence distinct $\theta_i's$, and so  distinct triples $(D_1,D_2,\psi)$.  
\vspace{1ex}\qed

Before finishing the current section, we give  an explicit example of an hyperbolic surface 
$s\colon L^2 \to \mathbb{S}^{m}_{1,1}$ whose associated subset $\mathcal{C}_s$ is nonempty.

Let us start by orthogonally decomposing  
$\mathbb{L}^{m+1} = \mathbb{R}^{m_1} \times \mathbb{L}^{m_2}$ 
and considering a curve $\alpha:I_1 \to \mathbb{S}^{m_1-1} \subset \mathbb{R}^{m_1}$  
parametrized by arc length.  Denote  $\tilde{\alpha} = i \circ \alpha$, 
where $i\colon \mathbb{R}^{m_1} \to \mathbb{L}^{m+1}$ is the inclusion, 
and consider the flat parallel vector subbundle 
$\mathcal{L} \subset  N_{\tilde{\alpha}}I$ of 
$\text{rank}\,\, k= m_2+1$ whose fiber at $v \in I_1$ is 
\begin{equation}\label{e:ortdescp}
\mathcal{L}(v) = \mathbb{R}\tilde{\alpha}(v)\oplus \mathbb{L}^{m_2}.
\end{equation}
If $\{\xi_1, \cdots, \xi_k\}$ is an orthonormal frame of  parallel sections 
of $\mathcal{L}$, with $\xi_1(v)=\tilde{\alpha}(v)$, then we can 
define a parallel vector bundle isometry  $\phi\colon I_1 \times \mathbb{L}^{k} \to \mathcal{L}$ by
$$\phi(v,Y) = \phi_v(Y) = \sum_{i=1}^k Y^i \xi_i(v).$$
 Let $e \in \mathbb{L}^{k}$ be such that $\phi_{v}(e) = \tilde{\alpha}(v)=\xi_1(v)$ for all $v \in I_1$, 
 and denote 
$$\Omega^0(\tilde{\alpha}) = \{Y \in \mathbb{L}^{k}: \left<Y,e\right> > 0\}.$$
Consider $\beta\colon I_0 \to \mathbb{S}^{k-1}_{1,1} \cap \Omega^0(\tilde{\alpha}) \subset \mathbb{L}^k$, 
another curve parametrized by arc length.  
Define $s\colon I_0 \times I_1 \to \mathbb{S}^m_{1,1} \subset \mathbb{L}^{m+1}$ 
by $s(u,v) = \phi_v(\beta(u))$. Then 
$$s_*\partial_u = \phi_v(\beta'(u)) \quad \text{and} \quad s_*\partial_v 
= \left<\beta(u), e\right>\tilde \alpha'(v),$$
hence $s$ is an immersion with induced metric
$ds^2=du^2+\rho^2(u)dv^2,$
where $\rho(u)=\left<\beta(u), e\right>$.  
Moreover, differentiating, say, the first of the preceding equations with respect
to $v$ gives that $\alpha^s(\partial_u, \partial_v)=0$.

By  a suitable change of coordinates $\tilde{u}=\gamma(u)$, 
we can pass to isothermal coordinates with respect to which the metric is written as 
$$ds^2=e^{2\lambda(\tilde{u})}(\text{d}\tilde{u}^2 + \text{d}v^2)$$
for some smooth function $\lambda=\lambda(\tilde{u})$,   
and we still have $\alpha^s(\partial_{\tilde{u}},\partial_v)=0$. 
Thus, the surface $s$ is an hyperbolic surface and $(\tilde{u},v)$ are real-conjugate coordinates.  
For simplicity, we rewrite $\tilde{u}$ by $u$.

  Let us show that, for the above surface 
  $s\colon I_0 \times I_1 \to \mathbb{S}^m_{1,1} \subset \mathbb{L}^{m+1}$, 
  the subset $\mathcal{C}_s$ is non-empty.  If we define 
$$E = \left<\partial_u,\partial_u\right> = e^{2\lambda(u)}, \quad F 
= \left<\partial_u,\partial_v\right> = 0 \quad \text{and} \quad G = \left<\partial_v,\partial_v\right> = e^{2\lambda(u)},$$
then the Christoffel symbols $\Gamma^1$ and $\Gamma^2$ defined by  \eqref{e:definitionchristsymbolsreal} satisfy
$$0 = E_v = 2\Gamma^1 E\;\;\;\mbox{and}\;\;\;2\lambda'e^{2\lambda} = G_u = 2\Gamma^2G.$$
  Hence $\Gamma^1 = 0$ and $\Gamma^2 = \lambda'$.
Given a pair of smooth functions $\tilde{U} = \tilde{U}(u)$ and $V = V(v)$, 
the functions $\varphi^{\tilde{U}}$ and $\varphi^V$  defined in the hyperbolic case by  \eqref{e:definitionvarphireal}
are given by $\varphi^{\tilde{U}} = \tilde{U}$ and $\varphi^V = Ve^{2\lambda}$. 
By suitably modifying $\tilde{U}$ we have 
$\varphi^{\tilde{U}} = e^{2\lambda}U$ and $\varphi^V = e^{2\lambda}V$, 
so, taking into account the definition of $\rho$ (see \eqref{e:definitionrhoreal}), we obtain
$$\rho = \rho^{\tilde{U}V} = \sqrt{2e^{2\lambda}(U + V)+1}.$$
From the expression of $\Gamma^1$ and $\Gamma^2$, the operator $Q$ in \eqref{e:definitionqreal} reduces to
$$Q(\theta) = \theta_{uv} - \Gamma^1\theta_u - \Gamma^2\theta_v + F\theta = \theta_{uv} - \lambda'\theta_v.$$
Now,
$$\rho_v = \frac{e^{2\lambda}V_v}{\sqrt{2e^{2\lambda}(U+V)+1}},$$
and so
$$\rho_{uv} = \frac{2\lambda'e^{2\lambda}V_v\left(2e^{2\lambda}(U+V)+1\right)
 - V_ve^{2\lambda}\left(2\lambda'e^{2\lambda}(U+V) + e^{2\lambda}U_u\right)}{(2e^{2\lambda}(U+V)+1)^{3/2}},$$
which implies that  $0 = Q(\rho) = \rho_{uv} - \lambda'\rho_v$  reduces to  $V_v(2\lambda' - U_u e^{2\lambda})= 0$. 
This equation is satisfied for $V = k\in \mathbb{R}$ or for $U = c - e^{-2\lambda}$. 
Thus  $\mathcal{C}_s$ is nonempty.

 We point out that other examples of surfaces  $s\colon L^2 \to \mathbb{S}^m_{1,1}$ as above can be obtained 
 by considering other types of orthogonal decompositions in  \eqref{e:ortdescp}.
 
 \section{The Classification} 
\label{ch:classification}

We are now in a position to state and prove the classification of hypersurfaces 
$f\colon M^n \to \mathbb{R}^{n+1}$ that carry a principal curvature of 
multiplicity $n-2$ and admit a genuine conformal deformation $\tilde{f}\colon M^n \to \mathbb{R}^{n+2}$.

\begin{theorem}\label{t:class}
Let $f\colon M^n \to \mathbb{R}^{n+1}$ be a  hypersurface with a principal 
curvature of multiplicity $n-2$.  Assume that $f$ is not a Cartan hypersurface 
on any open subset of $M^n$ and that it   admits a genuine conformal deformation 
$\tilde{f}:M^n \to \mathbb{R}^{n+2}$. Then, on each connected component 
of an open dense subset of $M^n$, it envelops a two-parameter congruence of 
hyperspheres $s\colon L^2 \to \mathbb{S}^{n+2}_{1,1}$ which is either an elliptic 
or hyperbolic surface with non-empty associated set $\mathcal{C}_s$.

Conversely, any simply connected hypersurface $f$ that envelops  a two parameter 
congruence of hyperspheres $s\colon L^2 \to \mathbb{S}^{n+2}_{1,1}$ that is either
 an elliptic or hyperbolic surface and is such that the set $\mathcal{C}_s$ 
 is non-empty admits  genuine conformal deformations in $\mathbb{R}^{n+2}$ 
 which are parametrized by $\mathcal{C}_s$.
\end{theorem}

\proof Composing $f$ with an inversion in $\mathbb{R}^{n+1}$, if necessary, 
we may assume that the principal curvature of $f$ with multiplicity $n-2$ is nowhere
vanishing. By Proposition~\ref{p:equivalentcartanhyp}, on an open dense subset 
of $M^n$, the hypersurface is either elliptic or hyperbolic and admits
a triple $(D_1,D_2,\psi)$ satisfying all 
conditions in the statement of that result. 
By Lemma \ref{l:reductionlemma}, the two-paramenter congruence of hyperspheres 
$s\colon L^2 \to \mathbb{S}^{n+2}_{1,1}$  that is enveloped by  $f$ is either 
an elliptic or hyperbolic surface, respectively, and the triple $(D_1,D_2,\psi)$
projects down to a  triple $(\bar{D}_1,\bar{D}_2,\bar{\psi})$ on $L^2$ 
satisfying all conditions in that lemma.
We conclude from Proposition \ref{p:paramet} that $(\bar{D}_1,\bar{D}_2,\bar{\psi})$ 
gives rise to an element of $\mathcal{C}_s$.  

Conversely, suppose $f\colon M^n \to \mathbb{R}^{n+1}$ is a simply connected 
hypersurface that  envelops a two-parameter congruence of hyperspheres 
$s\colon L^2 \to \mathbb{S}_{1,1}^{n+2}$ that is either an elliptic or 
hyperbolic surface, and is such that the set $\mathcal{C}_s$ is non-empty.
  By Proposition \ref{p:paramet}, each element of $\mathcal{C}_s$ gives rise to 
a  triple $(\bar{D}_1,\bar{D}_2,\bar{\psi})$ on $L^2$ satisfying 
all conditions in Lemma \ref{l:reductionlemma}.  
Then, it follows from Lemma \ref{l:reductionlemma} that $f$ is either 
elliptic or hyperbolic, respectively, and that $(\bar{D}_1,\bar{D}_2,\bar{\psi})$ 
can be lifted to a  triple $(D_1,D_2,\psi)$ on $M^n$ satisfying 
all conditions in Proposition \ref{p:equivalentcartanhyp}.  
Proposition \ref{p:equivalentcartanhyp} then implies that each such triple yields a 
genuine conformal deformation $\tilde{f}\colon M^n \to \mathbb{R}^{n+2}$ of $f$.

Finally,  by Proposition \ref{p:equivalentcartanhyp}, Lemma \ref{l:reductionlemma}
and Proposition \ref{p:paramet}, there are one-to-one correspondences between 
(congruence classes of) genuine conformal deformations of $f$ in $\mathbb{R}^{n+2}$, 
triples $(D_1,D_2,\psi)$ on $M^n$ as in 
Proposition \ref{p:equivalentcartanhyp},  triples 
$(\bar{D}_1, \bar{D}_2,\bar{\psi})$ on $L^2$ as 
in Lemma \ref{l:reductionlemma}, and elements 
of $\mathcal{C}_s$. In summary, genuine conformal deformations of $f$ 
in $\mathbb{R}^{n+2}$ are parametrized by $\mathcal{C}_s$.\qed

{\renewcommand{\baselinestretch}{1}
\hspace*{-30ex}\begin{tabbing}
\indent \= IMPA  \hspace{30ex} Universidade Federal de S\~ao Carlos \\
\>  Estrada Dona Castorina, 110 \hspace{6.5ex}
Via Washington Luiz km 235 \\
\> 22460-320 --- Rio de Janeiro
\hspace{7.5ex} 13565-905 --- S\~ao Carlos  \\
\> Brazil\hspace{31ex} Brazil\\
\> sergio.chion@impa.br  \hspace{15ex}
tojeiro@dm.ufscar.br
\end{tabbing}}

\end{document}